\newcommand{\prof}{\noindent\(\clubsuit\)\ }
\newcommand{\qed}{\hspace*{3pt}\hfill\(\clubsuit\)}

\newcommand{\bN}{\mathbb{N}}

\newcommand{\bR}{\mathbb{R}}
\newcommand{\bZ}{\mathbb{Z}}

\newcommand{\fh}{\mathfrak{h}}
\newcommand{\Bar}{\textrm{Bar}}
\newcommand{\Cobar}{\textrm{Cobar}}

\newcommand{\rmid}{\textrm{id}}
\newcommand{\rmim}{\textrm{im}}
\newcommand{\rrdc}{\Rightarrow\hspace{-0.7em}\Rightarrow}
\newcommand{\sND}{^{\textrm{\scriptsize ND}}}
\newcommand{\sD}{^{\textrm{\scriptsize D}}}
\newcommand{\st}{\ \textbf{\underline{st}}\ }
\newcommand{\uD}{\underline{\Delta}}
\newcommand{\up}{\underline{p}}
\newcommand{\uq}{\underline{q}}
\newcommand{\wC}{\widehat{C}}
\newcommand{\hd}{\widehat{d}}
\newcommand{\wdl}{\widehat{\delta}}

\newenvironment{ctikzpicture}{%
\begin{equation}\begin{tikzpicture}}
{\end{tikzpicture}\end{equation}}

\newenvironment{fenum}{%
\vspace{-2\parskip}\begin{enumerate} \itemsep = -0.2em%
}{%
\end{enumerate}
}
\newlength{\lengtha}

\newtheorem{thr}{Theorem}

\newtheorem{dfn}[thr]{Definition}
\newtheorem{prp}[thr]{Proposition}

\documentclass[a4paper,12pt]{article}
\usepackage{amsfonts}
\usepackage{amssymb}
\usepackage{verbatim}
\usepackage{graphicx}
\usepackage[usenames,dvipsnames]{color}
\usepackage{tikz}
\usetikzlibrary{calc,matrix,fit}
\usepackage{hyperref}

\definecolor{green}{rgb}{0,0.75,0}
\definecolor{violet}{rgb}{0.8225,0.2855,0.5225}
\definecolor{pink}{rgb}{1,0.5,0.5}
\definecolor{orange}{rgb}{0.8,0.6,0.3}

\title{Discrete Vector Fields\\and Fundamental Algebraic Topology.}
\author{{\small\emph{Ana Romero, Francis
Sergeraert\footnote{Both authors partially supported by Ministerio de Ciencia e
Innovación, Spain, Project MTM2009-13842-C02-01.}} }}
\date{[Version 5.6, May 28, 2010.]\footnotesize}

\begin{document}

\voffset=-2.5cm \hoffset=-0.9cm \sloppy


\maketitle

\tableofcontents

\section{Introduction.}

We show in this text how the most important homology equivalences of
fundamental Algebraic Topology can be obtained as \emph{reductions} associated
to \emph{discrete vector fields}. Mainly the homology equivalences whose
\emph{existence} ---~most often non-constructive~-- is proved by the main
spectral sequences, the Serre and Eilenberg-Moore spectral sequences. On the
contrary, the \emph{constructive} existence is here systematically looked for
and obtained.

Algebraic topology consists in applying \emph{algebraic} methods to study
\emph{topological} objects. Algebra is assumed to be more tractable than
Topology, and the motivation of the method is clear.

Algebraic Topology sometimes reduces non-trivial topological problems to some
algebraic problems which, in favourable cases, can be solved. For example, the
Brouwer theorem is reduced to the impossibility of factorizing the identity
\(\bZ \rightarrow \bZ\) into a composition \(\bZ \rightarrow 0 \rightarrow
\bZ\). Magic!

Algebraic topology leads to more and more sophisticated algebraic translations,
think for example of the spectral sequences, derived categories,
\(E_\infty\)-operads, \ldots\ Almost all the topologists are somewhat
transformed into algebraists, specialists in \emph{Homological Algebra}. This
is so true that sometimes it happens some elementary topological methods are
neglected. These ``elementary'' methods can however be very powerful, in
particular when reexamining Algebraic Topology with a \emph{constructive} view.

The \emph{constructive} point of view, so new in this area, forces the
algebraic topologist to carefully reexamine the very bases of his methods.
Because of the deep connection ``constructive'' = ``something which can be
processed on a computer'' and because a computer is unable to solve
sophisticated problems by itself, a constructive version of some theory must be
split into elementary steps, elementary enough to be in the scope of a computer
programming language.

For example in the paper~\cite{rbsr06}, the classical Serre and Eilenberg-Moore
spectral sequences have been replaced by a more elementary tool, the
\emph{Homological Perturbation Theorem}\footnote{Usually called the \emph{Basic
Perturbation Lemma}, a strange terminology when we observe its importance now,
fifty years after its discovery by Shih Weishu~\cite{shih}.}. Elementary enough
to be easily installed on a computer and it was so possible to obtain homology
and homotopy groups so far unreachable. Working with this elementary result
allowed us to replace the \mbox{--~non}-constructive~-- Serre and
Eilenberg-Moore spectral sequences by a \emph{constructive} process, more
elementary but constructive, and finally more powerful.

This paper goes along the same line. The most elementary tool in homotopy, the
Whitehead contraction, is systematically studied to obtain the main results of
fundamental algebraic topology.

The basic tool is a direct adaptation of the so-called \emph{discrete Morse
theory}. More precisely the notion of discrete vector field is described here
in an algebraic setting; it is nothing but a rewriting of the main part of
Robin Forman's wonderful paper~\cite{frmn}. A rewriting taking account of the
long experience learned when designing our methods of \emph{effective
homology}~\cite{srgr3, drss, rbsr06, rbsr07, rbsr08, rbsr09}. The main result
of this part is Theorem~\ref{81450}, the Vector-Field Reduction theorem; it is
implicitly contained in~\cite{frmn}, but we hope the presentation given here
through the essential notion of (homological) reduction, see
Section~\ref{86769}, should interest the reader. We use again here the
Homological Perturbation Theorem to obtain a very direct proof of this result,
leading in an elegant way to interesting satellite results, in particular when
the naturality of the process must be studied.

Then the various homology equivalences which are the very bases of Algebraic
Topology are restudied and proved to be in fact direct consequences of this
elementary theorem. Mainly, the normalization theorems, the Eilenberg-Zilber
theorem --~the ordinary one and the \emph{twisted} one as well~--, the Bar and
Cobar reductions which are essential in the Eilenberg-Moore spectral sequences
to compute the \emph{effective} homology of classifying spaces and loop spaces.
Obvious applications to the homological analysis of digital images are also
included.

This gives a very clear and simple understanding of all these results under a
form of combinatorial game playing with the \emph{degeneracy operators}. Since
the remarkable works by Sam Eilenberg and Saunders MacLane in the fifties, by
Daniel Kan in the sixties, these operators have somewhat been neglected. We
hope the results obtained here show the work around the combinatorial nature of
these operators is far from being finished.

This is of course interesting for our favourite theory, but other applications
are expected: this simple way to understand our main homology equivalences
gives also new methods of programming: the heart of the method is extremely
simple and allows the programmer to carefully concentrate his work on the very
kernel of the method. Compare for example the method which was used up to now
to program the Eilenberg-Moore spectral sequences, for
example~\cite[Section~4]{brrs}, with the \emph{direct reductions} which are now
very simply obtained in Sections~\ref{00896} and~\ref{40626}. The programming
experiences already undertaken with respect to the Eilenberg-Zilber reduction,
in particular in the twisted case, are very encouraging.

\section{Discrete vector fields.}

\subsection{W-contractions.}

\begin{dfn} --- \label{41027}
\emph{An \emph{elementary W-contraction} is a pair \((X,A)\) of simplicial
sets, satisfying the following conditions:
\begin{fenum}
\item
The component \(A\) is a \emph{simplicial subset} of the simplicial set \(X\).
\item
The difference \(X-A\) is made of exactly \emph{two} non-degenerate simplices
\(\tau \in X_{n}\) and \(\sigma \in X_{n-1}\), the second one \(\sigma\) being
a \emph{face} of the first one \(\tau\).
\item \label{09884}
The incidence relation \(\sigma = \partial_i \tau\) holds for a \emph{unique}
index \(i \in 0...n\).
\end{fenum}
\hspace*{\parindent}It is then said \(A\) is obtained from \(X\) by an
elementary W-contraction, and \(X\) is obtained from \(A\) by an elementary
W-extension.\qed}
\end{dfn}

For example, \(X\) could be made of three triangles and \(A\) of two only as in
the next figure.
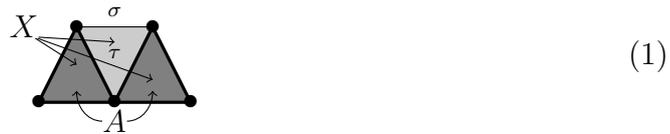
\begin{equation} \label{13789}
\begin{tikzpicture} [scale = 0.5, baseline = 0.5cm]
 \coordinate (0) at (0,0) ;
 \coordinate (1) at (2,0) ;
 \coordinate (2) at (4,0) ;
 \coordinate (3) at (1,2) ;
 \coordinate (4) at (3,2) ;
 \draw [fill = black!20] (1) -- (3) -- node [above] {\scriptsize\(\sigma\)} (4) ;
 \draw [very thick, fill = black!50] (0) -- (1) --(2) -- (4) -- (1) -- (3) -- (0) ;
 \foreach \i in {0,...,4}
   {\node at (\i) {\(\bullet\)} ;}
 \node at (2,1.3) {\scriptsize\(\tau\)} ;
 \node [inner sep = 0pt] (A) at (2,-0.5) {\(A\)} ;
 \draw [->] (A.0) to [out = 0, in = -90] (3,0.3) ;
 \draw [->] (A.180) to [out = 180, in = -90] (1,0.3) ;
 \node [anchor = 0, inner sep = 0pt] (X) at (0,2) {\(X\)} ;
 \draw [->] (X.-45) -- (1,1) ;
 \draw [->] (X.-40) -- (3,0.6) ;
 \draw [->] (X.-35) -- (2,1.6) ;
\end{tikzpicture}
\end{equation}

The condition~\ref{09884} is necessary -- and sufficient -- for the existence
of a topological contraction of \(X\) on \(A\). Think for example of the
minimal triangulation of the real projective plane \(P^2\bR\) as a simplicial
set \(X\), see the next figure: one vertex \(\ast\), one edge \(\sigma\) and
one triangle \(\tau\); no choice for the faces of \(\sigma\); the faces of
\(\tau\) must be \(\partial_0 \tau = \partial_2 \tau = \sigma\) and
\(\partial_1 \tau = \eta_0 \ast\) is the degeneracy of the base point. The
realization of \(X\) is homeomorphic to \(P^2\bR\). If you omit the
condition~\ref{09884} in the definition of W-contraction, then \((X, \ast)\)
would be a W-contraction, but \(P^2\bR\) is not contractible.
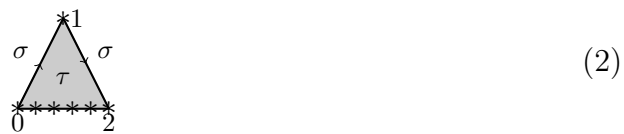
\begin{equation} \label{58160}
\begin{tikzpicture} [scale = 0.6, baseline = 0.5cm]
 \coordinate (0) at (0,0) ;
 \coordinate (1) at (2,0) ;
 \coordinate (2) at (1,2) ;
 \draw [thick, fill = black!20] (0) -- (1) -- node [anchor = -150]
   {\footnotesize\(\sigma\)} (2) -- node [anchor = -30]
   {\footnotesize\(\sigma\)} (0) ;
 \node at (1,0.67) {\footnotesize\(\tau\)} ;
 \draw [->] (0) -- (0.5,1) ;
 \draw [->] (2) -- (1.5,1) ;
 \foreach \i in {0,1,2} {\node at (\i) {\(\ast\)} ;}
 \foreach \i in {0.4, 0.8, 1.2, 1.6} {\node at (\i,0) {\(\ast\)} ;}
 \begin{scope} [font = \footnotesize]
 \node at (0,-0.3) {0} ;
 \node at (2,-0.3) {2} ;
 \node at (1.3,2) {1} ;
 \end{scope}
\end{tikzpicture}\end{equation}

\begin{dfn} --- \label{11769}
 \emph{A \emph{W-contraction} is a pair \((X,A)\) of simplicial sets satisfying
 the following conditions:
 \begin{fenum}
 \item
 The component \(A\) is a \emph{simplicial subset} of the simplicial set \(X\).
 \item There exists a sequence \((A_i)_{0 \leq i \leq m}\) with:
 \begin{fenum}
 \item
 \(A_0 = A\) and \(A_m = X\).
 \item
 For every \(0 < i \leq m\), the pair \((A_i, A_{i-1})\) is an elementary
 W-contraction.\qed
 \end{fenum}
 \end{fenum}}
\end{dfn}

In other words, a W-contraction is a finite sequence of elementary
contractions. If \((X,A)\) is a W-contraction, then a topological contraction
\(X \rightarrow A\) can be defined.

`W' stands for J.H.C. Whitehead, who undertook~\cite{whthj} a systematic study
of the notion of \emph{simple homotopy type}, defining two simplicial objects
\(X\) and \(Y\) as having the same simple homotopy type if they are equivalent
modulo the equivalence relation generated by the elementary W-contractions and
W-extensions.

Another kind of modification when examining a topological object can be
studied. Let us consider the usual triangulation of the square with two
triangles, the square cut by a diagonal. Then it is tempting to modify this
triangulation by pushing the diagonal onto two sides as roughly described in
this figure.
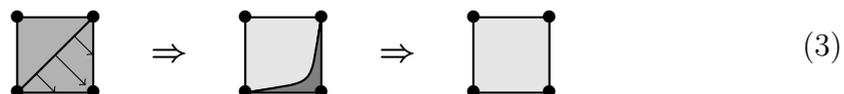
\begin{equation}\begin{tikzpicture} [baseline  = 0.5cm]
 \coordinate (0) at (0,0) ; \coordinate (1) at (1,0) ;
 \coordinate (2) at (1,1) ; \coordinate (3) at (0,1) ;
 \draw [thick, fill = black!30] (0) \foreach \i in {1,2,3,0,2} {-- (\i)} ;
 \draw [->] (0.75,0.75) -- (1,0.5) ;
 \draw [->] (0.5,0.5) -- (0.9,0.1) ;
 \draw [->] (0.25,0.25) -- (0.5,0) ;
 \foreach \i in {0,...,3} {\node at (\i) {\(\bullet\)} ;}
 \begin{scope} [xshift = 3cm]
 \coordinate (0) at (0,0) ; \coordinate (1) at (1,0) ;
 \coordinate (2) at (1,1) ; \coordinate (3) at (0,1) ;
 \draw [thick, fill = black!10] (0) .. controls +(10:1.2) and +(-100:1.2) .. (2)
 -- (3) -- (0) ;
 \draw [thick, fill = black!50] (0) .. controls +(10:1) and +(-100:1) .. (2)
 -- (1) -- (0) ;
 \foreach \i in {0,...,3} {\node at (\i) {\(\bullet\)} ;}
 \end{scope}
 \begin{scope} [xshift = 6cm]
 \coordinate (0) at (0,0) ; \coordinate (1) at (1,0) ;
 \coordinate (2) at (1,1) ; \coordinate (3) at (0,1) ;
 \draw [thick, fill = black!10] (0) \foreach \i in {1,2,3,0} {-- (\i)} ;
 \foreach \i in {0,...,3} {\node at (\i) {\(\bullet\)} ;}
 \end{scope}
 \node at (2,0.5) {{\boldmath\(\Rightarrow\)}} ;
 \node at (5,0.5) {{\boldmath\(\Rightarrow\)}} ;
\end{tikzpicture}
\end{equation}
Why not, but this needs other kinds of cells, here a square with \emph{four}
edges, while in a simplicial framework, the only objects of dimension 2 that
are provided are the triangles~\(\Delta^2\). Trying to overcome this essential
obstacle leads to two major subjects:
\begin{fenum}
\item
The Eilenberg-Zilber theorem, an algebraic translation of this idea, which
consists in \emph{algebraically} allowing the use of simplex products.
\item
The CW-complex theory, where the added cells are attached to the previously
constructed object through \emph{arbitrary} attaching maps.
\end{fenum}

This paper systematically reconsiders these essential ideas through the notion
of \emph{discrete vector field}.

\subsection{Algebraic discrete vector fields.}

The notion of discrete vector field (DVF) is due to Robin Forman~\cite{frmn};
it is an essential component of the so-called \emph{discrete Morse theory}. It
happens the notion of DVF will be here the major tool to treat the
\emph{fundamental} problems of algebraic topology, more precisely to treat the
general problem of \emph{constructive} algebraic topology.

This notion is usually described and used in combinatorial \emph{topology}, but
a purely algebraic version can also be given; we prefer this context.

\begin{dfn} --- \label{11799}
\emph{An \emph{algebraic cellular complex} (ACC) is a family:
\[
C = (C_p, d_p, \beta_p)_{p \in \bZ}
\]
of \emph{free} \(\bZ\)-modules and boundary maps.  Every \(C_p\) is called a
\emph{chain group} and is provided with a \emph{distinguished} \(\bZ\)-basis
\(\beta_p\); every basis component \(\sigma \in \beta_p\) is a
\emph{\(p\)-cell}. The boundary map \(d_p: C_p \rightarrow C_{p-1}\) is a
\(\bZ\)-linear map connecting two consecutive chain groups. The usual boundary
condition \(d_{p-1} d_p = 0\) is satisfied for every \(p \in \bZ\).}\qed
\end{dfn}

Most often we omit the index of the differential, so that the last condition
can be denoted by \(d^2 = 0\). The notation is redundant: necessarily, \(C_p =
\bZ[\beta_p]\), but the standard notation \(C_p\) for the group of \(p\)-chains
is convenient.

The chain complex associated to any sort of topological cellular complex is an
ACC. We are specially interested in the chain complexes associated to
simplicial sets.

Important: we \emph{do not} assume \emph{finite} the distinguished bases
\(\beta_p\), the chain groups are not necessarily of \emph{finite type}. This
is not an artificial extension to the traditional Morse theory: this point will
be often essential, but this extension is obvious.

\begin{dfn} ---
\emph{Let \(C\) be an ACC. A \((p-1)\)-cell \(\sigma\) is said to be a
\emph{face} of a \(p\)-cell~\(\tau\) if the coefficient of \(\sigma\) in \(d
\tau\) is non-null. It is a \emph{regular face} if this coefficient is +1 or
-1.}\qed
\end{dfn}

If \(\Delta^p\) is the standard simplex, every face of every subsimplex is a
regular face of this subsimplex. We gave after Definition~\ref{41027} an
example of triangulation of the real projective plane as a simplicial set; the
unique non-degenerate 1-simplex \(\sigma\) is \emph{not} a regular face of the
triangle \(\tau\), for \(d \tau = 2 \sigma\).

Note also the \emph{regular} property is \emph{relative}: \(\sigma\) can be a
regular face of \(\tau\) but also a non-regular face of another simplex
\(\tau'\).

\begin{dfn} --- \label{09885}
\emph{A \emph{discrete vector field} \(V\) on an algebraic cellular complex \(C
= (C_p, d_p, \beta_p)_{p \in \bZ}\) is a collection of pairs \(V = \{(\sigma_i,
\tau_i)\}_{i \in \beta}\) satisfying the conditions:
\begin{fenum}
\item
Every \(\sigma_i\) is some \(p\)-cell, in which case the other corresponding
component \(\tau_i\) is a \((p+1)\)-cell.  The degree \(p\) depends on \(i\)
and in general is not constant.
\item
Every component \(\sigma_i\) is a \emph{regular} face of the corresponding
component \(\tau_i\).
\item \label{18618}
A cell of \(C\) appears \emph{at most one time} in the vector field:  if \(i
\in \beta\) is fixed, then \(\sigma_i \neq \sigma_j\), \(\sigma_i \neq
\tau_j\), \(\tau_i \neq \sigma_j\) and \(\tau_i \neq \tau_j\) for every \(i
\neq j \in \beta\).\qed
\end{fenum}}
\end{dfn}

It is not required all the cells of \(C\) appear in the vector field \(V\). In
particular the void vector field is allowed. In a sense the remaining cells are
the most important.

\begin{dfn} ---
\emph{A cell \(\chi\) which does not appear in a discrete vector field \(V =
\{(\sigma_i, \tau_i)\}_{i \in \beta}\) is called a \emph{critical} cell. A
component \((\sigma_i, \tau_i)\) of the vector field \(V\) is a
\emph{\(p\)-vector} if \(\sigma_i\) is a \(p\)-cell. }\qed
\end{dfn}

 We do not consider in this
paper the traditional vector fields of differential geometry, which allows us
to call simply a \emph{vector field} which should be called a \emph{discrete
vector field}.

In case of an ACC coming from a topological cellular complex, a vector field is
a recipe to cancel ``useless'' cells in the underlying space, useless with
respect to the homotopy type. A component \((\sigma_i, \tau_i)\) of a vector
field can vaguely be thought of as a ``vector'' starting from the center of
\(\sigma_i\), going to the center of \(\tau_i\). For example \(\partial
\Delta^2\) and the circle have the same homotopy type, which is described by
the following scheme:\vspace{-30pt}
\begin{equation}\begin{tikzpicture} [scale = 0.5, baseline = (baseline)]
 \coordinate (baseline) at (0,0.5) ;
 \coordinate (0) at (0,0) ;
 \coordinate (1) at (1,2) ;
 \coordinate (2) at (2,0) ;
 \begin{scope} [font = \footnotesize]
 \node [anchor = 0] at (0) {0} ;
 \node [anchor = 0] at (1) {1} ;
 \node [anchor = 180] at (2) {2} ;
 \node [anchor = 0] at (6,0) {0} ;
 \node [anchor = 0] at (0.5,1) {01} ;
 \node [anchor = 180] at (1.5,1) {12} ;
 \node [anchor = 90] at (1,0) {02} ;
 \draw (6,0) .. controls +(60:4) and +(0:4) .. node [right] {12} (6,0) ;
 \end{scope}
 \foreach \i in {0,1,2} {\node at (\i) {\(\bullet\)} ;}
 \draw (0) \foreach \i in {1,2,0} {-- (\i)} ;
 \draw [thick, ->] (1) -- (0.5,1) ;
 \draw [thick, ->] (2) -- (1,0) ;
 \node at (4,1) {\(\Rightarrow\)} ;
 \node at (6,0) {\(\bullet\)} ;
\end{tikzpicture}
\end{equation}
The initial simplicial complex is made of three 0-cells 0, 1 and 2, and three
1-cells 01, 02 and 12. The drawn vector field is \(V = \{(1,01), (2,02)\}\),
and this vector field defines a homotopy equivalence between
\(\partial\Delta^2\) and the minimal triangulation of the circle as a
simplicial set. The last triangulation is made of the critical cells 0 and 12,
attached according to a process which deserves to be seriously studied in the
general case. This paper is devoted to a systematic use of this idea.

\subsection{V-paths and admissible vector fields.}

By the way, what about this vector field in \(\partial \Delta^2\)? \hfill
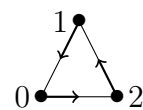
\begin{tikzpicture} [scale = 0.5, baseline = 0cm]
 \draw (0,0) node [anchor = 0] {0} node {\(\bullet\)}
      -- (2,0) node [anchor = 180] {2} node {\(\bullet\)}
      -- (1,2) node [anchor = 0] {1} node {\(\bullet\)} -- (0,0) ;
 \draw [thick, ->] (0,0) -- (1,0) ;
 \draw [thick, ->] (2,0) -- (1.5,1) ;
 \draw [thick, ->] (1,2) -- (0.5,1) ;
\end{tikzpicture}\hfill\hspace{0pt}

No critical cell and yet \(\partial \Delta^2\) does not have the homotopy type
of the void object. We must forbid possible \emph{loops}. This is not enough.
Do not forget the infinite case must be also covered; but look at this picture:
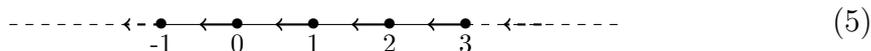
\begin{equation}\begin{tikzpicture} [font = \footnotesize, baseline = -0.5ex]
 \draw [dashed, thick, ->] (-1,0) -- (-1.5,0) ;
 \draw [dashed, thick, ->] (4,0) -- (3.5,0) ;
\foreach \i in {-1,0,...,3}
  {\node at (\i,0) {\(\bullet\)} ;
   \node [anchor = 90] at (\i,0) {\i} ;}
\foreach \i in {0,1,2,3}
  {\draw (\i,0) -- +(-1,0) ;
   \draw [thick, ->] (\i,0) -- +(-0.5,0) ;}
 \draw [dashed] (-3,0) -- (-1,0) ; \draw [dashed] (5,0) -- (3,0) ;
\end{tikzpicture}
\end{equation}
representing an infinite vector field on the real line triangulated as an
infinite union of 1-cells connecting successive integers. No critical cell and
yet the real line does not have the homotopy type of the void set. We must also
forbid the possible \emph{infinite} paths.

The notions of V-paths and admissible vector fields are the appropriate tools
to define the necessary restrictions.

\begin{dfn} --- \label{76422}
\emph{If \(V = \{(\sigma_i, \tau_i)\}_{i \in \beta}\) is a vector field on an
algebraic cellular complex \(C = (C_p, d_p, \beta_p)_p\), a \(V\)-path of
degree \(p\) is a sequence \(\pi = ((\sigma_{i_k}, \tau_{i_k}))_{0 \leq k <
m}\) satisfying:
\begin{fenum}
\item
Every pair \(((\sigma_{i_k}, \tau_{i_k}))\) is a component of the vector field
\(V\) and the cell \(\tau_{i_k}\) is a \(p\)-cell.
\item
For every \(0 < k < m\), the component \(\sigma_{i_k}\) is a face of
\(\tau_{i_{k-1}}\), non necessarily regular, but different from
\(\sigma_{i_{k-1}}\).
\end{fenum}
If \(\pi = ((\sigma_{i_k}, \tau_{i_k}))_{0 \leq k < m}\) is a \(V\)-path, and
if \(\sigma\) is a face of \(\tau_{i_{m-1}}\) different from
\(\sigma_{i_{m-1}}\), then \(\pi\) \emph{connects} \(\sigma_{i_0}\) and
\(\sigma\) through the vector field \(V\).} \qed\end{dfn}
\begin{equation}\begin{tikzpicture}
 [scale=0.5, font = \scriptsize, inner sep = 1pt, baseline = (baseline)]
 \coordinate (baseline) at (0,1) ;
 \foreach \i in {1,3,5} {\coordinate (0\i) at (\i,0) ;}
 \foreach \i in {0,2,4,6} {\coordinate (1\i) at (\i,2) ;}
 \fill [black!20] (01) -- (05) -- (16) -- (10) ;
 \foreach \i in {1,3,5} {\node [label = -90:\i] at (0\i) {\(\bullet\)} ;}
 \foreach \i in {0,2,4,6} {\node [label = 90:\i] at (1\i) {\(\bullet\)} ;}
 \draw [thick] (01) -- (05) (10) -- (16) ;
 \foreach \i in {1,3,5} {\draw [thick] (0\i) +(-1,2) -- +(0,0) -- +(1,2) ;}
 \foreach \i in {1,3,5} {\draw [->, thick] (0\i) ++(-0.5,1) -- +(0.5,0.25) ;}
 \foreach \i in {1,3} {\draw [->, thick] (0\i) ++(0.5,1) -- +(0.5,-0.25) ;}
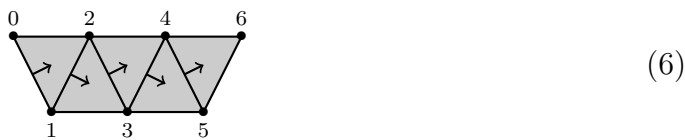
 \node [font = \footnotesize] at (3,-1.5) {A \(V\)-path connecting the edges 01 and 56.} ;
\end{tikzpicture}\end{equation}

In a \(V\)-path \(\pi = ((\sigma_{i_k}, \tau_{i_k}))_{0 \leq k < m}\) of degree
\(p\), a \((p-1)\)-cell \(\sigma_{i_k}\) is a regular face of \(\tau_{i_k}\),
for the pair \((\sigma_{i_{k}}, \tau_{i_k})\) is a component of the vector
field \(V\), but the same \(\sigma_{i_k}\) is non-necessarily a regular face of
\(\tau_{i_{k-1}}\).

\begin{dfn} ---
\emph{The \emph{length} of the path \(\pi = ((\sigma_{i_k}, \tau_{i_k}))_{0
\leq k < m}\) is \(m\).}\qed
\end{dfn}

If \((\sigma, \tau)\) is a component of a vector field, in general the cell
\(\tau\) has \emph{several} faces different from \(\sigma\), so that the
possible paths starting from a cell generate an oriented graph.


\begin{dfn} ---
\label{96757} \emph{A discrete vector field \(V\) on an algebraic cellular
complex \(C = (C_p, d_p, \beta_p)_{p \in \bZ}\) is \emph{admissible} if for
every \(p \in \bZ\), a function \(\lambda_p: \beta_p \rightarrow \bN\) is
provided satisfying the following property: every \(V\)-path starting from
\(\sigma \in \beta_p\) has a length bounded by \(\lambda_p(\sigma)\).}\qed
\end{dfn}

Excluding infinite paths is almost equivalent. The difference between both
possibilities is measured by Markov's principle; we prefer our more
constructive statement.

A circular path would generate an infinite path and is therefore excluded.

The next diagram, an \emph{oriented bipartite graph}, can help to understand
this notion of admissibility for some vector field \(V\). This notion makes
sense degree by degree. Between the degrees \(p\) and \(p-1\), organize the
\emph{source} \((p-1)\)-cells (resp. \emph{target} \(p\)-cells) as a lefthand
(resp. righthand) column of cells. Then every vector \((\sigma, \tau) \in V\)
produces an oriented edge \(\sigma \rightarrowtail \tau\). In the reverse
direction, if \(\tau\) is a target \(p\)-cell, the boundary \(d\tau\) is a
finite linear combination \(d\tau = \sum \alpha_i \sigma_i\), and some of these
\(\sigma_i\)'s are source cells, in particular certainly the corresponding
\(V\)-source cell~\(\sigma\). For every such source component \(\sigma_i\), be
careful \emph{except} for the corresponding source~\(\sigma\), you install an
oriented edge \(\sigma_i \leftarrowtail \tau\).

Then the vector field is admissible between the degrees \(p-1\) and \(p\) if
and only if, starting from some source cell \(\sigma\), all the (oriented)
paths have a length bounded by some integer \(\lambda_p(\sigma)\). In
particular, the loops are excluded. We draw the two simplest examples of vector
fields non-admissible. The lefthand one has an infinite path, the righthand one
has a loop, a particular case of infinite path.

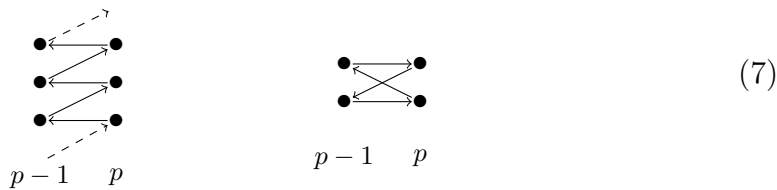
\begin{equation} \label{97760}
\begin{tikzpicture} [yscale=0.5, inner sep = 0pt, baseline = (baseline)]
 \coordinate (baseline) at (0,1) ;
 \foreach \i in {0,1} \foreach \j in {0,1,2}
   {\node (\i\j) at (\i,\j) {\(\bullet\)} ; }
 \foreach \j in {0,1,2}
   {\draw [->] (1\j) -- (0\j) ;}
 \foreach \j/\jj in {0/1,1/2}
   {\draw [->] (0\j) -- (1\jj) ;}
 \draw [->, dashed] (0.1,-1) -- (10) ;
 \draw [->, dashed] (02) -- (0.9,2.9) ;
 \node at (0,-1.5) {\footnotesize\(p-1\)} ;
 \node at (1,-1.5) {\footnotesize\(p\)} ;
 \begin{scope} [xshift = 4cm, yshift = 0.5cm]
 \foreach \i in {0,1} \foreach \j in {0,1}
   {\node (\i\j) at (\i,\j) {\(\bullet\)} ;}
 \foreach \i/\j/\ii/\jj in {0/0/1/0,1/0/0/1,0/1/1/1,1/1/0/0}
   {\draw [->] (\i\j) -- (\ii\jj) ;}
 \node at (0,-1.5) {\footnotesize\(p-1\)} ;
 \node at (1,-1.5) {\footnotesize\(p\)} ;
 \end{scope}

\end{tikzpicture}
\end{equation}


\begin{dfn} --- \label{92078}
\emph{Let \(V = \{(\sigma_i, \tau_i)_{i \in \beta}\}\) be a vector field on an
ACC. A \emph{Lyapunov function} for \(V\) is a function \(L: \beta \rightarrow
\bN\) satisfying the following condition: if \(\sigma_j\) is a face of
\(\tau_i\) different from \(\sigma_i\), then \(L(j) < L(i)\).}\qed
\end{dfn}

It is the natural translation in our discrete framework of the traditional
notion of Lyapunov function in differential geometry. It is clear such a
Lyapunov function proves the admissibility of the studied vector field. Obvious
generalizations to ordered sets more general than \(\bN\) are possible.

\subsection{Reductions.} \label{86769}

\subsubsection{Definition.}

\begin{dfn} --- \label{93254}
\emph{A (homology) \emph{reduction\footnote{Often called \emph{contraction},
but this terminology is not appropriate: it is important to understand such a
reduction has an \emph{algebraic} nature, like when you \emph{reduce} 6/4
\(\mapsto\) 3/2. When reducing a rational number, you cancel the opposite roles
of a common factor in numerator and denominator; in our homology reductions we
cancel the opposite roles of the \(A_\ast\) and \(B_\ast\) components in the
big chain complex to obtain the small one, see Figure~\ref{36743}.}} \(\rho\)
is a diagram:
\begin{equation}
 \begin{tikzpicture} [baseline = (baseline)]
 \begin{scope} [anchor = base, inner sep = 2pt]
 \coordinate (baseline) at (0,0) ;
 \node (1) at (0,0) {\(\rho =\)} ;
 \node (2) at (2.2,0) {\(\wC_\ast\)} ;
 \node (3) at (4.2,0) {\(C_\ast\)} ;
 \begin{scope} [font = \footnotesize, ->]
 \draw (2.base east) -- node [anchor = 90] {\(f\)} (3.base west) ;
 \draw ($(3.base west) + (0,5pt) $) -- node [anchor = -90] {\(g\)}
       ($(2.base east) + (0,5pt)$)  ;
 \draw ($(2.base west) + (0,7pt)$) .. controls +(160:1.2) and +(-160:1.2) ..
       node [left] {\(h\)} ($(2.base west) + (0,2pt)$) ;
 \draw (0.5,-0.55) rectangle (4.7,0.8) ;
 \end{scope}
 \end{scope}
 \end{tikzpicture}
\end{equation}
 where:
 \begin{fenum}
 \item
 The nodes \(\wC_\ast\) and \(C_\ast\) are chain complexes;
 \item
 The arrows \(f\) and \(g\) are chain complex morphisms;
 \item
 The self-arrow \(h\) is a homotopy operator, of degree +1;
 \item
 The following relations are satisfied:
 \begin{equation}
 \begin{array}{rcl}
 fg &=& \rmid_{C_\ast}
 \\
 gf + dh + hd &=& \rmid_{\wC_\ast}
 \\
 fh &=& 0
 \\
 hg &=& 0
 \\
 hh &=& 0
 \end{array}
 \end{equation}
 \end{fenum}}
\end{dfn}\vspace{-15pt}\qed

The relation \(fg = \rmid\) implies that \(g\) identifies the \emph{small}
chain complex \(C_\ast\) with the subcomplex \(C'_\ast := g(C_\ast)\) of the
\emph{big} chain complex \(\wC_\ast\). Furthermore the last one gets a
canonical decomposition \(\wC_\ast = \rmim(g) \oplus \ker(f)\). The relations
\(hg = 0\) and \(fh = 0\) imply the homotopy operator \(h\) is null on
\(\rmim(g) = C'_\ast\) and its image is entirely in \(\ker(f)\): the \(h\) map
is in fact defined on \(\ker(f)\), extended by the zero map on \(C'_\ast\).
Finally \(dh + hd\) is the identity map on \(\ker(f)\). Also, because of the
relation \(h^2 = 0\), the homotopy \(h\) is a codifferential and the pair
\((d,h)\) defines a Hodge decomposition \(\ker(f) = A_\ast \oplus B_\ast\) with
\(A_\ast = \rmim(h) = \ker(f) \cap \ker(h)\) and \(B_\ast = \ker(f) \cap
\ker(d) = \ker(f) \cap \rmim(d)\). The direct sum \(\ker(f) = A_\ast \oplus
B_\ast\) is a subcomplex of \(\wC_\ast\), but both components \(A_\ast\) and
\(B_\ast\) are only graded modules. These properties are illustrated in this
diagram.
\begin{equation} \label{36743}
  \begin{tikzpicture} [xscale = 2.2, yscale = 1.7, baseline = (n22)]
  \foreach \i/\ic in {0,1,2,3,4/4.25}
    {\foreach \j/\jc in {0,1,2,3,4/3.7} {\coordinate (\i\j) at (\ic,\jc) ;}}
  \foreach \j in {0,1,2,3,4} {\node (n0\j) at (0\j)
     {\(\{\cdots\)} ;}
  \foreach \i/\it in {1/p-1, 2/p, 3/p+1}
    {\foreach \j/\jt in {0/C, 1/C', 2/B, 3/A, 4/\wC}
      {\node (n\i\j) at (\i\j) {\(\jt_{\it}\)} ;}}
  \foreach \j/\jt in {0/C, 1/C', 2/B, 3/A, 4/\wC}
    {\node (n4\j) at (4\j) {\(\cdots\} =
    \hspace{5pt} \jt_\ast\)} ;}
  \begin{scope} [font = \scriptsize]
    \path (n14.base) coordinate (b) +(0,0.05) coordinate (bb) ;
    \foreach \i/\ii in {0/1, 1/2, 2/3, 3/4}
      {\draw [->] (b -| n\i4.east) -- node [below] {\(h\)} (b -| n\ii4.west) ;
       \draw [<-] (bb -| n\i4.east) -- node [above] {\(d\)} (bb -| n\ii4.west)
       ;}
    \foreach \i/\ii in {02/13, 12/23, 22/33, 32/43}
      {
      \path (n\i.north east) +(135:0.03) coordinate (1l) +(-45:0.03) coordinate (1r) ;
      \path (n\ii.south west) +(135:0.03) coordinate (2l) +(-45:0.03) coordinate (2r) ;
      \draw [->] (2l) -- (1l) ; \draw [->] (1r) -- (2r) ;
      \path (n\i.north east) --
      node [fill = white, inner sep = 1, label=135:\(d\),
      label = -45:\(h\)] {\(\cong\)} (n\ii.south west) ;
      }
    \foreach \i in {1,2,3}
      {
       \draw (n\i1.south west) rectangle (n\i3.north east) ;
       \path (\i1) -- node {\large\(\oplus\)} (\i2) -- node {\large\(\oplus\)} (\i3) ;
      }
    \foreach \i/\j in {1/2, 2/3}
      {
       \path (4\i) -- coordinate (c) (4\j) (c) +(0:0.32) node {\large\(\oplus\)};
      }
    \path (n43.north east) +(180:0.35) coordinate (c) ;
    \draw (n43.north east) rectangle (c |- n41.south) ;
    \foreach \i/\ii in {1/0,2/1,3/2,4/3} \foreach \j in {0,1}
      {\draw [->] (n\i\j) -- node [above] {\(d\)} (n\ii\j) ; }
    \foreach \i in {1,2,3}
      {
       \path (n\i0.north) +(180:0.03) coordinate (c1) +(0:0.03) coordinate (c2) ;
       \path (n\i1.south) +(180:0.03) coordinate (c3) +(0:0.03) coordinate (c4) ;
       \draw [->] (c2) -- (c4) ; \draw [->] (c3) -- (c1) ;
       \path (\i0) -- node [fill = white, inner sep = 1] {\(f \cong g\)} (\i1) ;
       \draw (c1 |- n\i3.north) +(90:0.1) -- (c1 |- n\i4.south)
             (c2 |- n\i3.north) +(90:0.1) -- (c2 |- n\i4.south) ;
      }
    \path (n40.north) ++(0:0.3) coordinate (c1) +(0:0.05) coordinate (c2) ;
    \draw [->] (c2) -- (c2 |- n41.south) ; \draw [->] (c1 |- n41.south) -- (c1) ;
    \path (c1) -- node [fill = white, inner sep = 1] {\(f \cong g\)} (c2 |- n41.south) ;
    \draw (c1 |- n43.north) +(90:0.1) -- (c1 |- n44.south)
          (c2 |- n43.north) +(90:0.1) -- (c2 |- n44.south) ;
\end{scope}
  \path (n00.south west) +(-135:0.2) coordinate (c1)
        (n44.north east) +(30:0.3) coordinate (c2) ;
  \draw [very thick] (c1) rectangle (c2) ;
  \path (c1) +(-90:0.4) coordinate (c1) ;
  \path (c1) -- node [draw, pos=0.167] {\(A_\ast = \ker f \cap \ker h\)}
                node [draw, pos = 0.5] {\(C'_\ast = \textrm{im}\,g\)}
                node [draw, pos= 0.833] {\(B_\ast = \ker f \cap \ker d\)} (c1 -|
                c2);
\end{tikzpicture}
\end{equation}

We will simply denote such a reduction by \(\rho = (f_\rho, g_\rho, h_\rho) :
\wC_\ast \rrdc C_\ast\) or simply by \(\rho: \wC_\ast \rrdc C_\ast\).

A reduction \(\rho = (f_\rho, g_\rho, h_\rho) : \wC_\ast \rrdc C_\ast\)
establishes a strong connection between the chain complexes \(\wC_\ast\) and
\(C_\ast\): it is a particular quasi-isomorphism between chain complexes
describing the big one \(\wC_\ast\) as the direct sum of the small one
\(C'_\ast = C_\ast\) and another chain-complex \(\ker(f)\), the latter being
provided with an \emph{explicit} null-reduction \(h\). The morphisms \(f\) and
\(g\) are inverse homology equivalences\footnote{Often called \emph{chain
equivalences}, yet they do not define an equivalence between \emph{chains} but
between \emph{homology} classes.}.

The main role of a reduction is the following. It often happens the big
complex~\(\wC_\ast\) is so enormous that the homology groups of this complex
are out of scope of computation; we will see striking examples where the big
complex is \emph{not of finite type}, so that the homology groups \(H_\ast
\wC_\ast\) are not computable from \(\wC_\ast\), even in theory. But if on the
contrary the small complex \(C_\ast\) has a reasonable size, then its homology
groups are computable and they are canonically isomorphic to those of the big
complex.

More precisely, if the \emph{homological problem} is solved in the small
complex, then the reduction produces a solution of the same problem for the big
complex. See~\cite[Section~4.4]{rbsr09} for the definition of the notion of
\emph{homological problem}. In particular, if \(z \in C_p\) is a cycle
representing the homology class \(\fh \in H_p(C_\ast)\), then \(gz\) is a cycle
representing the corresponding class in \(H_p(\wC_\ast)\). Conversely, if \(z\)
is a cycle of \(\wC_p\), then the homology class of this cycle corresponds to
the homology class of \(fz\) in \(H_p(C_\ast)\). If ever this homology class is
null and if \(c \in C_{p+1}\) is a boundary-preimage of \(fz\) in the small
complex, then \(gc + hz\) is a boundary-preimage of \(z\) in the big complex.

\begin{thr} ---
If \(\rho = (f,g,h): \wC_\ast \rrdc C_\ast\) is a reduction between the chain
complexes \(\wC_\ast\) and \(C_\ast\), then any homological problem in the big
complex \(\wC_\ast\) can be solved through a solution of the same problem in
the small complex \(C_\ast\).\qed
\end{thr}

\subsubsection{The Homological Perturbation Theorem.}

This theorem is often called the Basic Perturbation ``Lemma''. It was
introduced in a (very important) particular case by Shih Weishu~\cite{shih},
allowing him to obtain a \emph{more effective} version of the Serre spectral
sequence: Shih so obtained an \emph{explicit} homology equivalence between the
chain complex of the total space of a fibration and the corresponding Hirsch
complex. The general scope of this lemma was signalled by Ronnie
Brown~\cite{brwnr1}. Our organization of \emph{Effective Homology},
see~\cite{srgr3, drss, rbsr06, rbsr07, rbsr08, rbsr09}, consists in combining
this lemma with functional programming, more precisely with the notion of
locally effective object.

\begin{thr}\label{07404}{\bf (Homological Perturbation Theorem)} ---
Let \(\rho = (f,g,h): (\wC_\ast, \hd) \rrdc (C_\ast, d)\) be a reduction and
let \(\wdl\) be a perturbation of the differential~\(\hd\) of the big
chain-complex. We assume the \emph{nilpotency hypothesis} is satisfied: for
every \(c \in \wC_n\), there exists \(\nu \in \bN\) satisfying \((\wdl
h)^\nu(c) = 0\). Then a perturbation \(\delta\) can be defined for the
differential \(d\) and a new reduction \(\rho' = (f', g, h'): (\wC_\ast, \hd +
\wdl) \rrdc (C_\ast, d + \delta)\) can be constructed.
\end{thr}

The nilpotency hypothesis states the composition \(\wdl h\) is pointwise
nilpotent. The process described by the Theorem perturbs the differential \(d\)
of the small complex, becoming \(d + \delta\), and the components \((f,g,h)\)
of the reduction which so becomes a new reduction \(\rho' = (f',g',h'):
(\wC_\ast, \hd + \wdl) \rrdc (C_\ast, d + \delta)\).

Which is magic in the BPL is the fact that a sometimes complicated perturbation
of the ``big'' differential can be accordingly reproduced in the ``small''
differential; in general it is not possible, unless the nilpotency hypothesis
is satisfied.

\prof See~\cite[\S 1 and 2]{shih} and \cite{brwnr1}. A detailed proof in the
present context is at~\cite[Theorem 50]{rbsr09}.\qed

This theorem\footnote{Usually called a ``lemma''!} is so important in
\emph{effective} Homological Algebra that Julio Rubio's team at Logro\~no
decided to write a proof in the language of the Isabelle theorem prover, and
succeeded~\cite{arbr}; it is the starting point to obtain \emph{proved}
programs using this crucial result, a fascinating challenge.

\subsection{A vector field generates a reduction.}

Let \(C_\ast = (C_\ast, d_\ast, \beta_\ast)\) be an algebraic cellular complex
provided with an \emph{admissible} discrete vector field \(V\). Then a
reduction \(\rho: C_\ast \rrdc C^c_\ast\) can be constructed where the small
chain complex \(C^c_\ast\) is the \emph{critical} complex; it is also a
cellular complex but generated only by the \emph{critical} cells of
\(\beta_\ast\), those which do not appear in the vector field \(V\), with a
differential appropriately defined, combining the initial differential~\(d\) of
the initial complex and the vector field.

This result is due to Robin Forman~\cite[Section 8]{frmn}. It is here extended
to the complexes not necessarily of finite type, but the extension is obvious.
Many (slightly) different proofs are possible. Forman's proof uses an
intermediate chain complex, the \emph{Morse} chain complex, which depends in
fact only on the vector field, it is made of the chains that are invariant for
the flow canonically associated to the vector field. We give here two other
organizations, each one having its own interest.

The first one reduces the result to the standard Gauss elimination process for
the linear systems. It has the advantage of identifying the heart of the
process, the component \(d_{2,1}^{-1}\) of the formulas~(\ref{76423}). It is
nothing but the homotopy component~\(h\) of the looked-for reduction. For
critical time applications, storing the values of this component when it is
calculated could be a good strategy, the rest being a direct and quick
consequence. Our \(d_{2,1}^{-1}\) is the operator \(L\) of~\cite[Page
122]{frmn}.

The difference between \emph{ordinary} and \emph{effective} homology is easy to
see here. A ``simple'' user of Forman's paper could skip the construction of
this operator, used only to prove the Morse complex gives the right homology
groups. On the contrary, in effective homology, this operator is crucial, from
a theoretical point of view and for a computational point of view as well: the
profiler experiments show most computing time is \emph{then} devoted to this
operator, which is \emph{not} the case in ordinary homology.

The second proof uses the Homological Perturbation Theorem~\ref{07404}. It is
much faster and highlights the decidedly wide scope of this result. The general
style of application of this theorem is again there: starting from a particular
case where the result is obvious, the general case is viewed as a
\emph{perturbation} of the particular case. The nilpotency condition must be
satisfied, which amounts to requiring the vector field is admissible. The
critical complex of Forman is the bottom chain complex of the obtained
reduction, and the Morse complex is its canonical image in the top complex. How
to be simpler?

\subsubsection{Using Gauss elimination.}

\begin{prp} {\textbf{\emph{(Hexagonal lemma)}}} --- \label{90055}
Let \(C = (C_p, d_p)_p\) be a chain complex. For some \(k \in \bZ\), the chain
groups \(C_k\) and \(C_{k+1}\) are given with decompositions \(C_k = C'_k
\oplus C''_k\) and \(C_{k+1} = C'_{k+1} \oplus C''_{k+1}\), so that between the
degrees \(k-1\) and \(k+2\) this chain complex is described by the diagram:
\begin{equation}\begin{tikzpicture} [xscale = 2, baseline = (01.base)]
 \node (10) at (1,0) {\(C'_k\)} ;
 \node (20) at (2,0) {\(C'_{k+1}\)} ;
 \node (01) at (0,1) {\(C_{k-1}\)} ;
 \node (31) at (3,1) {\(C_{k+2}\)} ;
 \node (12) at (1,2) {\(C''_k\)} ;
 \node (22) at (2,2) {\(C''_{k+1}\)} ;
 \begin{scope} [->, font = \scriptsize]
 \draw (10) -- node [below] {\(\alpha\)} (01) ;
 \draw (20) -- node [below] {\(\beta\)} (10) ;
 \draw (31) -- node [below] {\(\gamma\)} (20) ;
 \draw (12) -- node [above] {\(\delta\)} (01) ;
 \draw (22) -- node [below] {\(\varepsilon\)} (12) ;
 \draw (31) -- node [above] {\(\eta\)} (22) ;
 \draw (20) -- node [above, pos = 0.7] {\(\varphi\)} (12) ;
 \draw (22) -- node [below, pos = 0.7] {\(\psi\)} (10) ;
 \draw ($(12.0)+(0,0.1)$) -- node [above] {\(\varepsilon^{-1}\)}
 ($(22.180)+(0,0.1)$);
 \end{scope}
 \draw [dotted] (10) -- node (11) [fill = white] {\(\oplus\)} (12) ;
 \draw [dotted] (20) -- node (21) [fill = white] {\(\oplus\)} (22) ;
 \draw [->, font = \scriptsize] (11) -- node [above, pos = 0.3]  {\(d\)} (01) ;
 \draw [->, font = \scriptsize] (31) -- node [above, pos = 0.7] {\(d\)} (21) ;
 \path (21) -- node [fill = white] {} (11) ;
 \draw [->, font = \scriptsize] (21) -- node [above, pos = 0.1] {\(d\)} (11) ;
\end{tikzpicture}
\end{equation}
The partial differential \(\varepsilon: C''_{k+1} \rightarrow C''_{k}\) is
assumed to be an isomorphism. Then a canonical reduction can be defined \(\rho:
C \rrdc C'\) where \(C'\) is the same chain complex as \(C\) except between the
degrees \(k-1\) and \(k+2\):
\begin{equation}\begin{tikzpicture} [xscale = 1.7, baseline = (0.base)]
 \node (-2) at (-2,0) {\(C_{k-2}\)} ;
 \node (-1) at (-1,0) {\(C_{k-1}\)} ;
 \node (0) at (0,0) {\(C'_k\)} ;
 \node (1) at (1.5,0) {\(C'_{k+1}\)} ;
 \node (2) at (2.5,0) {\(C_{k+2}\)} ;
 \node (3) at (3.5,0) {\(C_{k+3}\)} ;
\begin{scope} [->, font = \scriptsize, above]
 \draw [dashed] (-2) -- (-2.7,0) ;
 \draw (-1) -- (-2) ;
 \draw (0) -- node {\(\alpha\)} (-1) ;
 \draw (1) -- node {\(\beta - \psi \varepsilon^{-1} \phi\)} (0) ;
 \draw (2) -- node {\(\gamma\)} (1) ;
 \draw (3) -- (2) ;
 \draw [dashed] (4.2,0) -- (3) ;
\end{scope}
\end{tikzpicture}
\end{equation}
\end{prp}

\prof An integer matrix {\footnotesize \(\left[\begin{array}{cc}\varepsilon &
\phi \\ \psi & \beta\end{array}\right]\)} is equivalent to the matrix
{\footnotesize\(\left[\begin{array}{cc} \varepsilon & 0 \\ 0 & \beta - \psi
\varepsilon^{-1} \phi\end{array}\right]\)} if \(|\varepsilon| = 1\), it is the
simplest case of Gauss' elimination. More generally, the following matrix
relation is always satisfied, even if the matrix entries are in turn coherent
linear maps:
\begin{equation}
\left[\begin{array}{cc} \varepsilon & \varphi \\ \psi & \beta\end{array}\right]
= \left[\begin{array}{cc} 1 & 0 \\ \psi \varepsilon^{-1} & 1 \end{array}\right]
\left[\begin{array}{cc} \varepsilon & 0 \\ 0 & \beta - \psi \varepsilon^{-1}
\varphi \end{array}\right] \left[\begin{array}{cc} 1 & \varepsilon^{-1} \varphi
\\ 0 & 1\end{array}\right]
\end{equation}
The lateral matrices of the right-hand term can be considered as basis changes.
These matrices define an isomorphism \(\rho': C \rightarrow \overline{C}\)
between the initial chain complex C and the chain complex \(\overline{C}\) made
of the same chain groups but the differentials displayed on this diagram:
\begin{equation}\begin{tikzpicture} [xscale = 2, baseline = (01.base)]
 \node (10) at (1,0) {\(C'_k\)} ;
 \node (20) at (2.2,0) {\(C'_{k+1}\)} ;
 \node (01) at (0,1) {\(C_{k-1}\)} ;
 \node (31) at (3.2,1) {\(C_{k+2}\)} ;
 \node (12) at (1,2) {\(C''_k\)} ;
 \node (22) at (2.2,2) {\(C''_{k+1}\)} ;
 \begin{scope} [->, font = \scriptsize]
 \draw (10) -- node [below] {\(\alpha\)} (01) ;
 \draw (20) -- node [below] {\(\beta - \psi \varepsilon^{-1} \varphi\)} (10) ;
 \draw (31) -- node [below] {\(\gamma\)} (20) ;
 \draw (12) -- node [above] {\(0\)} (01) ;
 \draw (22) -- node [above] {\(\varepsilon\)} (12) ;
 \draw (31) -- node [above] {\(0\)} (22) ;
 \draw (20) -- node [above, pos = 0.7] {\(0\)} (12) ;
 \draw (22) -- node [below, pos = 0.7] {\(0\)} (10) ;
 \end{scope}
 \draw [dotted] (10) -- node (11) [fill = white] {\(\oplus\)} (12) ;
 \draw [dotted] (20) -- node (21) [fill = white] {\(\oplus\)} (22) ;
 \draw [->, font = \scriptsize] (11) -- node [above, pos = 0.3]  {\(d\)} (01) ;
 \draw [->, font = \scriptsize] (31) -- node [above, pos = 0.7] {\(d\)} (21) ;
 \path (21) -- node [fill = white] {} (11) ;
 \draw [->, font = \scriptsize] (21) -- node [above, pos = 0.1] {\(d\)} (11) ;
\end{tikzpicture}
\end{equation}

Throwing away the component \(\varepsilon: C''_{k+1} \rightarrow C''_{k}\) from
this chain complex \(\overline{C}\) produces a reduction \(\rho'': \overline{C}
\rrdc C'\) to the announced chain complex \(C'\). The desired reduction is
\(\rho = \rho'' \rho': C \rrdc C'\) with an obvious interpretation of the
composition \(\rho'' \rho'\). Finally \(\rho = (f,g,h)\) with:
\begin{fenum}
\item
The morphism \(f\) is the identity except:
\begin{equation}
f_{k} = \left[\begin{array}{cc} - \psi \varepsilon^{-1} & 1
\end{array}\right] \hspace{1cm} f_{k+1} = \left[\begin{array}{cc} 0 & 1 \end{array}
\right]
\end{equation}
\item
The morphism \(g\) is the identity except:
\begin{equation}
g_{k} = \left[\begin{array}{c}0 \\ 1\end{array}\right] \hspace{1cm} g_{k+1} =
\left[\begin{array}{c} -\varepsilon^{-1} \varphi \\ 1\end{array}\right]
\end{equation}
\item
The homotopy operator \(h\) is the null operator except:
\begin{equation}
h_k = \left[\begin{array}{cc} \varepsilon^{-1} & 0 \\ 0 & 0\end{array}\right]
\end{equation}
\end{fenum}
matrices to be interpreted via appropriate block decompositions.\qed

Note the boundary components \(\alpha\) and \(\gamma\) are not modified by the
reduction process. So that if independent ``hexagonal'' decompositions are
given for every degree, the process can be applied to every degree
simultaneously.

\begin{thr} --- \label{09803}
Let \(C = (C_p, d_p)_p\) be a chain complex. We assume every chain group is
decomposed \(C_p = D_p \oplus E_p \oplus F_p\). The boundary maps \(d_p\) are
then decomposed in \(3 \times 3\) block matrices \([d_{p,i,j}]_{1 \leq i,j \leq
3}\). If every component \(d_{p,2,1}: D_p \rightarrow E_{p-1}\) is an
isomorphism, then the chain complex can be canonically reduced to a chain
complex \((F_p, d'_p)\).
\end{thr}

\prof Simultaneously applying the formulas produced by the hexagonal lemma
gives the desired reduction, the components of which are:

\begin{equation} \label{76423}
 \begin{array}{cc}
 d'_p = d_{p,3,3} - d_{p,3,1} d_{p,2,1}^{-1} d_{p,2,3}
 &
 f_p = \left[\begin{array}{ccc} 0 & - d_{p,3,1} d_{p,2,1}^{-1} & 1 \end{array}\right]
 \\[10pt]
 g_p = \left[\begin{array}{c} - d_{p,2,1}^{-1} d_{p,2,3} \\ 0 \\ 1
 \end{array}\right]
 &
 h_{p-1} = \left[\begin{array}{ccc} 0 & d_{p,2,1}^{-1} & 0 \\ 0 & 0 & 0
 \\ 0 & 0 & 0 \end{array}\right]
 \end{array}
\end{equation}
\hspace*{0pt}\hfill\(\clubsuit\)

It is an amusing exercise to check the displayed formulas satisfy the required
relations: \(d'f = fd\), \(dg = gd'\), \(fg = 1\), \(dh + hd + fg = 1\), \(fh =
0\), \(hg = 0\), \(hh = 0\), stated in Definition~\ref{93254}. The components
\(d_{p,1,1}\), \(d_{p,1,2}\), \(d_{p,1,3}\), \(d_{p,22}\) and \(d_{p,3,2}\) do
not play any role in the homological nature of \(C\), but these components are
not independent of the others, because of the relation \(d_{p-1}d_p = 0\).

\subsubsection{A vector field generates a reduction, first proof.}

Let \(C = (C_p, d_p, \beta_p)_p\) be an algebraic cellular complex and \(V =
\{(\sigma_i, \tau_i)\}_{i \in \beta}\) an \emph{admissible} discrete vector
field on \(C\). These data are fixed in this section. We intend to apply the
hexagonal lemma to obtain a canonical reduction of the initial chain complex
\(C\) to a reduced chain complex \(C^c\) where the generators are the critical
cells of \(C\) with respect to~\(V\).

\begin{dfn} ---
\emph{If \(\sigma\) (resp. \(\tau\)) is a \((p-1)\)-cell (resp. a \(p\)-cell)
of \(C\), then the \emph{incidence number} \(\varepsilon(\sigma, \tau)\) is the
coefficient of \(\sigma\) in the differential \(d\tau\). This incidence number
is non-null if and only if \(\sigma\) is a face of \(\tau\); it is \(\pm 1\) if
and only if \(\sigma\) is a \emph{regular} face of \(\tau\).}\qed
\end{dfn}

In particular, for the pairs \((\sigma_i, \tau_i)\) of our vector field, the
relation \(|\varepsilon(\sigma_i, \tau_i)| = 1\) is satisfied.

\begin{dfn} --- \label{68436}
\emph{If \(v = (\sigma_i, \tau_i)\) is a component of our vector field \(V\),
we call \(\sigma_i\) the \emph{source} of \(v\), we call \(\tau_i\) the
\emph{target} of \(v\). We also write \(\tau_i = V(\sigma_i)\) and \(\sigma_i =
V^{-1}(\tau_i)\).}\qed
\end{dfn}

A cell basis \(\beta_p\) is canonically divided by the vector field \(V\) into
three components \(\beta_p = \beta^t_p + \beta^s_p + \beta^c_p\) where
\(\beta^t_p\) (resp. \(\beta^s_p\), \(\beta^c_p\)) is made of the target (resp.
source, critical) cells. For the condition~\ref{18618} of
Definition~\ref{09885} implies a cell cannot be simultaneously a source cell
and a target cell.

The vector field \(V\) defines a bijection between \(\beta^s_{p-1}\) and
\(\beta^t_{p}\). The decompositions of the bases \(\beta_\ast\) induce a
corresponding decomposition of the chain groups \(C_p = C^t_p \oplus C^s_p
\oplus C^c_p\), so that every differential \(d_p\) can be viewed as a \(3
\times 3\) matrix\footnote{Our choice to put \(C^t_p\) before \(C^s_p\) is
intended to produce a situation close to which is described in
Figure~(\ref{36743}). Think the homotopy operator starts from source cells and
go to target cells, while the differential goes in the opposite direction.}.

\begin{prp} --- \label{83464}
Let \(d_{p,2,1}: C^t_{p} \rightarrow C^s_{p-1}\) be the component of the
differential \(d_{p}\) starting from the target cells, going to the source
cells. Then \(d_{p,2,1}\) is an isomorphism.
\end{prp}

\prof If \(\sigma \in \beta^s_{p-1}\), the length of the \(V\)-paths between
\(\sigma\) and the critical cells is bounded. Let us call \(\lambda(\sigma)
\geq 1\) the maximal length of such a path. This length function \(\lambda\) is
a grading: \(\beta^{s,\ell}_{p-1} = \{\sigma \in \beta^s_{p-1} \st
\lambda(\sigma) = \ell\}\). In the same way \(C^{s,\ell}_{p-1} =
\bZ[\beta^{s,\ell}_{p-1}]\), which in turn defines a filtration
\(\widehat{C}^{s,\ell}_{p-1} = \oplus_{k \leq \ell} C^{s,k}_{p-1}\).

The \(V\)-bijection between \(\beta^s_{p-1}\) and \(\beta^t_p\) defines an
isomorphic grading on \(\beta^t_p\) and an analogous filtration on \(C^t_p =
\bZ[\beta^t_p]\); we denote in the same way \(C^{t,\ell}_p =
\bZ[\beta^{t,\ell}_p]\) and \(\widehat{C}^{t,\ell}_p = \oplus_{k \leq \ell}
C^{t,k}_p\).

If \(\sigma \in \beta^{s,\ell}_{p-1}\), that is, if \(\lambda(\sigma) = \ell\),
then the corresponding \(V\)-image \(\tau \in \beta^{t,\ell}_{p}\) has
certainly~\(\sigma\) as a \emph{regular} face and in general other
\emph{source} faces; the grading for these last faces is \(< \ell\) and if
\(\ell > 1\), one of these source faces has the grading \(\ell - 1\), if the
longest \(V\)-path is followed.

A consequence of this description is the following. The partial differential
\(d_{p,2,1}\) can in particular be restricted to \(d_{p,2,1}:
\widehat{C}^{t,\ell}_p \rightarrow \widehat{C}^{s,\ell}_{p-1}\); dividing by
the same map between degrees \((\ell - 1)\) then produces a quotient map
\(\overline{d}_{p,2,1}: C^{t,\ell}_{p} = \bZ[\beta^{t,\ell}_p] \rightarrow
\bZ[\beta^{s,\ell}_{p-1}] = C^{s,\ell}_{p-1}\), which is a \(V\)-isomorphism of
\(\bZ\)-modules \emph{with distinguished basis}.

The standard recursive argument then proves \(d_{p,2,1}\) is an
isomorphism.\qed

\begin{thr} \textbf{\emph{(Vector-Field Reduction Theorem)}} --- \label{81450}
Let \(C = (C_p,d_p,\beta_p)_p\) be an algebraic cellular complex and \(V =
\{(\sigma_i, \beta_i)\}_{i \in \beta}\) be an \emph{admissible} discrete vector
field on \(C\). Then the vector field \(V\) defines a canonical reduction
\(\rho = (f,g,h): (C_p, d_p) \rrdc (C^c_p, d'_p)\) where \(C^c_p =
\bZ[\beta^c_p]\) is the free \(\bZ\)-module generated by the critical
\(p\)-cells.
\end{thr}

\prof Theorem~\ref{09803} + Proposition~\ref{83464}.\qed

The bottom chain complex \((C^c_\ast, d'_\ast)\) and its image \(g(C^c_\ast)
\subset C_\ast\) are both versions of the Morse complex described in Forman's
paper~\cite[Section 7]{frmn}.

The inverse \(d_{p,2,1}^{-1}\) is crucial when concretely using this
proposition to solve the homological problem for \((C_p, d_p)\) thanks to a
solution of the same problem for the \emph{critical complex} \((C^c_p, d'_p)\):
this matrix can be -- and will sometimes be in this text -- not of finite type!

If \(\tau \in \beta^t_p\), the value \(d_{p,2,1}(\tau)\) is defined by:
\begin{equation}
d_{p,2,1}(\tau) = \sum_{\sigma \in \beta^s_{p-1}} \varepsilon(\sigma, \tau)
\sigma
\end{equation}
where \(\varepsilon(V^{-1}(\tau), \tau) = \pm 1\), but the other coefficients
are arbitrary. We so obtain the recursive formula:
\begin{equation} \label{08522}
d_{p,2,1}^{-1}(\sigma) = \varepsilon(\sigma, V(\sigma))\,\left(V(\sigma) -
\sum_{\sigma' \in \beta^s_{p-1} - \{\sigma\}} \varepsilon(\sigma', V(\sigma))
\,d_{p,2,1}^{-1}(\sigma')\right)
\end{equation}
with \(\varepsilon(\sigma, V(\sigma)) = \pm 1\), the other incidence numbers
being arbitrary.

 This formula is easily recursively programmed.

\subsubsection{Using the Homological Perturbation Theorem.} \label{32626}

We give now a direct proof of Theorem~\ref{81450} based on the Homological
Perturbation Theorem~\ref{07404}.

\prof Instead of considering the right differential \(d\) of our chain complex
\((C_\ast, d)\), we start with a different chain complex \((C_\ast, \delta)\);
the underlying graded module \(C_\ast\) is the same, but the differential is
``simplified''. The new differential \(\delta_p: C_p \rightarrow C_{p-1}\) is
roughly defined by the formula \(\delta = \varepsilon V^{-1}\).

  More precisely, we take account of the canonical decomposition defined by the
vector field \(V\):
\begin{equation}
 C_p = C_p^t \oplus C_p^s \oplus C_p^c = \bZ[\beta_p^t] \oplus \bZ[\beta_p^s]
 \oplus \bZ[\beta_p^c]
\end{equation}
as follows. If \(\tau \in \beta_p\) is an element of \(\beta^s_p\) or
\(\beta^c_p\), then we decide \(\delta_p(\tau) = 0\). If \(\tau \in
\beta^t_p\), that is, if \(\tau\) is a target cell, then there is a unique
vector \((\sigma, \tau) \in V\) with \(\tau\) as the second component and we
define \(\delta_p(\tau) = \varepsilon(\sigma, \tau) \sigma\ \mbox{``=''}\
(\varepsilon V^{-1})(\tau)\).

The homotopy operator \(h = \varepsilon V\) defined in the same way in the
reverse direction obviously defines an initial \emph{reduction} \((C_p,
\delta_p) \rrdc (C_p^c, 0)\) of this initial chain complex over the chain
complex with a null differential generated by the critical cells. Look again at
the diagram~(\ref{36743}): the reverse diagonal arrows \(d\) and \(h\) of the
diagram correspond to our simplified differential \(\delta\) and the homotopy
operator \(h\), while the horizontal arrows \(d\) are null.

Restoring the right differential \((d_p)_p\) on \((C_p, d_p)_p\) can be
considered as introducing a \emph{perturbation} \((d_p - \delta_p)_p\) of the
differential \((\delta_p)_p\). The basic perturbation lemma can be applied if
the \emph{nilpotency condition} is satisfied, that is, if \((d - \delta) \circ
h\) is pointwise nilpotent. The perturbation is nothing but the right
differential except for the target cells; if \(\tau\) is such a target cell,
then \((d - \delta)(\tau) = d(\tau) - \varepsilon(\sigma, \tau) \sigma\) if
\((\sigma, \tau)\) is the corresponding vector.

Let us examine this composition \((d - \delta)h\). It is non-trivial only for a
source cell \(\sigma\) and \(h(\sigma) = V(\sigma)\). Because of the definition
of \(\delta\), the value of \((d - \delta)h(\sigma)\) is \(d(\tau) - \sigma\),
it is made of the faces of \(\tau\) \emph{except the starting cell} \(\sigma\).
Considering next \(h(d - \delta)h(\sigma)\) amounts to selecting the source
cells of \(d(\tau) - \sigma\) and applying again the vector field~\(V\). You
understand we are just following all the \emph{\(V\)-paths} starting from the
initial source cell \(\sigma\), see Definition~\ref{76422}. If the vector field
is admissible, the length of these V-paths is bounded, the process terminates,
\mbox{\((d - \delta)h\)} is pointwise nilpotent and the \emph{nilpotency
condition}, required to apply the Homological Perturbation Theorem, is
satisfied.

It happens the explicit formulas for the new reduction produced by this theorem
are exactly the formulas~(\ref{76423}) combined with the recursive
definition~(\ref{08522}) of the homotopy operator \(h'\) of the perturbed
reduction.\qed

This proof is not only more direct, but also much easier to program, at least
if the Homological Perturbation Theorem is implemented in your programming
environment.

\begin{dfn} ---
\emph{A \emph{W-reduction} is the reduction produced by Theorem~\ref{81450}
when an algebraic cellular complex is provided with an admissible discrete
vector field.}\qed
\end{dfn}

\section{Computing vector fields.}

\subsection{Introduction.}

A vector field is most often an \emph{initial} tool to construct an interesting
\emph{reduction}, deduced from the vector field by the Vector-Field Reduction
Theorem. For example we will see the so fundamental Eilenberg-Zilber reduction
is a direct consequence of a simple vector field. The plan is then:
\begin{equation}
\mbox{Vector field} \Longrightarrow \mbox{Reduction}
\end{equation}

But it happens most interesting results of this text in fact have been obtained
with the reverse plan:
\begin{equation}
\mbox{Reduction} \Longrightarrow \mbox{Vector field}
\end{equation}

The story is the following. The \emph{discrete vector field} point of view is
relatively recent. In many situations, the \emph{reductions} we are interested
in are known for a long time. The problem is then the following: is it possible
in fact to obtain such a reduction thanks to an appropriate discrete vector
field?  We explain in this section there exists a systematic method to
\emph{deduce} such a vector field from the reduction, of course \emph{if it
exists}, which is not necessarily the case.

The reader can then wonder which can motivate such a research, when the
corresponding reduction is already known! The point is the following: The
\emph{already known} reduction has therefore been defined by a different
process, often relatively sophisticated. These reductions are important in
\emph{constructive} algebraic homology, and the algorithms implemented them can
be complex, therefore time and space consuming.

Several experiments then show that if a vector field giving directly such a
reduction is finally found, then the algorithm computing this reduction can be
much more efficient, in time \emph{and} space complexity, in particular for the
terrible homotopy operator, so important from a constructive point of view,
almost always neglected by ``classical'' topologists. The algorithm explained
here to identify a vector field possibly defining a \emph{known} reduction can
therefore be considered as a program systematically \emph{improving} other
programs \emph{previously written}.

There are also some cases where some reduction is \emph{conjectured}, but not
yet proved. We obtain in such a context a relatively spectacular result:
Eilenberg and MacLane sixty years ago conjectured\footnote{Eilenberg and
MacLane's terminology is a little confusing: what \emph{we} call a
\emph{reduction} is a \emph{contraction} in \emph{their} paper, while
\emph{their reduction} is a chain complex morphism inducing a weak homology
equivalence.} the existence of a direct reduction \(C_\ast(BG) \rrdc
\Bar(C_\ast(G))\) for a simplicial group \(G\), see the comments after the
statement of~\cite[Theorem~20.1]{elmc1}. Numerous calculations using repeatedly
the Homological Perturbation Theorem led us to suspect we had obtained such a
direct reduction, but we were unable to prove it. Using the below algorithm, we
easily obtained a vector field producing this \emph{possible} reduction. It
happens the structure of this vector field is very easy to understand, which
this time produces a \emph{proof} of the Eilenberg-MacLane conjecture.
Furthermore, the calculation of the \emph{effective} homology of the
Eilenberg-MacLane spaces is so much more efficient, a key point for the quick
calculation of homotopy groups of \emph{arbitrary} simply connected spaces.

\subsection{The unique possible vector field.}

We consider in this section a fixed \emph{given} reduction:
\begin{equation}
\rho = (f,g,h): \wC_\ast \rrdc C_\ast
\end{equation}
between two \emph{cellular} chain complexes. And we want to study whether some
admissible discrete vector field could ``explain'' this reduction.

The problem is the following: we have to divide the cells elements of
\(\beta_\ast\) in \(\wC_\ast\) in \emph{source}, \emph{target} and
\emph{critical} cells, and we have also to define a coherent \emph{pairing}
between source and target cells. In such a way the reduction then obtained by
the Vector-Field Reduction Theorem~\ref{81450} is just the given reduction
\(\rho\).

The key point is the collection of formulas~(\ref{76423}).

Observe firstly the homotopy operator \(h\) is non-null only for source cells.
More precisely, a cell \(\sigma \in \beta_p\) is a source cell if and only if
\(h(\sigma) \neq 0\). No circularity in a V-path, for the vector field must be
\emph{admissible}. So that if \((\sigma, \tau)\) is the vector starting from
\(\sigma\), the occurence of \(\tau = V(\sigma)\) in the formula~(\ref{08522})
cannot be cancelled by other terms. So the reduction \(\rho\) gives a decision
algorithm for the source property of a cell.

Let us now consider a cell \(\sigma\) which has been proved being a source
cell. What is the corresponding target cell? The value \(h(\sigma)\) is a
linear combination of cells, and at least one of them has \(\sigma\) in its
boundary. In fact only one, for again, otherwise there will be a circular
V-path. So no choice for the corresponding target cell \(\tau\). If the
incidence number \(\varepsilon(\sigma, \tau)\) is not \(\pm 1\), possible, no
vector field is possible. If such a vector field on the contrary exists, we
have an algorithm defining the pairing source cell \(\mapsto\) target cell.

In the finite case, this is enough. The first algorithm gives the list of the
source cells, the second one gives the list of the target cells, and the
remaining cells are the critical ones. There remains to verify the
formulas~(\ref{08522}).

In the infinite case, in general there is no algorithm allowing one to decide
whether a given reduction comes from a vector field. Because of the usual
obstacle forbidding for example to solve the halting problem.

Nevertheless you can also decide whether \emph{some} cell \(\tau\) is a target
cell. You examine firstly if it is a source cell, and if it is, it is not a
target cell. If \(\tau\) is not a source cell, then you consider the boundary
\(d\tau = \sum \alpha_i \sigma_i\), it is a finite linear combination, and you
just have to examine if some \(\sigma_i\) in the boundary is a source cell, and
if the corresponding target cell is this cell \(\tau\).

Finally if a cell is neither a source cell nor a target cell, it is a critical
cell. Of course if the rank of your chain groups are infinite, you cannot
verify the formulas~(\ref{08522}) for all the cells.

Another point is useful. In most interesting reductions, the small chain
complex \(C_\ast\) is of finite type, even if the big one \(\wC_\ast\) is not.
It is then important to identify and to understand the ``nature'' of the
critical cells. These cells can be easily obtained thanks to the
formula~(\ref{76423}) for \(g_p\). Take some generator \(a\) of the small
complex \(C_\ast\), apply the lifting operator \(g: C_\ast \rightarrow
\wC_\ast\) to obtain \(g(a)\); the formula for \(g\) proves this image is made
of exactly \emph{one} critical cell and an arbitrary combination of source
cells; so that running the generators present in \(g(a)\) allows you to
identify the critical cell of \(\wC_\ast\) corresponding to the generator \(a\)
of \(C_\ast\).

Once these simple tricks are applied, the game becomes the following, when it
makes sense. Examining the ``general style'' of as many source, target and
critical cells as reasonably possible, you have to guess the unknown vector
field in fact your reduction \(\rho\) has ``followed'' to be defined.
Experience shows it can be quite amusing. We will see this rule is quite simple
for the Eilenberg-Zilber reduction, at least if the good point of view is
understood, a little less simple for the direct reduction \(C_\ast(BG) \rrdc
\Bar(C_\ast(G))\), very sophisticated for the direct reduction \(C_\ast(\Omega
X) \rrdc \Cobar(C_\ast(X))\).

Before defining the Eilenberg-Zilber vector field in Section~\ref{08795}, we
will play this game with the Eilenberg-Zilber reduction obtained sixty years
ago.

\section{Two simple examples.}

This section is devoted to two classical results of algebraic topology proved
here thanks to simple discrete vector fields. The first one is the
\emph{normalization theorem} for simplicial homology. The second one is devoted
to the ``simplest'' Eilenberg-MacLane space, \(K(\bZ,1)\), simple but
unfortunately not of finite type! So that a homological reduction \(K(\bZ,1)
\rrdc S^1\) on the \emph{very small} standard circle, one vertex and one edge,
is quite useful.

\subsection{Simplicial sets and their chain complexes.} \label{84511}

The reader is assumed to be familiar with the elementary definitions,
properties and results about the \emph{simplicial sets}. Not to be confused
with simplicial \emph{complexes}. The main references are maybe~\cite{may,
grjr}; the notes~\cite[Section~7]{rbsr09} or~\cite{srgr01} can also be useful.

The basic category in this context is the category \(\uD\). An object
\(\underline{p}\) of \(\uD\) is the set of the integers \(0 \leq i \leq p\),
also denoted by \([0\ldots p]\). A \(\uD\)-morphism \(\alpha: \underline{p}
\rightarrow \underline{q}\) is an increasing function: \(i \leq j \Rightarrow
\alpha(i) \leq \alpha(j)\).

Let \(X\) be a simplicial set. For every natural number \(p \in \bN\) the set
of \(p\)-simplices~\(X_p\) is defined. For every \(\uD\)-morphism \(\alpha:
\underline{p} \rightarrow \underline{q}\), a corresponding map \(\alpha^\ast:
X_q \rightarrow X_p\) is defined. The simplicial set \(X\) can be viewed as a
\emph{contravariant} functor \(X: \uD \rightarrow \underline{\textrm{Set}}\).

The face \(\uD\)-morphisms \(\partial^p_i: [0\ldots(p-1)] \rightarrow [0\ldots
p]\) are defined for \(p \geq 1\) and \(0 \leq i \leq p\). The (elementary)
degeneracy \(\uD\)-morphisms \(\eta^p_i: |0\ldots(p+1)] \rightarrow [0\ldots
p]\) are defined for \(0 \leq i \leq p\).
\begin{equation}
\begin{tikzpicture} [yscale = 0.3, baseline = (b5)] \scriptsize
 \fill [black!20, rounded corners] (-1.3,-1) rectangle (2.2,10) ;

 \node (a9) at (0,9) [label=left:\(0\)] {\(\bullet\)} ;
 \node (b9) at (1,9) [label = right: \(0\)] {\(\bullet\)} ;
 \draw [->] (a9) -- (b9) ;
 \node (a8) at (0,8) [label=left:\(1\)] {\(\bullet\)} ;
 \node (b8) at (1,8) [label = right: \(1\)] {\(\bullet\)} ;
 \draw [->] (a8) -- (b8) ;
 \node (a5) at (0,5) [label=left:\(i-1\)] {\(\bullet\)} ;
 \node (b5) at (1,5) [label = right: \(i-1\)] {\(\bullet\)} ;
 \draw [->] (a5) -- (b5) ;
 \draw [very thick, dotted] (a8) -- (a5) ;
 \draw [very thick, dotted] (b8) -- (b5) ;
 \node (a4) at (0,4) [label=left:\(i\)] {\(\bullet\)} ;
 \node (b4) at (1,4) [label = right: \(i\)] {\(\bullet\)} ;
 \node (b3) at (1,3) [label = right: \(i+1\)] {\(\bullet\)} ;
 \draw [->] (a4) -- (b3) ;
 \node (a1) at (0,1) [label=left:\(p-1\)] {\(\bullet\)} ;
 \node (b0) at (1,0) [label = right: \(p\)] {\(\bullet\)} ;
 \draw [very thick, dotted] (a4) -- (a1) ;
 \draw [very thick, dotted] (b3) -- (b0) ;
 \draw [->] (a1) -- (b0) ;
 \normalsize
 \node at (-2,5) {\(\partial^p_i =\)} ;
\end{tikzpicture}
\hspace{30pt}
\begin{tikzpicture} [yscale = 0.3, baseline = (b5)] \scriptsize
 \fill [black!20, rounded corners] (-1.3,-1) rectangle (2.2,10) ;
 \node (a9) at (0,9) [label=left:\(0\)] {\(\bullet\)} ;
 \node (b9) at (1,9) [label = right: \(0\)] {\(\bullet\)} ;
 \draw [->] (a9) -- (b9) ;
 \node (a8) at (0,8) [label=left:\(1\)] {\(\bullet\)} ;
 \node (b8) at (1,8) [label = right: \(1\)] {\(\bullet\)} ;
 \draw [->] (a8) -- (b8) ;
 \node (a5) at (0,5) [label=left:\(i-1\)] {\(\bullet\)} ;
 \node (b5) at (1,5) [label = right: \(i-1\)] {\(\bullet\)} ;
 \draw [->] (a5) -- (b5) ;
 \draw [very thick, dotted] (a8) -- (a5) ;
 \draw [very thick, dotted] (b8) -- (b5) ;
 \node (a4) at (0,4) [label=left:\(i\)] {\(\bullet\)} ;
 \node (b4) at (1,4) [label = right: \(i\)] {\(\bullet\)} ;
 \node (a3) at (0,3) [label = left: \(i+1\)] {\(\bullet\)} ;
 \draw [->] (a4) -- (b4) ;
 \draw [->] (a3) -- (b4) ;
 \node (a0) at (0,0) [label=left:\(p+1\)] {\(\bullet\)} ;
 \node (b1) at (1,1) [label = right: \(p\)] {\(\bullet\)} ;
 \draw [very thick, dotted] (a3) -- (a0) ;
 \draw [very thick, dotted] (b4) -- (b1) ;
 \draw [->] (a0) -- (b1) ;
 \normalsize
 \node at (-2,5) {\(\eta^p_i =\)} ;
\end{tikzpicture}
\end{equation}
An arbitrary composition of (elementary) degeneracies \(\eta = \eta^q_{i_{p-q}}
\eta^{q+1}_{i_{p-q-1}} \cdots \eta^{p-1}_{i_1}\) is a \emph{degeneracy}; this
expression is unique if the inequalities \(i_{p-q} < \cdots < i_1\) are
required. The degeneracies are nothing but the \emph{surjective}
\(\uD\)-morphisms.

Most often the sup-index of the face and degeneracy operators are omitted, we
write simply \(\partial_i\) (resp. \(\eta_i\)) instead of \(\partial^p_i\)
(resp. \(\eta^p_i\)).

For every simplicial set \(X\), every \(p > 0\) and every \(0 \leq i \leq p\),
a face operator \(\partial_i: X_p \rightarrow X_{p-1}\) is defined, this is the
map which applies a simplex to its \mbox{\(i\)-th face}, thought of as opposite
to the \(i\)-th vertex; an elementary degeneracy operator \(\eta_i: X_p
\rightarrow X_{p+1}\) is defined for \(0 \leq i \leq p\), it is the map which
applies a simplex to the same but the \(i\)-th vertex is ``repeated'', replaced
by a degenerate edge, increasing the dimension by 1. Every \(\uD\)-morphism is
a composition of face and (elementary) degeneracy operators, so that the face
and degeneracy operators between the simplex sets \(X_p\) are enough to define
the simplicial structure, at least if they satisfy appropriate compatibility
conditions.

\begin{dfn} ---
\emph{A \(p\)-simplex \(\sigma \in X_p\) is \emph{degenerate} if there exists
an integer \(q < p\), a \(\uD\)-morphism \(\alpha: \underline{p} \rightarrow
\underline{q}\) and a \(q\)-simplex \(\tau \in X_q\) satisfying \(\sigma =
\alpha^\ast \tau\). We denote by \(X\sD_p\) (resp. \(X\sND_p\)) the set of
degenerate (resp. non-degenerate) \(p\)-simplices.}\qed
\end{dfn}

The Eilenberg-Zilber lemma gives for every simplex a canonical expression from
a unique non-degenerate simplex.

\begin{thr} \emph{\textbf{(Eilenberg-Zilber lemma)}} --- \label{34608}
  Let \(\sigma\) be a \(p\)-simplex of a simplicial set \(X\). Then there
  exists a unique triple \((q, \eta, \tau)\), the \emph{Eilenberg triple} of
  \(\sigma\), satisfying:
\begin{fenum}
\item
\(0 \leq q \leq p\).
\item
\(\eta: \underline{p} \rightarrow \underline{q}\) is a surjective
\(\uD\)-morphism.
\item
\(\tau \in X_q\) is non-degenerate and \(\alpha^\ast \tau = \sigma\).
\end{fenum}
\end{thr}

\prof \cite[(8.3)]{elzl1} \qed

In short, every simplex comes from a unique non-degenerate simplex, by a unique
degeneracy, and every non-degenerate simplex \(\sigma \in X_p\) generates a
collection of degenerate simplices in any dimension \(> p\): \(\{\eta^\ast
\sigma\}_{\eta \in \uD^{\textrm{\tiny surj}}(-,p) - \{\textrm{\scriptsize
id}_{[0\ldots p]}\}}\).

In a sense the degenerate simplices are somewhat \emph{redundant} and the
\emph{normalization theorem} explains you can neglect them when defining -- and
\emph{computing} -- the homology of the underlying simplicial set.

\begin{dfn} ---
\emph{Let \(X\) be a simplicial set. The (non-normalized) chain complex
\(C_\ast X\) associated to \(X\) is the algebraic cellular complex \((\bZ[X_p],
d_p, X_p)_p\) where the differential \(d_p\) is defined by \(d_p \sigma =
\sum_{i = 0}^p (-1)^i
\partial_i \sigma\).\\[\parskip]
\hspace*{\parindent}The normalized chain complex \(C\sND_\ast X\) is the
algebraic cellular complex \((\bZ[X\sND_p], d_p, X\sND_p)_p\), using only the
non-degenerate cells; the differential is defined by the same formula, except
every degenerate face is cancelled.}\qed
\end{dfn}

It is a little better to firstly observe the degenerate simplices generate a
chain \emph{sub}complex \(C\sD_\ast X\) and the normalized chain complex is
nothing but the quotient \(C\sND_\ast X = C_\ast X / C\sD_\ast X\). The
following result is classical and fundamental.

\begin{thr} \emph{\textbf{(Normalization Theorem)}} ---
The projection \(C_\ast X \rightarrow C\sND_\ast X\) is a homology
equivalence.\qed
\end{thr}

The proof is not difficult, see for example~\cite[Section VIII.6]{mcln2}.

\subsection{A discrete vector field proof of the Normalization Theorem.}

The standard proof of the Normalization Theorem is not difficult but however
requires some lucidity. We intend to give a slightly different proof based on a
discrete vector field. Consider this as an opportunity to illustrate our
technique in a very simple case. Maybe you will find this proof simpler and,
why not, funnier. The announced homology equivalence is a direct consequence of
the next proposition.

\begin{prp} --- \label{77213}
The degenerate cellular complex \(C_\ast\sD\) is acyclic.
\end{prp}

\prof Let us recall a degenerate simplex \(\sigma \in X\sD_p\) of a simplicial
set \(X\) has a unique expression: \(\sigma = \eta \tau\) with \(\tau \in
X_q\sND\) a non-degenerate \(q\)-simplex, with \(q < p\), and \(\eta:
\underline{p} \rightarrow \underline{q}\) a \(\uD\)-surjection. We represent
such a surjection as an increasing sequence of \((p+1)\) integers of
\(\underline{q}\). For example 0122234 represents the unique surjection \(\eta:
\underline{6} \rightarrow \underline{4}\) satisfying \(\eta(2) = \eta(3) =
\eta(4) = 2\).

We will construct an \emph{admissible} discrete vector field on \(C\sD_\ast\)
where \emph{every} cell is source or target of the vector field: no remaining
cell, that is, no critical cell, which at once proves our complex is acyclic:
homology equivalent to the null complex.

We must divide the set of cells in two ``equal'' parts in bijection by the
vector field to construct. The process is simple. Necessarily, because a
degeneracy operator \(\eta\) is a non-injective surjection, some images are
repeated. We take the first occurence of a repetition, which has a
\emph{multiplicity}: for example the multiplicity of the first repetition of
012222345556 is 4. If this multiplicity is \emph{even}, we decide the
corresponding cell is a source of the vector field to be defined, and the
corresponding target is the cell obtained by adding an extra repetition to this
first one. For example if \(\sigma = \eta \tau\) with \(\eta = 012222345556\),
we decide this \(\sigma\) is a source and the corresponding target is \(\sigma'
= \eta' \tau\) with \(\eta' = 0122222345556\), first multiplicity 5; this new
simplex is also degenerate.

Symmetrically, if the multiplicity of the first repetition is odd, therefore at
least~3, we decide the cell is a target, and removing one repetition gives back
the corresponding source. The division of the degenerate cells in two ``equal''
parts is clear, how to be simpler?

If a cell is a source, then applying the boundary formula shows this source
actually is a \emph{regular} face of the corresponding target. Observe you
cannot exchange the choices of odd and even when processing the number of
repetitions, otherwise the regularity condition would not hold.

The most important remains to do: we must prove the \emph{admissibility} of
this vector field. We construct a Lyapunov function, see
Definition~\ref{92078}. Instead of a general definition for the desired
Lyapunov function, we prefer a unique example which illustrates all the cases
to be considered. Let~\(\tau\) be some non-degenerate 4-simplex and \(\sigma =
\eta \tau\) with \(\eta = 01123334 = \eta_1\eta_4\eta_5\). We define the value
\(L_p(\sigma) = 4\), the ``genuine'' dimension of \(\sigma\), that is the
dimension of the non-degenerate simplex~\(\tau\) associated to \(\sigma\) by
the Eilenberg-Zilber lemma.

We defined \(V(\sigma) = \sigma' = (011123334)\tau = \eta_6(\eta_5(\eta_2(
\eta_1(\tau)))) = (\eta_1 \eta_2 \eta_5 \eta_6)\tau\), do not forget a
simplicial set is a \emph{contravariant} functor. We must consider the faces
\(\partial_i \sigma'\) for \(0 \leq i \leq 8\). In fact the faces 1, 2 and 3
are not to be considered, for such a face is the initial cell~\(\sigma\), the
incidence number \(\varepsilon(\sigma, \sigma')\) being in this case -1: from
the \(algebraic\) point of view, it is a unique face, a regular one. The faces
5, 6, and 7 are \((01112334)\tau\), not to be considered either, for this is
not a source cell: the multiplicity of 1 is 3, odd. There remains the faces 0,
4 and 8.

For example \(\partial_4 \sigma' = \partial_4 \eta_6 \eta_5 \eta_2 \eta_1 \tau
= \eta_5 \eta_4 \eta_2 \eta_1 \partial_2 \tau\) and the genuine dimension of
this simplex is \(\leq 3\). It is possible this face is not a source cell, but
anyway, if it is, the value of the Lyapunov function has decreased. The same
for any face operator which is not swallowed by a degeneracy, in this case, for
the faces 0, 4 and 8. We have proved that for any face \(\sigma''\) of
\(V(\sigma) = \sigma'\), except the initial face \(\sigma\), if \(\sigma''\) is
a source cell, then the inequality \(L_p(\sigma'') < L_p(\sigma)\) holds.

The Lyapunov function \(L_p\) can be used for the bounding function
\(\lambda_p\) required by Definition~\ref{96757}. Our vector field is
admissible and the proposition is proved.\qed

\subsection
  [\texorpdfstring{The Eilenberg-MacLane space $K(\bZ,1)$.}
  {The Eilenberg-MacLane space K(Z,1).}]
  {The Eilenberg-MacLane space \boldmath\(K(\bZ,1)\).}

  It is an abelian simplicial group defined as follows. The set of
\(n\)-simplices is simply the abelian group \(\bZ^n\), and a simplex is
traditionally denoted as a \emph{bar}-object \(\sigma = [a_1|\cdots|a_n]\). In
particular only one vertex (0-cell), the void bar object~\([\,]\), while for
every positive \(n\), the simplex set \(\bZ^n\) is infinite. The face operators
\(\partial_i: \bZ^n \rightarrow \bZ^{n-1}\) are defined as follows:
\begin{equation}
\begin{array}{rcl}
\partial_0([a_1|\cdots|a_n]) &=& [a_2|\cdots|a_n]
\\
\partial_i([a_1|\cdots|a_n]) &=&
[a_1|\cdots|a_{i-1}|(a_i+a_{i+1})|a_{i+2}|\cdots|a_n]\ \ \mbox{if \(1 < i <
n\)}
\\
\partial_n([a_1|\cdots|a_n]) &=& [a_2|\cdots|a_{n-1}]
\end{array}
\end{equation}

The degeneracy operator \(\eta_i: \bZ^n \rightarrow \bZ^{n+1}\) consists in
putting an extra null component in position \(i\):
\begin{equation}
\eta_i([a_1|\cdots|a_n]) = [a_1|\cdots|a_{i}|\,0\,|a_{i+1}|\cdots|a_n]
\end{equation}
so that the collection of non-degenerate simplices \(K^{ND}_{n}\) is
\(\bZ_\ast^n\): no component of the bar expression must be null.

This definition can be generalized to any group \(G\), non-necessarily abelian,
then producing the simplicial set \(K(G,1)\). The topological
\emph{realization} \(X = |K(G, 1)|\) of this simplicial object has a
characteristic property up to homomotopy equivalence: all its homotopy groups
are null except \(\pi_1 X = G\), but this fact is not used here. In general,
\(K(G,1)\) \emph{is not} a simplicial \emph{group}; but if \(G\) is abelian,
then \(K(G,1)\) is an abelian simplicial group.

It happens \(K(\bZ, 1)\) is the minimal model of the circle \(S^1\), minimal in
the sense of \mbox{Kan~\cite[\S 9]{may}}. Not really concretely minimal in the
ordinary sense: as many non-degenerate \(n\)-cells as elements in
\(\bZ_\ast^n\). On one hand such a big model is unavoidable if you want to
\emph{compute} the homotopy groups of an \emph{arbitrary} simplicial set; on
the other hand a \emph{reduction} is easily obtained \(C_\ast K(\bZ, 1) \rrdc
C_\ast S^1\) where this time \(S^1\) is the usual model of the circle with one
vertex \(\ast\) and one (non-degenerate) edge \(s^1\), so that the (normalized)
chain complex of \(C_\ast S^1\) is simply: \( 0 \leftarrow \bZ[\ast]
\stackrel{0}{\leftarrow} \bZ[s^1] \leftarrow 0\).

The \emph{ordinary} homology of \(K(\bZ,1)\) is obvious: \(\bZ\) in degrees 0
and 1, nothing else. But in \emph{effective} homology, it is \emph{mandatory}
to keep in one's environment the reduction \(C_\ast K(\bZ,1) \rrdc C_\ast
S^1\). It happens this reduction can be deduced from a discrete vector field.

\begin{prp} ---
An \emph{admissible discrete vector field} \(V\) can be defined on
\(C_\ast^{ND}(K(\bZ, 1))\) as follows:
\begin{fenum}
\item
The critical cells are the unique \(0\)-simplex \([\,]\) and the \(1\)-simplex
\([1]\).
\item
The source cells are the non-critical cells \([a_1 | \cdots]\) satisfying \(a_1
\neq 1\).
\item
The target cells are the cells \([a_1 | a_2|\cdots]\) of dimension \(\geq 2\)
satisfying \(a_1 = 1\).
\item
The pairing \emph{[source cell \(\leftrightarrow\) target cell]} associates to
the source cell \([a_1|a_2|\cdots]\) the target cell \([1 | a_1-1 | a_2 |
\cdots]\) if \(a_1 > 1\) and \([1 | a_1 |a_2 | \cdots]\) if \(a_1 < 0\).
\end{fenum}
\end{prp}

\prof The admissibility property comes from the following observation: any
V-path decreases the absolute value of the first component \(a_1\). Two
examples:
\begin{equation}
\begin{array}{l}
 [3|6] \mapsto [1|2|6] \mapsto [2|6] \mapsto [1|1|6] \mapsto \mbox{halt!}
 \\
 {}[-3|6] \mapsto [1|-3|6] \mapsto [-2|6] \mapsto [1|-2|6] \mapsto [-1|6] \mapsto
 [1|-1|6] \mapsto \mbox{halt!}
\end{array}
\end{equation}

For example the four faces of \([1|2|6]\) are \([2|6]\), a source cell which is
the continuation of the path, \([3|6]\), the source cell matching \([1|2|6]\),
not to be considered, \([1|8]\) not a source cell, and \([1|2]\) not a source
cell. For the last target cell of this path, only one face \([2|6]\) is a
source cell but it just matches \([1|1|6]\), so that the path cannot be
continued. We let the reader play the same game with the second example.

The critical chain complex is therefore \(0 \leftarrow \bZ[\,]
\stackrel{d'}{\leftarrow} \bZ[1] \leftarrow 0\), but what is the differential
\(d'\)? We must use the first formula~(\ref{76423}). The differential \(d[1] =
\partial_0[1] - \partial_1[1] = [\,] - [\,] = 0\) is null, so that \(d_{p,3,3}\)
and \(d_{p,2,3}\) are null as well, and finally \(d' = 0\).\qed

Our two examples are in a sense symmetric. In the first case, for the
normalization theorem, the studied complex is an acyclic subcomplex, so that
the reduction \(C_\ast X \rrdc C^{ND}_\ast X\) has a quotient complex as the
small complex. In the second case \(K(\bZ, 1)\), the critical cells \([\,]\)
and \([1]\) generate a subcomplex whcih then is always the critical complex,
for the components \(d_{p,2,3}\) of the boundary matrices are null; it is the
quotient by this subcomplex which is acyclic.

It is interesting to compare the previous proposition with the formula (14.4)
of~\cite{elmc2}: which shows Eilenberg and MacLane already had a perfect
knowledge of our vector field, sixty years ago, even if they did not find
useful to introduce the corresponding terminology.

\section{W-reductions of digital images.}

\subsection{Introduction.}

\begin{dfn} ---
\emph{A \emph{digital image}, an \emph{image} in short, is a \emph{finite}
algebraic cellullar complex \((C_p, d_p, \beta_p)_{p \in \bZ}\): every
\(\beta_p\) is finite and furthermore every \(\beta_p\) is empty outside an
interval \([0\ldots P]\), the smallest possible \(P\) being the
\emph{dimension} of the image.}\qed
\end{dfn}

For example the various techniques of scientific imaging produce \emph{images},
some finite objects, typically finite sets of pixels. If you are interested in
some \emph{homological} analysis of such an image, you associate to it a
\emph{geometrical} cellular complex, again various techniques can be used, and
finally this defines an algebraic cellular complex. There remains to compute
the homology groups of this complex; in fact computing the \emph{effective}
homology~\cite{rbsr06} of this complex is much better: think for example of
this interesting notion of \emph{persistent} homology where the homology groups
are not enough, you must \emph{exhibit} cycles representing the homology
classes, for example to study if such a cycle ``remains alive'' when the image
is modified; it so happens effective homology is exactly designed to construct
such cycles.

The Vector-Field Reduction Theorem (Theorem~\ref{81450}) is then particularly
welcome. The bases \(\beta_p\) of the initial complex can be enormous, but
appropriately choosing a vector field can produce, applying this Reduction
Theorem, a new complex which is homology equivalent, with small \emph{critical}
bases \(\beta^c_p\); so that the homology computations are then fast. More
important: because this result produces a \emph{reduction} between both
complexes, the fast solution of the homological problem for the critical
complex gives also at once a solution for the same problem in the big initial
complex. Let us recall these methods of \emph{effective} homology were designed
for initial complexes \emph{not of finite type}; which is successful for
infinite objects should also work for big ones, big but finite! He who can do
more can do less.

The simplest case is the case of a chain complex with only two consecutive
chain groups, and the general case is easily reduced to this one. It is the
problem known as the reduction of an integer matrix to the Smith form about
which much work has already been done, often impressive and fascinating. The
matter here is just of taking account of the frequent presence of terms \(\pm
1\) in the matrices produced by imaging, and to examine if this method of
vector fields can be used. It is nothing but the following game: let \(M\) be a
given integer matrix; please extract a square triangular submatrix, as large as
possible, where all the diagonal terms are \(\pm 1\); triangular with respect
to some order of rows and columns to be appropriately chosen. This has
certainly already been extensively studied by the specialists in imaging; we
are not such specialists and we do not pretend invent anything. But this point
of view of \emph{effective} homology could after all be of some usefulness.

\subsection{Vector fields and integer matrices.}

Let \(M\) be an initial matrix \(M \in \textrm{Mat}_{m,n}(\bZ)\), with \(m\)
rows and \(n\) columns. Think of \(M\) as the unique non-null boundary matrix
of the chain complex:
\begin{equation}
\cdots \leftarrow 0 \leftarrow \bZ^m \stackrel{M}{\longleftarrow} \bZ^n
\leftarrow 0 \leftarrow \cdots
\end{equation}

A \emph{vector field} \(V\) for this matrix is nothing but a set of integer
pairs \(\{(a_i, b_i)\}_i\) satisfying these conditions:
\begin{fenum}
\item
\(1 \leq a_i \leq m\) and \(1 \leq b_i \leq n\).
\item
The entry \(M[a_i, b_i]\) of the matrix is \(\pm 1\).
\item
The indices \(a_i\) (resp. \(b_i\)) are pairwise different.
\end{fenum}

This clearly corresponds to a vector field, and constructing such a vector
field is very easy. But there remains as usual the main problem: is this vector
field \emph{admissible}? An interesting but serious problem!

Because the context is finite, it is just a matter of avoiding \emph{loops}. If
the vector field is admissible, it defines a \emph{partial} order between
source cells: the relation \(a > a'\) is satisfied between source cells if and
only if a \(V\)-path goes from \(a\) to \(a'\). The non-existence of loops
guarantees this is actually a partial order.

Conversely, let \(V\) be a vector field for our matrix \(M\). If \(1 \leq a, a'
\leq m\), with \(a \neq a'\), we can decide \(a > a'\) if there is an
\emph{elementary} V-path from \(a\) to \(a'\), that is, if a vector \((a,b)\)
is present in \(V\) and the entry \(M[a', b]\) is non-null; for this
corresponds to a cell \(b\) with in particular \(a\) as regular face and \(a'\)
as an arbitrary face. We so obtain a binary relation. Then the vector field
\(V\) is admissible if and only if this binary relation actually transitively
generates a partial order, that is, if again there is no loop \(a_1 > a_2 >
\cdots > a_k = a_1\).

We can therefore summarize the reduction problem by a vector field as follows:
given the matrix \(M\), what process could produce a vector field as large as
possible, but admissible, that is, without any loops? Finding such a vector
field of maximal size seems much too difficult in real applications. Finding a
maximal admissible vector field, not the same problem, is more reasonable but
still serious. We start with some simple heuristic strategies to obtain
significant admissible vector fields, large but most often not maximal.

\subsubsection{Using a predefined order.}

{A direct way to quickly construct an admissible vector field consists in
\emph{pre}defining an order between row indices, and to collect all the indices
for which some column is ``above this index''. Let us play with this toy-matrix
given by our random generator:}\vspace{-10pt}

 {\footnotesize\begin{equation} M =
\left[\begin{array}{ccccc}
  0 & 0 & -1 & -1 & 0
  \\
  0 & -1 & 0 & 0 & 1
  \\
  0 & 0 & 0 & 1 & 1
  \\
  0 & -1 & 1 & 0 & -1
  \\
  -1 & 1 & -1 & 0 & 0
\end{array} \right]
\end{equation}}

If we take simply the index order between row indices, we see the columns 1, 4
and 5 can be selected, giving the vector field \(\{(5,1), (3,4), (4,5)\}\).
This leads to reorder rows and colums in the respective orders \((5,4,3,1,2)\)
and \((1,5,4,2,3)\), rewriting the matrix \(M\) as:\vspace{-10pt}

 {\footnotesize{\begin{equation}
 M = \left[\begin{array}{ccccc}
  -1 & 0 & 0 & 1 & -1
  \\
  0 & -1 & 0 & -1 & 1
  \\
  0 & 1 & 1 & 0 & 0
  \\
  0 & 0 & -1 & 0 & -1
  \\
  0 & 1 & 0 & -1 & 0
 \end{array}\right]
\end{equation}}}%
where the \(3 \times 3\) top left-hand block is triangular. The block
decomposition which follows then corresponds to the one used in the Hexagonal
Lemma~\ref{90055}:\vspace{-10pt}

{\footnotesize\begin{equation}
 \varepsilon = \left[\begin{array}{ccc}
 -1 & 0 & 0 \\ 0 & -1 & 0 \\ 0 & 1 & 1
 \end{array}\right],\
 \varphi = \left[\begin{array}{cc}
 1 & -1 \\ -1 & 1 \\ 0 & 0
 \end{array}\right],\
 \psi = \left[\begin{array}{ccc}
 0 & 0 & -1 \\ 0 & 1 & 0
 \end{array}\right],\
 \beta = \left[\begin{array}{cc}
 0 & -1 \\ -1 & 0
 \end{array}\right].
\end{equation}}

The Hexagonal Lemma produces a reduction of \(M: \bZ^5 \leftarrow \bZ^5\) to
\(M': \bZ^2 \leftarrow \bZ^2\) with \(M' = \beta - \psi \varepsilon^{-1}
\phi\), where the computation of \(\varepsilon^{-1}\) is easy, for
\(\varepsilon\) is triangular unimodular, a computation directly given by the
formula~\ref{08522}. This gives:\vspace{-10pt}

{\footnotesize\begin{equation} M' = \left[\begin{array}{cc}
 -1 & 0 \\ -2 & 1
\end{array}\right]
\end{equation}}%
which is at once triangular unimodular, so that in fact our random matrix \(M\)
is an automorphism of \(\bZ^5\).

\subsubsection{Geometric orders.}

For chain complexes coming from actual digital images, partial orders coming
from geometrical properties can be very convenient. We take again a toy
example, a screen with a \(3 \times 3\) ``resolution'' and this image, eight
pixels black and one white.
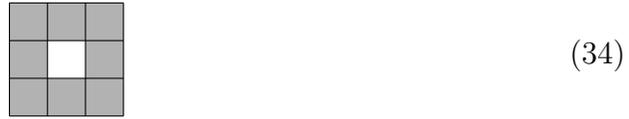
\begin{equation}\begin{tikzpicture} [scale = 0.5, baseline = 0.7cm]
 \fill [black!30] (0,0) rectangle (3,3) ;
 \fill [white] (1,1) rectangle (2,2) ;
 \foreach \i in {0,...,3}
   {\draw (\i, 0) -- +(0,3) ;
    \draw (0, \i) -- +(3,0) ;
   }
\end{tikzpicture}\end{equation}

The bases of the corresponding cellular complex are made of 16 vertices, 24
edges and 8 squares.
\begin{equation}
 0 \leftarrow \bZ^{16} \leftarrow \bZ^{24} \leftarrow \bZ^8 \leftarrow 0
\end{equation}

The boundary matrices are a little complicated. To design an \emph{admissible}
vector field, we can decide the only allowed vectors are oriented leftward or
downward, this is enough to avoid \emph{loops}. Various systematic methods are
possible. Such a method could for example produce this vector field:
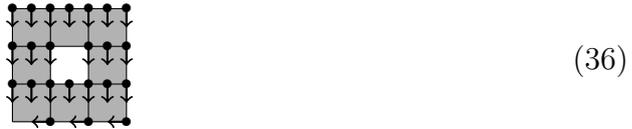
\begin{equation}\begin{tikzpicture} [scale = 0.5, baseline = 0.7cm]
 \fill [black!30] (0,0) rectangle (3,3) ;
 \fill [white] (1,1) rectangle (2,2) ;
 \foreach \i in {0,...,3}
   {\draw (\i, 0) -- +(0,3) ;
    \draw (0, \i) -- +(3,0) ;
   }
 \foreach \i in {0,0.5,...,3}
   \foreach \j in {1,3}
     {\draw [->, thick] (\i,\j) node {\scriptsize\(\bullet\)} -- +(0,-0.5) ;}
 \foreach \i in {0,0.5,1,2,2.5,3}
   {\draw [->, thick] (\i,2) node {\scriptsize\(\bullet\)} -- +(0,-0.5) ;}
 \foreach \i in {1,2,3}
   {\draw [->, thick] (\i,0) node {\scriptsize\(\bullet\)} -- +(-0.5,0) ;}
\end{tikzpicture}\end{equation}

There remains only two critical cells, one vertex and one edge:
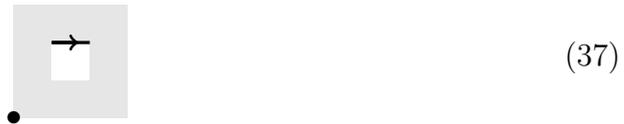
\begin{equation}\begin{tikzpicture} [scale = 0.5, baseline = 0.7cm]
 \fill [black!10] (0,0) rectangle (3,3) ;
 \fill [white] (1,1) rectangle (2,2) ;
 \draw [very thick] (1,2) -- (2,2) ;
 \draw [very thick, ->] (1,2) -- (1.7,2) ;
 \node at (0,0) {\(\bullet\)} ;
\end{tikzpicture}\end{equation}

The reduced chain complex is \(0 \leftarrow \bZ \leftarrow \bZ \leftarrow 0\)
with the null map between both copies of \(\bZ\). The reduction described by
Theorem~\ref{81450} gives in particular a chain map \(g: C^c_\ast \rightarrow
C_\ast\) from the critical chain complex to the initial one; the generator of
\(H_1(C_\ast^c)\) is the only critical edge and its image in the initial chain
complex is this cycle:
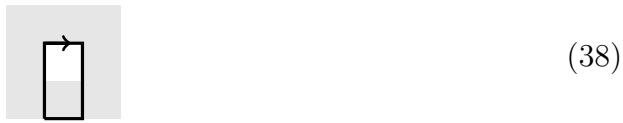
\begin{equation}\begin{tikzpicture} [scale = 0.5, baseline = 0.7cm]
  \fill [black!10] (0,0) rectangle (3,3) ;
  \fill [white] (1,1) rectangle (2,2) ;
  \draw [very thick] (1,0) -- (1,2) -- (2,2) -- (2,0) -- (1,0) ;
 \draw [very thick, ->] (1,2) -- (1.7,2) ;
\end{tikzpicture}
\end{equation}

This reduction informs that \(H_1(C_\ast) = \bZ\) and produces a representant
for the generating homology class. The \(f_1\)-component of the reduction
\(f_1: C_1 \rightarrow C^c_1\) is null except for these two edges where the
image is the unique critical edge:
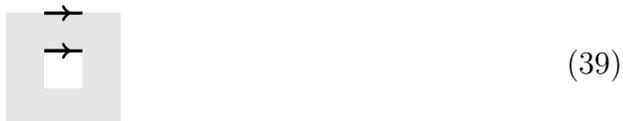
\begin{equation}\begin{tikzpicture} [scale = 0.5, baseline = 0.7cm]
  \fill [black!10] (0,0) rectangle (3,3) ;
  \fill [white] (1,1) rectangle (2,2) ;
  \foreach \i in {2,3} {\draw [very thick] (1,\i) -- (2,\i) ;
                        \draw [very thick, ->] (1,\i) -- (1.7,\i) ;
}
\end{tikzpicture}
\end{equation}
So that if some cycle is given in the initial chain complex:
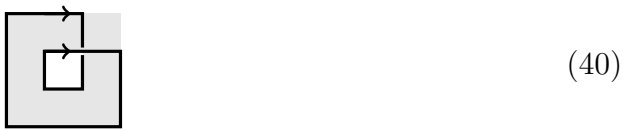
\begin{equation}\begin{tikzpicture} [scale = 0.5, baseline = 0.7cm]
  \fill [black!10] (0,0) rectangle (3,3) ;
  \fill [white] (1,1) rectangle (2,2) ;
  \draw [very thick] (2,2.1) -- (2,3) -- (0,3) -- (0,0) -- (3,0) -- (3,2)
                -- (1,2) -- (1,1) -- (2,1) -- (2,1.9) ;
  \draw [very thick, ->] (1,2) -- (1.7,2) ;
  \draw [very thick, ->] (1,3) -- (1.7,3) ;
\end{tikzpicture}
\end{equation}
the \(f\)-image gives its homology class, here twice the generator of
\(H_1(C_\ast^c)\). If ever the homology class so calculated is null, again the
reduction produces a boundary preimage. For example the homology class of this
rectangle 1-cycle is null and the image shows the boundary preimage, the sum of
two squares, \emph{computed} by the vector field.
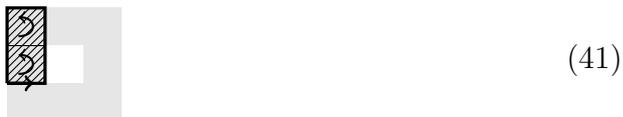
\begin{equation}\begin{tikzpicture} [scale = 0.5, baseline = 0.7cm]
  \fill [black!10] (0,0) rectangle (3,3) ;
  \fill [white] (1,1) rectangle (2,2) ;
  \draw [very thick] (0,1) -- (1,1) -- (1,3) -- (0,3) -- (0,1) ;
 \draw [very thick, ->] (0,1) -- (0.7,1) ;
 \draw [->, thick] (0.3,2.25) .. controls +(0:0.5) and +(0:0.5) .. (0.3,2.75) ;
 \draw [->, thick] (0.3,1.25) .. controls +(0:0.5) and +(0:0.5) .. (0.3,1.75) ;
  \begin{scope}
  \clip (0,1) rectangle (1,3) ;
  \foreach \i in {0,0.2,...,3}
    {\draw (0,\i) -- +(1,1) ;}
  \draw (0,2) -- (1,2) ;
  \end{scope}
\end{tikzpicture}
\end{equation}

All these comments can here be directly read from the vector field, for this
``image'' is very small. But this process can easily be made automatic for the
actual images produced by computers, describing the \emph{effective} homology
of this image, with the same amount of information: the \emph{homological
problem} for this image is \emph{solved}.

\subsubsection{Constructing an appropriate order.} \label{62238}

A more sophisticated strategy consists, given an \emph{admissible} vector field
already constructed, in trying to add a new vector to obtain a better
reduction. The already available vector field defines a partial order between
the source cells with respect to this vector field and the game now is to
search a new vector to be added, but keeping the admissibility property. This
process is applied by  starting from the void vector field.

Let us try to apply this process to the same matrix as before:\vspace{-10pt}

 {\footnotesize\begin{equation} M =
\left[\begin{array}{ccccc}
  0 & 0 & -1 & -1 & 0
  \\
  0 & -1 & 0 & 0 & 1
  \\
  0 & 0 & 0 & 1 & 1
  \\
  0 & -1 & 1 & 0 & -1
  \\
  -1 & 1 & -1 & 0 & 0
\end{array} \right]
\end{equation}}

We start with the void vector field \(V_0 = \{\}\). Running the successive rows
in the usual reading order, we find \(M[1,3] = -1\), and we add the vector
(1,3), obtaining \(V_1 = \{(1,3)\}\). Only one source cell 1, but we must note
that it is from now on forbidden to add a vector which would produce the
relation \(4 > 1\) or \(5 > 1\): this will generate a loop \(1 > 4 > 1 > \cdots
\) and the same for 5. In other words, the partial order to be recorded is \(1
> 4\) and \(1 > 5\), even if 4 and 5 are not yet source cells. Also the row 1
and the column 3 are now used and cannot be used anymore.

We read the row 2 and find \(M[2,2] = -1\), which suggests to add the vector
\((2,2)\), possible, with the same restrictions as before. Now \(V_2 = \{(1,3),
(2,2)\}\).

Reading the row 3 suggests to add the vector \((3,4)\) where 4 has 1 as a face,
because \(M[1,4] = -1\). This does not create any cycle, and we define \(V_3 =
\{(1,3), (2,2), (3,4)\}\). We note also that \(3 > 1\).

Reading the row 4, the only possibility would be the new vector \((4,5)\), but
2 is a face of 5 and this would generate the loop \(4 > 2 > 4 > \cdots\),
forbidden. It is impossible to add a vector \((4,-)\).

Finally we can add the vector \((5,1)\), convenient, for 1 has no other face
than~5; adding this vector certainly keeps the admissibility property.

This leads to the vector field \(V_4 = \{(1,3), (2,2), (3,4), (5,1)\}\), which
generates the partial order on \(1,2,3,5\) where the only non-trivial relations
are \(3 > 1 > 5\) and \(2 > 5\). In particular \(\lambda(3) = 3\), \(\lambda(1)
= \lambda(2) = 2\) and \(\lambda(5) = 1\), see the comments preparing
Proposition~\ref{83464}. Reordering the rows and columns in the respective
orders \((3,1,2,5,4)\) and \((4,3,2,1,5)\) gives the new form for our
matrix:\vspace{-10pt}

{\footnotesize\begin{equation}
 M = \left[\begin{array}{ccccc}
 1 & 0 & 0 & 0 & 1
 \\
 -1 & -1 & 0 & 0 & 0
 \\
 0 & 0 & -1 & 0 & 1
 \\
 0 & -1 & 1 & -1 & 0
 \\
 0 & 1 & -1 & 0 & -1
 \end{array}\right]
\end{equation}}

The vector field has 4 components, and the \(4 \times 4\) top left-hand
submatrix is triangular unimodular. The reduction is better. The corresponding
blocks are:\vspace{-10pt}

{\footnotesize\begin{equation}
 \varepsilon = \left[\begin{array}{cccc}
 1 & 0 & 0 & 0 \\ -1 & -1 & 0 & 0 \\ 0 & 0 & -1 & 0 \\ 0 & -1 & 1 & -1
 \end{array}\right],\
 \varphi = \left[\begin{array}{c}
 1 \\ 0 \\ 1 \\ 0
 \end{array}\right],\
 \psi = \left[\begin{array}{cccc}
 0 & 1 & -1 & 0
 \end{array}\right],\
 \beta = \left[\begin{array}{c}
 -1
 \end{array}\right].
\end{equation}}

The same formula as before gives a reduction producing the matrix
\(\left[-1\right]\) which of course can be reduced to the void matrix.

Both examples show a first step of reduction produces a smaller matrix which in
turn can sometimes be also reduced, even if the used vector field is maximal.

\subsubsection{In more realistic situations.}

The toy example of the previous section in fact is misleading. For such a small
matrix, it is easy to maintain the state of the situation on one's draft sheet,
but when the matrix is for example \(100,\!000 \times 100,\!000\) with 10
non-null entries on every column, the work becomes harder.

Let \(V_k\) be a vector field already obtained, which we intend to enrich by a
new vector \(v\). This vector field generates an order graph, such as:
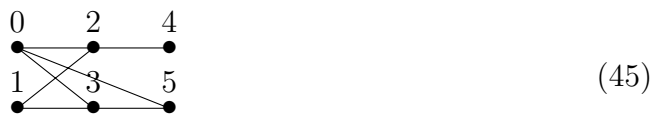
\begin{equation}\begin{tikzpicture} [inner sep = 1.5pt, yscale = 0.8, baseline = 0.3cm]
 \foreach \i in {0,1,2} \foreach \j in {0,1}
   {\coordinate (\i\j) at (\i,\j) ;}
 \foreach \i/\j in {00/1, 01/0, 10/3, 11/2, 20/5, 21/4}
   {\node [label = above:\j] at (\i) {\(\bullet\)} ;}
 \draw (21) -- (01) -- (20) -- (00) -- (11) (10) -- (01) ;
\end{tikzpicture}\end{equation}
to be read as follows: with respect to \(V_k\), the cells 0, 1, 2 and 3 are
source cells, 4 and 5 are not source cells, and \(a > b\) if \(a\) is connected
to \(b\) by a path going rightward. For example \(0 > 4\) and~5 are true, but
\(2 > 5\) is false, which of course in this context does not imply \(2 \leq
5\)!

The cells of such a graph are divided in two parts: the \emph{source} cells, 0,
1, 2 and~3 in our example; the \emph{minimal} cells, here 4 and 5, certainly
not source cells with respect to \(V_k\).

A new vector to be added is something like \(a_0 > \{a_1,a_2,a_3\}\). We decide
adding this vector is allowed if one of these conditions is satisfied:
\begin{fenum}
\item
\(a_0\) is a minimal cell and \(a_1\), \(a_2\) and \(a_3\) are not source cells
of the previous graph, otherwise a loop \emph{could} be generated.
\item
\(a_0\) is not present in our graph order.
\end{fenum}

The new vector to be added cannot have a source already used, but it can be a
minimal cell of the previous graph. With respect to our example, adding the
vector \(4 > \{5,6\}\) is correct, generating the new order graph:
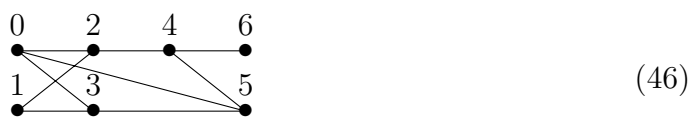
\begin{equation}\begin{tikzpicture} [inner sep = 1.5pt, yscale = 0.8, baseline = 0.3cm]
 \foreach \i in {0,1,2,3} \foreach \j in {0,1}
   {\coordinate (\i\j) at (\i,\j) ;}
 \foreach \i/\j in {00/1, 01/0, 10/3, 11/2, 30/5, 21/4, 31/6}
   {\node [label = above:\j] at (\i) {\(\bullet\)} ;}
 \draw (31) -- (01) -- (30) -- (00) -- (11) (10) -- (01) (21) -- (30) ;
\end{tikzpicture}\end{equation}
The status of 4 is changed from minimal to source, and the new cell 6 gets the
minimal status.

It would be illegal to add the vector \(4 > 1\), because this would generate a
loop. It would be legal to add the vector \(4 > 3\), this would not generate
any loop, but such a possibility is missed by our simplified method, for this
would need a complete analysis of the order relation, too time consuming for
real examples.

Also, the vector \(6 > \{3,4,7\}\) can certainly be added, for its source 6 is
not in the previous graph.
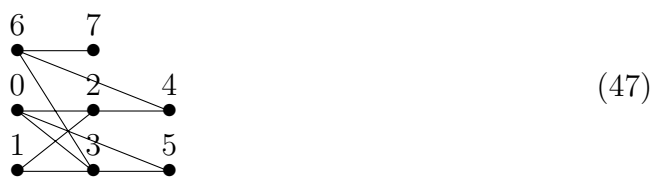
\begin{equation}\begin{tikzpicture} [inner sep = 1.5pt, yscale = 0.8, baseline = 1cm]
 \foreach \i in {0,1,2} \foreach \j in {0,1,2}
   {\coordinate (\i\j) at (\i,\j) ;}
 \foreach \i/\j in {00/1, 01/0, 10/3, 11/2, 20/5, 21/4, 02/6, 12/7}
   {\node [label = above:\j] at (\i) {\(\bullet\)} ;}
 \draw (21) -- (01) -- (20) -- (00) -- (11) (10) -- (01)
       (10) -- (02) -- (21) (02) -- (12) ;
\end{tikzpicture}\end{equation}

If the reader would like a toy example illustrating this example, it is enough
to reexamine the example of Section~\ref{62238}. The vector field which was
obtained by a careful examination of the situation is in fact also produced by
the above automatic heuristic method! The corresponding order graph for this
example is this one:
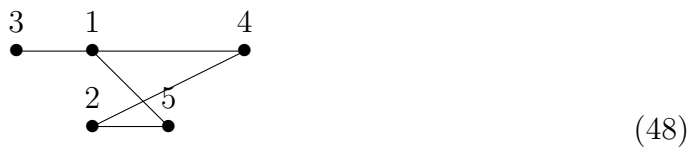
\begin{equation}\begin{tikzpicture} [inner sep = 2pt]
 \foreach \i in {0,1,2,3} \foreach \j in {0,1}
   {\coordinate (\i\j) at (\i,\j) ;}
 \foreach \i/\j in {01/3, 10/2, 11/1, 20/5, 31/4}
   {\node [label = above:\j] at (\i) {\(\bullet\)} ;}
 \draw (01) -- (31) -- (10) -- (20) -- (11) ;
\end{tikzpicture}
\end{equation}
where 4 is the unique remaining minimal cell.

\subsubsection{The corresponding graph problem.}

To finish this short study of the chain complexes connected to images, let
\(M\) be an integer \(r \times c\) matrix. We consider firstly the particular
case where all the entries are elements of \(\{-1,0,+1\}\). Such a matrix
generates a graph \(G\) with nodes of two colors: the row nodes, indexed by
\([1\ldots r]\) and the column nodes indexed by \([1\ldots c]\). Every matrix
entry \(\pm 1\) generates an edge connecting the corresponding row and column
nodes. The coloring of the graph is compatible with the incidence relations.
The game is then the following: construct a collection of edges \(V =
\{(r_i,c_i)\}\) satisfying the following requirements. Every first component
\(r_i\) is a row index, every second component \(c_i\) is a column index. The
\(r_i\)'s (resp. the \(c_i\)'s) are pairwise different. Consider the subgraph
made of all these \(c_i\)'s and their \emph{neighbouring} row nodes. Orient the
edge \((r_i,c_i) \in V\) from the row node to the column node for every
selected edge, \emph{and} for all the \emph{other} edges starting from \(c_i\),
orient these edges from \(c_i\) to the row node. Then you must not have any
loop in this graph.

Let us take again our toy matrix:\vspace{-10pt}

 {\footnotesize\begin{equation} M =
\left[\begin{array}{ccccc}
  0 & 0 & -1 & -1 & 0
  \\
  0 & -1 & 0 & 0 & 1
  \\
  0 & 0 & 0 & 1 & 1
  \\
  0 & -1 & 1 & 0 & -1
  \\
  -1 & 1 & -1 & 0 & 0
\end{array} \right]
\end{equation}}

The graph corresponding to the vector field constructed in Section~\ref{62238}
is:
\begin{equation}\begin{tikzpicture} [draw, inner sep = 1pt, scale = 0.5, baseline=1cm]
 \foreach \i in {-1,0,...,5} \foreach \j in {0,...,4}
   {\coordinate (\i\j) at (\i,\j) ;}
 \foreach \i/\j in {-14/1, 12/3, 30/4, 34/2, 52/5}
   {\node [draw, inner sep = 2pt] (c\j) at (\i) {\j} ;}
 \foreach \i/\j in {14/5, 10/1, 32/4, 50/3, 54/2}
   {\node [draw, circle] (r\j) at (\i) {\j} ;}
 \foreach \i/\j in {5/1, 2/2, 1/3, 3/4}
    {\draw [->, very thick] (r\i) -- (c\j) ;}
 \foreach \i/\j in {2/4, 2/5, 3/4, 3/5, 4/1}
    {\draw [->] (c\i) -- (r\j) ;}
 \foreach \j in {2,3,4} {\draw [dotted] (c5) -- (r\j) ;}
\end{tikzpicture}
\end{equation}
where the row nodes are circles and the column nodes are squares. The thick
arrows are the elements of our vector field; the thin arrows go from one column
to the neighbouring row nodes, when the column has been selected, except the
component of the vector field itself, oriented in the reverse direction.
Finally the column node 5 is not used.

It is then not difficult to prove that in our example there \emph{does not
exist} any admissible vector field of size 5. It would be necessary to ``use''
every node. You cannot select the vector \((4,2)\), for the last possibility
for the row node 2 would be \((2,5)\) generating a loop:
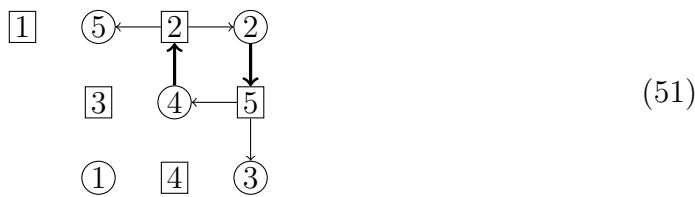
\begin{equation}\begin{tikzpicture} [draw, inner sep = 1pt, scale = 0.5, baseline=1cm]
 \foreach \i in {-1,0,...,5} \foreach \j in {0,...,4}
   {\coordinate (\i\j) at (\i,\j) ;}
 \foreach \i/\j in {-14/1, 12/3, 30/4, 34/2, 52/5}
   {\node [draw, inner sep = 2pt] (c\j) at (\i) {\j} ;}
 \foreach \i/\j in {14/5, 10/1, 32/4, 50/3, 54/2}
   {\node [draw, circle] (r\j) at (\i) {\j} ;}
 \foreach \i/\j in {4/2, 2/5}
    {\draw [->, very thick] (r\i) -- (c\j) ;}
 \foreach \i/\j in {2/2, 5/4, 2/5, 5/3}
    {\draw [->] (c\i) -- (r\j) ;}
\end{tikzpicture}
\end{equation}

In the same way, you cannot use the vector \((4,5)\). Necessarily, the vector
\((4,3)\) must be used. But it would be necessary to use the vectors \((1,4)\)
and \((3,5)\), again generating a loop of period 3. In other words, no reorder
of rows and columns can make our matrix triangular. Which does not prevent our
random matrix from being an isomorphism.

The vector field of size 4 is in fact in this example the vector field of
maximal size. Solving such a problem for giant graphs of this sort is probably
intractable.

For matrices with arbitrary integer entries, not necessarily in
\(\{-1,0,+1\}\), the only difference is the following. If an entry \(a_{r,c}\)
of the matrix has an absolute value \(> 1\), then the vector \((r,c)\) is
forbidden. If the column \(c\) is used in a vector field, the orientation of
the edge \((r,c)\) is necessarily \([c \rightarrow r]\).

\section{The Eilenberg-Zilber W-reduction.} \label{25115}

\subsection{Introduction.}

The Eilenberg-Zilber theorem is essential in combinatorial topology, usually
presented as a result of algebraic topology. The systematic use of
\emph{discrete vector fields} gives a more precise analysis: this result is in
fact a particular case of deformations \`a la Whitehead, which gives as a
\emph{by-product} the usual homological result.

The subject is the following. Let \(\Delta^p\) and \(\Delta^q\) be two standard
simplices of respective dimensions \(p\) and \(q\). What about the product
\(\Delta^{p,q} = \Delta^p \times \Delta^q\)? The lazy solution consists in
enriching your collection of \emph{elementary} models by all the possible
products of simplices, so that \(\Delta^{p,q}\) is then viewed as
``elementary''. Why not, but the penalty is not far: the numerous results
patiently obtained in \emph{simplicial} topology are nomore valid for the
spaces made of these less elementary models.

Another solution consists in \emph{triangulating} this product
\(\Delta^{p,q}\), feasible but not so easy. With a drawback, it is then
difficult to read the product structure in this triangulation.

The Eilenberg-Zilber theorem settles the right connection between both
solutions, certainly one of the most important results in Algebraic Topology:
it is in fact the heart of the fundamental Serre and Eilenberg-Moore spectral
sequences, and our analysis based on discrete vector fields gives a precise
description of this interpretation.

\subsection{Triangulations.}

We have to work in the simplicial complex \(\Delta^{p,q} = \Delta^p \times
\Delta^q\). A vertex of \(\Delta^p\) is an integer in \(\up = [0\ldots p]\), a
(non-degenerate) \(d\)-simplex of \(\Delta^p\) is a strictly increasing
sequence of integers \(0 \leq v_0 < \ldots < v_d \leq p\). The same for our
second factor \(\Delta^q\).

The canonical triangulation of \(\Delta^p \times \Delta^q\) is made of
(non-degenerate) simplices \(((v_0,v'_0), \ldots, (v_d,v'_d))\) satisfying the
relations:
\begin{itemize}\itemsep = -3pt
\item
\(0 \leq v_0 \leq v_1 \leq \cdots \leq v_d \leq p\).
\item
\(0 \leq v'_0 \leq v'_1 \leq \cdots \leq v'_d \leq q\).
\item
\((v_i,v'_i) \neq (v_{i-1},v'_{i-1})\) for \(1 \leq i \leq d\).
\end{itemize}
In other words, the canonical triangulation of \(\Delta^{p,q} = \Delta^p \times
\Delta^q\) is associated to the poset \(\up \times \uq\) endowed with the
\emph{product order} of the factors. For example the three maximal simplices of
\(\Delta^{2,1} = \Delta^2 \times \Delta^1\) are:
\begin{itemize}
\item
\(((0,0), (0,1), (1,1), (2,1))\).
\item
\(((0,0), (1,0), (1,1), (2,1))\).
\begin{picture}(0,0)(0,0)
\makebox{\begin{tikzpicture} [baseline = 1cm, scale = 0.6] \node at (-5,0) {} ;
\draw [thick] (0,0) -- (1,-1) -- (3,0) -- (3,4) -- (1,3) -- (1,-1) ; \draw
[thick] (0,0) -- (0,4) -- (3,4) -- (1,3) -- (0,4) ; \draw [dashed] (0,0) --
(1,3) (1,-1) -- (3,4) ; \draw [dotted] (0,0) -- (3,0) (0,0) -- (3,4) ;
\begin{scope} [font = \scriptsize]
\node [left] at (0,0) {(0,0)} ; \node [below] at (1,-1) {(1,0)} ; \node [right]
at (3,0) {(2,0)} ; \node [left] at (0,4) {(0,1)} ; \node [anchor = 160] at
(1,3) {(1,1)} ; \node [right] at (3,4) {(2,1)} ;
\end{scope}
\end{tikzpicture}}
\end{picture}
\item
\(((0,0), (1,0), (2,0), (2,1))\).
\end{itemize}

\subsection{Simplex = s-path.}

We can see the poset \(\up \times \uq\) as a lattice where we arrange the first
factor \(\up\) in the horizontal direction and the second factor \(\uq\) in the
vertical direction. The first figure below is the lattice \(\underline{2}
\times \underline{1}\) while the other figures are representations of the
maximal simplices of \(\Delta^{2,1} = \Delta^2 \times \Delta^1\) as
\emph{increasing} paths in the lattice.
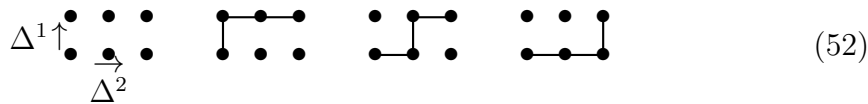
\begin{equation}\begin{tikzpicture} [scale=0.5, baseline = 0cm]
\foreach \i in {0,1,2} {\foreach \j in {0,1} {\node (\i\j) at (\i,\j)
{\(\bullet\)} ;} ;} ; \node (ua) at (-0.3,0.5) {\(\uparrow\)} ; \node [left] at
(ua) {\(\Delta^1\)} ; \node (ra) at (1,-0.3) {\(\rightarrow\)} ; \node [below]
at (ra) {\(\Delta^2\)} ;
\begin{scope} [xshift = 4cm]
\foreach \i in {0,1,2} {\foreach \j in {0,1} {\node (\i\j) at (\i,\j)
{\(\bullet\)} ;} ;} ; \draw [thick] (0,0) -- (0,1) -- (1,1) -- (2,1) ;
\end{scope}
\begin{scope} [xshift = 8cm]
\foreach \i in {0,1,2} {\foreach \j in {0,1} {\node (\i\j) at (\i,\j)
{\(\bullet\)} ;} ;} ; \draw [thick] (0,0) -- (1,0) -- (1,1) -- (2,1) ;
\end{scope}
\begin{scope} [xshift = 12cm]
\foreach \i in {0,1,2} {\foreach \j in {0,1} {\node (\i\j) at (\i,\j)
{\(\bullet\)} ;} ;} ; \draw [thick] (0,0) -- (1,0) -- (2,0) -- (2,1) ;
\end{scope}
\end{tikzpicture}\end{equation}

\begin{dfn} ---
\emph{An \emph{s-path} \(\pi\) of the lattice \(\up \times \uq\) is a finite
sequence \(\pi = ((a_i,b_i))_{0 \leq i \leq d}\) of elements of \(\up \times
\uq\) satisfying \( (a_{i-1},b_{i-1}) < (a_i,b_i)\) for every \(1 \leq i \leq
d\) with respect to the product order. The \(d\)-simplex \(\sigma_\pi\)
represented by the path \(\pi\) is the convex hull of the points \((a_i,b_i)\)
in the prism \(\Delta^{p,q}\).}\qed
\end{dfn}

The simplices \(\Delta^p\) and \(\Delta^q\) have affine structures which define
a product affine structure on \(\Delta^{p,q}\), and the notion of convex hull
is well defined on \(\Delta^{p,q}\).

``S-path'' stands for ``path representing a simplex'', more precisely a
non-degenerate simplex. Replacing the strict inequality between two successive
vertices by a non-strict inequality would lead to analogous representations for
degenerate simplices, but such simplices are not to be considered in this
section.

This representation of a simplex as an s-path running in a lattice is the key
point to master the relatively complex structure of the canonical prism
triangulations.

\begin{dfn}\label{07444} ---
 \emph{The \emph{last simplex} \(\lambda_{p,q}\) of the prism \(\Delta^{p,q}\) is the
 \((p+q)\)-simplex defined by the path:
 \begin{equation}
 \lambda_{p,q} = ((0,0), (1,0), \ldots, (p,0), (p,1), \ldots, (p,q)).
 \end{equation}}
\end{dfn}\vspace{-15pt}
\qed

The path runs some edges of \(\Delta^p \times 0\), visiting all the
corresponding vertices in the right order; next it runs some edges of \(p
\times \Delta^q\), visiting all the corresponding vertices also in the right
order. Geometrically, the last simplex is the convex hull of the visited
vertices. The last simplex of the prism \(P_{1,2} = \Delta^1 \times \Delta^2\)
is shown in the figure below. The path generating the last simplex is drawn in
full lines, the other edges of this last simplex are dashed lines, and the
other edges of the prism are in dotted lines.

\begin{equation}\begin{tikzpicture} [baseline = 0.5cm]
 \begin{scope} [xshift = -1cm]
    \node at (-1,1) {\(\Delta^2\)} ;
    \coordinate (0) at (0,0) ;
    \coordinate (1) at (-0.5,0.9) ;
    \coordinate (2) at (0.3,1.5) ;
    \draw [thick, fill = black!60] (1.center) -- (0.center) -- (2.center) -- (1.center) ;
    \node at (0,0) {\(\bullet\)} ;
    \node at (-0.5,0.9) {\(\bullet\)} ;
    \node at (0.3,1.5) {\(\bullet\)} ;
 \end{scope}
 \coordinate (00) at (0,0) ;
 \coordinate (01) at (-0.5,0.9) ;
 \coordinate (02) at (0.3,1.5) ;
 \coordinate (10) at (3,0) ;
 \coordinate (11) at (2.5,0.9) ;
 \coordinate (12) at (3.3,1.5) ;
 \foreach \i in {00,01,02,10,11,12}
   {\node at (\i) {\(\bullet\)} ; } ;
 \draw [thick] (00) -- (10) -- (11) -- (12) ;
 \draw [dashed] (11) -- (00) -- (12) -- (10) ;
 \draw [dotted, thin] (01) -- (00) -- (02) -- (01) -- (11)
                      (01) -- (12) -- (02) ;
 \draw [thick] (0,-0.7) node {\(\bullet\)} -- node [above] {\(\Delta^1\)} (3,-0.7)
    node {\(\bullet\)} ;
\end{tikzpicture}\end{equation}

\subsection{Subcomplexes.}

\begin{dfn} ---
\emph{The \emph{hollowed prism} \(H\Delta^{p,q} \subset \Delta^{p,q}\) is the
difference:
\begin{equation}
H\Delta^{p,q} := \Delta^{p,q} - \textrm{int(last simplex).}
\end{equation} }
\end{dfn}\vspace{-15pt}
\qed

The faces of the last simplex are retained, but the interior of this simplex is
removed.

\begin{dfn} ---
\emph{The \emph{boundary} \(\partial \Delta^{p,q}\) of the prism
\(\Delta^{p,q}\) is defined by:
\begin{equation}
 \partial \Delta^{p,q} := (\partial\Delta^p \times \Delta^q) \cup
 (\Delta^p \times \partial\Delta^q)
\end{equation}}
\end{dfn}\vspace{-15pt}
\qed

It is the geometrical Leibniz formula.

We will give a detailed description of the pair \((H\Delta^{p,q}, \partial
\Delta^{p,q})\) as a \emph{W-contraction}, cf. Definition~\ref{41027}; it is a
combinatorial version of the well-known topological contractibility of
\(\Delta^{p,q} - \{\ast\}\) on \(\partial \Delta^{p,q} \) for every point
\(\ast\) of the interior of the prism. A very simple admissible vector field
will be given to homologically annihilate the difference \(H\Delta^{p,q} -
\partial\Delta^{p,q}\). In fact, carefully ordering the components of this
vector field will give the desired W-contraction.

\subsection{Interior and exterior simplices of a prism.}

\begin{dfn} ---
\emph{A simplex \(\sigma\) of the prism \(\Delta^{p,q}\) is said
\emph{exterior} if it is included in the boundary of the prism: \(\sigma
\subset \partial \Delta^{p,q}\). Otherwise the simplex is said \emph{interior}.
We use the same terminology for the s-paths, implicitly referring to the
simplices coded by these paths.}\qed
\end{dfn}

The faces of an exterior simplex are also exterior, but an interior simplex can
have faces of both sorts.

\begin{prp} ---
An s-path \(\pi\) in \(\up \times \uq\) is interior if and only if the
projection-paths \(\pi_1\) on \(\up\) and \(\pi_2\) on \(\uq\) run all the
respective vertices of \(\up\) and \(\uq\).
\end{prp}

The first s-path \(\pi\) in the figure below represents a 1-simplex \emph{in}
\(\partial P_{1,2}\), for the point~1 is missing in the projection \(\pi_2\) on
the second factor \underline{2}: \(\pi\) is an \emph{exterior} simplex. The
second s-path \(\pi'\) represents an \emph{interior} 2-simplex of \(P_{1,2}\),
for both projections are surjective.

\begin{equation}\label{29457}\begin{tikzpicture} [scale = 0.5, baseline = 0.5cm]
 \foreach \i in {0,1,2}
   {\foreach \j in {0,1}
     {\node (\i\j) at (\j,\i) {\(\bullet\)} ;}}
 \draw [thick] (00.center) -- (21.center) ;
 \node [left] at (-0.3,1) {\(\pi =\)} ;
 \node at (4.2,0) {\(\pi = \partial_1 \pi'\)} ;
\begin{scope} [xshift = 9cm]
 \foreach \i in {0,1,2}
   {\foreach \j in {0,1}
     {\node (\i\j) at (\j,\i) {\(\bullet\)} ;}}
 \draw [thick] (00.center) -- (10.center) -- (21.center) ;
 \node [left] at (-0.3,1) {\(\pi' =\)} ;
\end{scope}\end{tikzpicture}
\end{equation}

In particular, if \(\pi = ((a_i,b_i))_{0 \leq i \leq d}\) is an interior
simplex of \(\Delta^{p,q}\), then necessarily \((a_0,b_0) = (0,0)\) and
\((a_d,b_d) = (p,q)\): an s-path representing an interior simplex of
\(\Delta^{p,q}\) starts from \((0,0)\) and arrives at \((p,q)\).

\prof If for example the first projection of \(\pi\) is not surjective, this
means the first projection of the generating path does not run all the vertices
of \(\Delta^p\), and therefore is included in one of the faces \(\partial_k
\Delta^p\) of \(\Delta^p\). This implies the simplex \(\sigma_\pi\) is included
in \(\partial_k \Delta_p \times \Delta_q \subset \partial \Delta^{p,q}\).\qed

We so obtain a simple descrition of an interior simplex \(((a_i,b_i))_{0 \leq i
\leq d}\): it starts from \((a_0,b_0) = (0,0)\) and arrives at \((a_d,b_d) =
(p,q)\); furthermore, for every \mbox{\(1 \leq i \leq d\)}, the difference
\((a_i,b_i) - (a_{i-1}, b_{i-1})\) is \((0,1)\) or \((1,0)\) or \((1,1)\): both
components of this difference are non-negative, and if one of these components
is \(\geq 2\), then the surjectivity property is not satisfied. In a
geometrical way, the only possible \emph{elementary steps} for an s-path
\(\pi\) describing an \emph{interior} simplex of \(\Delta^{p,q}\) are:
\begin{equation}\begin{tikzpicture} [scale = 0.5, thick, baseline = 0cm]
 \foreach \i in {0,1} {\foreach \j in {0,1}
   {\node (\i\j) at (\j,\i) {\(\bullet\)} ;}}
 \draw (00.center) -- (01.center) ;
 \draw [dotted] (00.center) -- +(-0.5,-0.5)
                (00.center) -- +(-0.7,0)
                (00.center) -- +(0,-0.7) ;
 \draw [dotted] (01.center) -- +(0.5,0.5)
                (01.center) -- +(0.7,0)
                (01.center) -- +(0,0.7) ;
 \begin{scope} [xshift = 6cm]
  \foreach \i in {0,1} {\foreach \j in {0,1}
   {\node (\i\j) at (\j,\i) {\(\bullet\)} ;}}
 \draw (00.center) -- (11.center) ;
 \draw [dotted] (00.center) -- +(-0.5,-0.5)
                (00.center) -- +(-0.7,0)
                (00.center) -- +(0,-0.7) ;
 \draw [dotted] (11.center) -- +(0.5,0.5)
                (11.center) -- +(0.7,0)
                (11.center) -- +(0,0.7) ;
 \end{scope}
 \begin{scope} [xshift = 12cm]
  \foreach \i in {0,1} {\foreach \j in {0,1}
   {\node (\i\j) at (\j,\i) {\(\bullet\)} ;}}
 \draw (00.center) -- (10.center) ;
 \draw [dotted] (00.center) -- +(-0.5,-0.5)
                (00.center) -- +(-0.7,0)
                (00.center) -- +(0,-0.7) ;
 \draw [dotted] (10.center) -- +(0.5,0.5)
                (10.center) -- +(0.7,0)
                (10.center) -- +(0,0.7) ;
 \end{scope}
\end{tikzpicture}
\end{equation}

\subsection{Faces of s-paths.}

If \(\pi = ((a_i,b_i))_{0 \leq i \leq d}\) represents a \(d\)-simplex
\(\sigma_\pi\) of \(\Delta^{p,q}\), the face \(\partial_k \sigma_\pi\) is
represented by the same s-path except the \(k\)-th component \((a_k, b_k)\)
which is removed: we could say this point of \(\up \times \uq\) is
\emph{skipped}. For example in Figure~(\ref{29457}) above, \(\partial_1 \pi' =
\pi\). In particular a face of an interior simplex is not necessarily interior.

\begin{prp} --- \label{75744}
Let \(\pi = ((a_i, b_i))_{0 \leq i \leq d}\) be an s-path representing an
interior \(d\)-simplex of \(\Delta^{p,q}\). The faces \(\partial_0 \pi\) and
\(\partial_d \pi\) are certainly \emph{not} interior. For \(1 \leq k \leq
d-1\), the face \(\partial_k \pi\) is interior if and only if the point
\((a_k,b_k)\) is a right-angle bend of the s-path \(\pi\) in the lattice \(\up
\times \uq\).
\end{prp}

\prof Removing the vertex \((a_0,b_0) = (0,0)\) certainly makes non-surjective
a projection \(\pi_1\) or \(\pi_2\) (or both if \((a_1,b_1) = (1,1)\)). The
same if the last point \((a_d,b_d)\) is removed.

If we examine now the case of \(\partial_k\pi\) for \(1 \leq k \leq d-1\), nine
possible configurations for two consecutive elementary steps before and after
the vertex \((a_k, b_k)\) to be removed:
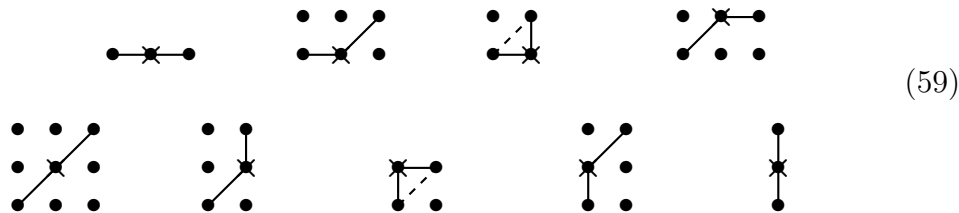
\begin{equation}\begin{tikzpicture} [scale = 0.5, thick, baseline = -0.5cm]
 \foreach \i in {0,1,2} {\node (\i) at (\i,0) {\(\bullet\)} ;} ;
 \draw (0.center) -- (1.center) -- (2.center) ;
 \node at (1,0) {{\boldmath\(\times\)}} ;
\begin{scope} [xshift = 5cm]
 \foreach \i in {0,1}
   {\foreach \j in {0,1,2}
     {\node (\i\j) at (\j,\i) {\(\bullet\)} ;} ;} ;
 \draw (00.center) -- (01.center) -- (12.center) ;
 \node at (1,0) {{\boldmath\(\times\)}} ;
\end{scope}
\begin{scope} [xshift = 10cm]
 \foreach \i in {0,1}
   {\foreach \j in {0,1}
     {\node (\i\j) at (\j,\i) {\(\bullet\)} ;} ;} ;
 \draw (00.center) -- (01.center) -- (11.center) ;
 \draw [dashed] (00.center) -- (11.center) ;
 \node at (1,0) {{\boldmath\(\times\)}} ;
\end{scope}
\begin{scope} [xshift = 15cm]
 \foreach \i in {0,1}
   {\foreach \j in {0,1,2}
     {\node (\i\j) at (\j,\i) {\(\bullet\)} ;} ;} ;
 \draw (00.center) -- (11.center) -- (12.center) ;
 \node at (1,1) {{\boldmath\(\times\)}} ;
\end{scope}
\begin{scope} [xshift = -2.5cm, yshift = -4cm]
 \foreach \i in {0,1,2}
   {\foreach \j in {0,1,2}
     {\node (\i\j) at (\j,\i) {\(\bullet\)} ;} ;} ;
 \draw (00.center) -- (11.center) -- (22.center) ;
 \node at (1,1) {{\boldmath\(\times\)}} ;
\end{scope}
\begin{scope} [xshift = 2.5cm, yshift = -4cm]
 \foreach \i in {0,1,2}
   {\foreach \j in {0,1}
     {\node (\i\j) at (\j,\i) {\(\bullet\)} ;} ;} ;
 \draw (00.center) -- (11.center) -- (21.center) ;
 \node at (1,1) {{\boldmath\(\times\)}} ;
\end{scope}
\begin{scope} [xshift = 7.5cm, yshift = -4cm]
 \foreach \i in {0,1}
   {\foreach \j in {0,1}
     {\node (\i\j) at (\j,\i) {\(\bullet\)} ;} ;} ;
 \draw (00.center) -- (10.center) -- (11.center) ;
 \draw [dashed] (0,0) -- (1,1) ;
 \node at (0,1) {{\boldmath\(\times\)}} ;
\end{scope}
\begin{scope} [xshift = 12.5cm, yshift = -4cm]
 \foreach \i in {0,1,2}
   {\foreach \j in {0,1}
     {\node (\i\j) at (\j,\i) {\(\bullet\)} ;} ;} ;
 \draw (00.center) -- (10.center) -- (21.center) ;
 \node at (0,1) {{\boldmath\(\times\)}} ;
\end{scope}
\begin{scope} [xshift = 17.5cm, yshift = -4cm]
 \foreach \i in {0,1,2} {\node (\i) at (0,\i) {\(\bullet\)} ;} ;
 \draw (0.center) -- (1.center) -- (2.center) ;
 \node at (0,1) {{\boldmath\(\times\)}} ;
\end{scope}
\end{tikzpicture}\end{equation}
In these figures, the intermediate point \(\bullet\hspace{-7.7pt}\times\) of
the displayed part of the considered s-path is assumed to be the point
\((a_k,b_k)\) of the lattice, to be removed to obtain the face
\(\partial_k\pi\). In the cases 1, 2, 4 and 5, skipping this point makes
non-surjective the first projection \(\pi_1\) on \(\up\ \mbox{`='}\ \Delta^p\).
In the cases 5, 6, 8 and 9, the second projection \(\pi_2\) on \(\uq\) becomes
non-surjective. There remain the cases 3 and 7 where the announced right-angle
bend is observed. \qed

\subsection{From vector fields to W-contractions.}

Let \((X,A)\) be an elementary W-contraction, cf Definition~\ref{41027}. The
difference \(X - A\) is made of two non-degenerate simplices \(\sigma\) and
\(\tau\), the first one being a face of the second one with a \emph{unique}
face index. The pair \((\sigma, \tau)\) is nothing but the unique vector of a
vector field \(V\), a vector field which, via the Vector-Field Reduction
Theorem~\ref{81450}, defines also the W-reduction \(C_\ast^{ND} X \rrdc
C_\ast^{ND} A\)

\begin{dfn} ---
\emph{A simplicial pair (X,A) is an elementary \emph{filling} if the difference
\(X-A\) is made of a \emph{unique} non-degenerate simplex \(\sigma\), all the
faces of which are therefore simplices of \(A\).}\qed
\end{dfn}

You might think \(\partial\sigma\) is the initial state of a decayed tooth in
the body \(A\), to be restored by adding \(\textrm{int}(\sigma)\), obtaining
\(X\).

\begin{dfn} ---
\emph{Let \((X,A)\) be a simplicial pair. A \emph{description by a filling
sequence} of this pair, more precisely of the difference \(X-A\), is an
ordering \((\sigma_i)_{0 < i \leq r}\) of the non-degenerate simplices of
\(X-A\) satisfying the following condition: if \(A_i = A \cup (\cup_{j=1}^i
\sigma_j)\), then every pair \((A_i, A_{i-1})\) is an elementary filling.}\qed
\end{dfn}

Every pair \((X,A)\) with a finite number of non-degenerate simplices in
\(X-A\) can be described by a filling sequence: order the missing simplices
according to their dimension. In particular, adding an extra vertex is a
particular filling.

It is convenient to describe the general W-contractions, see
Definition~\ref{11769}, by \emph{special} filling sequences.

\begin{prp} ---
Let (X,A) be a simplicial pair. This pair is a \emph{W-contraction} if and only
if it admits a description by a  filling sequence \(F = (\sigma_i)_{0 < i \leq
2r}\) satisfying the extra condition: for every \emph{even} index \(2i\), the
simplex \(\sigma_{2i-1}\) is a face of \(\sigma_{2i}\) with a unique face
index.
\end{prp}

\prof Such a description is nothing but the vector field \(V =
\{(\sigma_{2i-1}, \sigma_{2i})_{0 < i \leq r}\}\) with an \emph{extra
information}: the vectors are \emph{ordered} in such a way they justify also
the W-contraction property. Such a vector field is necessarily
\emph{admissible}: all the \(V\)-paths go to \(A\) and cannot loop.\qed

This extra information given by the order on the elements of the vector field
is an avatar of the traditional difference between homotopy and homology.

\subsection{The theorem of the hollowed prism.} \label{68696}

\begin{thr} --- \label{79523}
The pair: \((H\Delta^{p,q}, \partial \Delta^{p,q})\) is a W-contraction.
\end{thr}

The hollowed prism can be W-contracted on the boundary of the same prism.

\noindent\(\clubsuit\) The proof is recursive with respect to the pair
\((p,q)\). If \(p = 0\), the boundary of \(\Delta^0 = \ast\) is void, so that
the boundary of \(\Delta^{0,q}\) is simply \(\partial\Delta^q\); the last
simplex is the unique \(q\)-simplex, the hollowed prism \(H\Delta^{0,q}\) is
also \(\partial\Delta^q\): the desired W-contraction is trivial, more precisely
the corresponding vector field is empty. The same if \(q = 0\) for the pair
\((H\Delta^{p,0}, \partial\Delta^{p,0})\).

Now we prove the general case \((p,q)\) with \(p,q>0\), assuming the proofs of
the cases \((p-1,q-1)\), \((p,q-1)\) and \((p-1,q)\) are available. Three
justifying filling sequences are available; it is more convenient to see the
sequences of simplices as sequences of \emph{s-paths}:
\begin{itemize}\itemsep = -3pt
\item
\(F_1 = (\pi^1_i)_{0 < i \leq 2r_1}\) for \(\Delta^{p-1,q-1} = \partial_0
\Delta^p \times \partial_0 \Delta^q\).
\item
\(F_2 = (\pi^2_i)_{0 < i \leq 2r_2}\) for \(\Delta^{p,q-1} = \Delta^p \times
\partial_0 \Delta^q\).
\item
\(F_3 = (\pi^3_i)_{0 < i \leq 2r_3}\) for \(\Delta^{p-1,q} = \partial_0
\Delta^p \times \Delta^q\).
\end{itemize}

All the components of these filling sequences can be viewed as s-paths starting
from \((1,1)\) (resp. \((0,1)\), \((1,0)\)), going to \((p,q)\).

These filling sequences are made of all the non-degenerate s-paths (simplices)
of the difference \(H\Delta^{\ast,\ast} - \partial \Delta^{\ast,\ast}\),
ordered in such a way every face of an s-path is either interior \emph{and}
present \emph{beforehand} in the list, or exterior; furthermore, for the
s-paths of even index, the previous one is one of its faces. Using these
sequences, we must construct an analogous sequence for the bidimension
\((p,q)\).

Every s-path \(\pi^j_i\) of dimension \(d\) can be completed into an interior
s-path \(\overline{\pi}^j_ i\) of dimension \(d+1\) in \(\up \times \uq\) in a
unique way, adding a first diagonal step \(((0,0),(1,1))\) if \(j=1\), or a
first vertical step \(((0,0),(0,1))\) if \(j = 2\), or a first horizontal step
\(((0,0),(1,0))\) if \(j=3\). Conversely, every interior s-path of
\(\Delta^{p,q}\) can be obtained from an interior s-path of
\(\Delta^{p-1,q-1}\), \(\Delta^{p,q-1}\) or \(\Delta^{p-1,q}\) in a unique way
by this completion process.

For example, in the next figure, we illustrate how an s-path \(\pi^1_ i\) of
\(\underline{3} \times \underline{2}\) can be completed into an s-path
\(\overline{\pi}^1_i\) of \(\underline{4} \times \underline{3}\):
\begin{equation}\begin{tikzpicture} [scale = 0.3, font = \scriptsize, baseline = 0.3cm]
 \begin{scope}
 \foreach \i in {0,1,2,3} {\foreach \j in {0,...,3,4} {\node (\i\j) at (\j,\i)
  {\(\bullet\)} ;} ;} ;
 \node [font = \normalsize, left] at (10) {\(\pi^1_ i =\hspace{5pt}\)} ;
 \draw [thick] (11.center) -- (12.center) -- (22.center) -- (33.center)
           -- (34.center) ;
 \end{scope}
\begin{scope} [xshift = 15cm]
 \foreach \i in {0,...,3} {\foreach \j in {0,...,4} {\node (\i\j) at (\j,\i)
  {\(\bullet\)} ;} ;} ;
 \node [font = \normalsize, left] at (10) {\(\overline{\pi}^1_ i =\hspace{5pt}\)} ;
 \draw [thick] (00.center) -- (11.center) -- (12.center) -- (22.center) -- (33.center)
           -- (34.center) ;
\end{scope}
\end{tikzpicture}\end{equation}
Adding such a last diagonal step \emph{does not add} any right-angle bend in
the s-path, so that the assumed incidence properties of the initial sequence
\(\Sigma^1 = (\pi^1_1, \ldots, \pi^1_{2r_1})\) are essentially preserved in the
completed sequence \(\overline{\Sigma}^1 = (\overline{\pi}^1_1, \ldots,
\overline{\pi}^1_{2r_1})\): the faces of each s-path are already present in the
sequence or are exterior; in the even case \(\pi^1_{2i} \in \Sigma^1\), the
previous s-path \(\pi^1_{2i-1}\) is a face of \(\pi^1_{2i}\) and hence
\(\overline{\pi}^1_{2i-1}\) is a face of \(\overline{\pi}^1_{2i}\). For example
in the illustration above, if \(i\) is even, certainly \(\partial_1 \pi^1_ i =
\pi^1_{i-1}\) (for this face is the only interior face) and this implies also
\(\partial_2 \overline{\pi}^1_ i = \overline{\pi}^1_{i-1}\).

On the contrary, in the case \(j=2\), the completion process \emph{can} add one
right-angle bend, nomore. For example, in this illustration:
\begin{equation}\begin{tikzpicture} [scale = 0.3, font = \scriptsize, baseline = 0.4cm]
 \foreach \i in {0,1,2,3} {\foreach \j in {0,...,4} {\node (\i\j) at (\j,\i)
  {\(\bullet\)} ;} ;} ;
 \node [font = \normalsize, left] at (10) {\(\pi^2_ i =\hspace{5pt}\)} ;
 \draw [thick] (10.center) -- (11.center) -- (12.center) -- (22.center)
           -- (33.center) -- (34.center) ;
\begin{scope} [xshift = 15cm]
 \foreach \i in {0,...,3} {\foreach \j in {0,...,4} {\node (\i\j) at (\j,\i)
  {\(\bullet\)} ;} ;} ;
 \node [font = \normalsize, left] at (10) {\(\overline{\pi}^2_ i =\hspace{5pt}\)} ;
 \draw [thick] (00.center) -- (10.center) -- (11.center) -- (12.center) -- (22.center)
           -- (33.center) -- (34.center) ;
\end{scope}
\begin{scope} [xshift = 30cm]
  \foreach \i in {0,...,3} {\foreach \j in {0,...,4} {\node (\i\j) at (\j,\i)
      {\(\bullet\)} ;} ;} ;
  \node [font = \normalsize, left] at (10)
 {\(\partial_1\overline{\pi}^2_ i =\hspace{5pt}\)} ;
 \draw [thick] (00.center) -- (11.center) -- (12.center) -- (22.center)
           -- (33.center) -- (34.center) ;
\end{scope}
\end{tikzpicture}\end{equation}
If the index \(i\) is even, then \(\partial_2 \pi^2_ i = \pi^2_{i-1}\) and the
relation \(\partial_3 \overline{\pi}^2_ i = \overline{\pi}^2_{i-1}\) is
satisfied as well. But another face of \(\overline{\pi}^2_ i\) is interior,
namely \(\partial_1 \overline{\pi}^2_i\), generated by the new right-angle
bend; because of the diagonal nature of the first step of this face, this face
is present in the list \(\overline{\Sigma}^1\), see the previous illustration.

Which is explained about \(\overline{\Sigma}^2\) with respect to
\(\overline{\Sigma}^1\) is valid as well for~\(\overline{\Sigma}^3\) with
respect to \(\overline{\Sigma}^1\).

The so-called last simplices, see Definition~\ref{07444}, must not be
forgotten! The last simplex \(\lambda_{p-1,q-1}\) (resp. \(\lambda_{p,q-1}\))
\emph{is not} in the list \(\Sigma_1\) (resp. \(\Sigma_2\)): these lists
describe the contractions of the \(hollowed\) prisms over the corresponding
boundaries: all the interior simplices are in these lists except the last ones.
The figure below gives these simplices in the case \((p,q) = (4,3)\):
\begin{equation}\begin{tikzpicture} [scale = 0.3, font = \scriptsize, thick,
   baseline = 0.2cm]
 \foreach \i in {0,...,2,3} {\foreach \j in {0,...,3,4}
   {\node (\i\j) at (\j,\i) {\(\bullet\)} ;} ;} ;
 \node [left, font = \normalsize] at (10) {\(\lambda_{3,2} = \hspace{5pt}\)} ;
 \draw (11.center) -- (14.center) -- (34.center) ;
\begin{scope} [xshift = 15cm]
 \foreach \i in {0,...,2,3} {\foreach \j in {0,...,4}
   {\node (\i\j) at (\j,\i) {\(\bullet\)} ;} ;} ;
 \node [left, font = \normalsize] at (10) {\(\lambda_{4,2} = \hspace{5pt}\)} ;
 \draw (10.center) -- (14.center) -- (34.center) ;
\end{scope}
\end{tikzpicture}\end{equation}
Examining now the respective completed paths:
\begin{equation}\begin{tikzpicture} [scale = 0.3, font = \scriptsize, thick,
   baseline = 0.2cm]
 \foreach \i in {0,...,2,3} {\foreach \j in {0,...,3,4}
   {\node (\i\j) at (\j,\i) {\(\bullet\)} ;} ;} ;
 \node [left, font = \normalsize] at (10) {\(\overline{\lambda}_{3,2} = \hspace{5pt}\)} ;
 \draw (00.center) -- (11.center) -- (14.center) -- (34.center) ;
\begin{scope} [xshift = 15cm]
 \foreach \i in {0,...,2,3} {\foreach \j in {0,...,4}
   {\node (\i\j) at (\j,\i) {\(\bullet\)} ;} ;} ;
 \node [left, font = \normalsize] at (10) {\(\overline{\lambda}_{4,2} = \hspace{5pt}\)} ;
 \draw (00.center) -- (10.center) -- (14.center) -- (34.center) ;
\end{scope}
\end{tikzpicture}\end{equation}
shows that \(\partial_1\overline{\lambda}_{p,q-1} =
\overline{\lambda}_{p-1,q-1}\); also the faces \(\partial_{p}
\overline{\lambda}_{p-1,q-1}\) and \(\partial_{p+1}
\overline{\lambda}_{p,q-1}\) are respectively in \(\overline{\Sigma}^1\) and
\(\overline{\Sigma}^2\).

Putting together all these facts leads to the conclusion: If \(\Sigma^1\),
\(\Sigma^2\) and \(\Sigma^3\) are respective filling sequences for
\((H\Delta^{(\ast,\ast)} - \partial \Delta^{(\ast,\ast)})\), with \((\ast,\ast)
= (p-1,q-1)\), \((p,q-1)\) and \((p-1,q)\) then the following list \emph{is a
filling sequence proving the desired W-contraction property} for the indices
\((p,q)\):
\begin{equation}
\overline{\Sigma}^1 \ ||\ \  \overline{\Sigma}^2 \ ||\ \
(\overline{\lambda}_{p-1,q-1}, \overline{\lambda}_{p-1,q}) \ ||\ \
\overline{\Sigma}^3
\end{equation}
where `\(||\)' is the list concatenation. Fortunately, the last simplex
\(\lambda_{p,q} = \overline{\lambda}_{p-1,q}\) is the only interior simplex of
\(\Delta^{p,q}\) missing in this list.\qed

\subsection{Examples.}

The reader can apply himself the above algorithm for the small dimensions. The
table below gives all the results for \((p,q) \leq (2,2)\).
\begin{equation}\begin{tikzpicture} [sc/.style = {font = \scriptsize}, thick,
                                  scale = 0.3333]
 \node at (23,0) {} ;
 \foreach \i in {0,1} {\foreach \j in {0,1}
   {\node [sc]  (\i\j) at (\j,\i) {\(\bullet\)} ;} ;} ;
 \node [left] at (-2,0.5) {\((p,q) = (1,1)\)} ;
 \draw (00.center) -- (11.center) ;
 \begin{scope} [xshift = 4cm]
 \foreach \i in {0,1} {\foreach \j in {0,1}
   {\node [sc]  (\i\j) at (\j,\i) {\(\bullet\)} ;} ;} ;
 \draw (00.center) -- (10.center) -- (11.center) ;
 \end{scope}
\end{tikzpicture}\end{equation}
\begin{equation}\begin{tikzpicture} [sc/.style = {font = \scriptsize}, thick,
                                  scale = 0.3333]
 \node at (23,0) {} ;
 \foreach \i in {0,1} {\foreach \j in {0,1,2}
   {\node [sc]  (\i\j) at (\j,\i) {\(\bullet\)} ;} ;} ;
 \node [left] at (-2,0.5) {\((p,q) = (2,1)\)} ;
 \draw (00.center) -- (11.center) -- (12.center) ;
 \begin{scope} [xshift = 4cm]
 \foreach \i in {0,1} {\foreach \j in {0,1,2}
   {\node [sc]  (\i\j) at (\j,\i) {\(\bullet\)} ;} ;} ;
 \draw (00.center) -- (10.center) -- (11.center) -- (12.center) ;
 \end{scope}
 \begin{scope} [xshift = 8cm]
 \foreach \i in {0,1} {\foreach \j in {0,1,2}
   {\node [sc]  (\i\j) at (\j,\i) {\(\bullet\)} ;} ;} ;
 \draw (00.center) -- (01.center) -- (12.center) ;
 \end{scope}
 \begin{scope} [xshift = 12cm]
 \foreach \i in {0,1} {\foreach \j in {0,1,2}
   {\node [sc]  (\i\j) at (\j,\i) {\(\bullet\)} ;} ;} ;
 \draw (00.center) -- (01.center) -- (11.center) -- (12.center) ;
 \end{scope}
\end{tikzpicture}\end{equation}
\begin{equation}\begin{tikzpicture} [sc/.style = {font = \scriptsize}, thick,
                                  scale = 0.3333]
 \node at (23,0) {} ;
 \foreach \i in {0,1,2} {\foreach \j in {0,1}
   {\node [sc]  (\i\j) at (\j,\i) {\(\bullet\)} ;} ;} ;
 \node [left] at (-2,0.5) {\((p,q) = (1,2)\)} ;
 \draw (00.center) -- (10.center) -- (21.center) ;
 \begin{scope} [xshift = 4cm]
 \foreach \i in {0,1,2} {\foreach \j in {0,1}
   {\node [sc]  (\i\j) at (\j,\i) {\(\bullet\)} ;} ;} ;
 \draw (00.center) -- (10.center) -- (20.center) -- (21.center) ;
 \end{scope}
 \begin{scope} [xshift = 8cm]
 \foreach \i in {0,1,2} {\foreach \j in {0,1}
   {\node [sc]  (\i\j) at (\j,\i) {\(\bullet\)} ;} ;} ;
 \draw (00.center) -- (11.center) -- (21.center) ;
 \end{scope}
 \begin{scope} [xshift = 12cm]
 \foreach \i in {0,1,2} {\foreach \j in {0,1}
   {\node [sc]  (\i\j) at (\j,\i) {\(\bullet\)} ;} ;} ;
 \draw (00.center) -- (10.center) -- (11.center) -- (21.center) ;
 \end{scope}
\end{tikzpicture}\end{equation}
\begin{equation}\begin{tikzpicture} [sc/.style = {font = \scriptsize}, thick,
                                  scale = 0.3333]
 \node at (23,0) {} ;
 \foreach \i in {0,1,2} {\foreach \j in {0,1,2}
   {\node [sc]  (\i\j) at (\j,\i) {\(\bullet\)} ;} ;} ;
 \node [left] at (-2,0.5) {\((p,q) = (2,2)\)} ;
 \draw (00.center) -- (11.center) -- (22.center) ;
 \begin{scope} [xshift = 4cm]
 \foreach \i in {0,1,2} {\foreach \j in {0,1,2}
   {\node [sc]  (\i\j) at (\j,\i) {\(\bullet\)} ;} ;} ;
 \draw (00.center) -- (11.center) -- (21.center) -- (22.center) ;
 \end{scope}
 \begin{scope} [xshift = 8cm]
 \foreach \i in {0,1,2} {\foreach \j in {0,1,2}
   {\node [sc]  (\i\j) at (\j,\i) {\(\bullet\)} ;} ;} ;
 \draw (00.center) -- (10.center) -- (21.center) -- (22.center) ;
 \end{scope}
 \begin{scope} [xshift = 12cm]
 \foreach \i in {0,1,2} {\foreach \j in {0,1,2}
   {\node [sc]  (\i\j) at (\j,\i) {\(\bullet\)} ;} ;} ;
 \draw (00.center) -- (10.center) -- (20.center) -- (21.center) -- (22.center) ;
 \end{scope}
 \begin{scope} [xshift = 16cm]
 \foreach \i in {0,1,2} {\foreach \j in {0,1,2}
   {\node [sc]  (\i\j) at (\j,\i) {\(\bullet\)} ;} ;} ;
 \draw (00.center) -- (10.center) -- (11.center) -- (22.center) ;
 \end{scope}
 \begin{scope} [xshift = 20cm]
 \foreach \i in {0,1,2} {\foreach \j in {0,1,2}
   {\node [sc]  (\i\j) at (\j,\i) {\(\bullet\)} ;} ;} ;
 \draw (00.center) -- (10.center) -- (11.center) -- (21.center) -- (22.center) ;
 \end{scope}
 \begin{scope} [yshift = -4.5cm]
 \foreach \i in {0,1,2} {\foreach \j in {0,1,2}
   {\node [sc]  (\i\j) at (\j,\i) {\(\bullet\)} ;} ;} ;
 \draw (00.center) -- (11.center) -- (12.center) -- (22.center) ;
 \begin{scope} [xshift = 4cm]
 \foreach \i in {0,1,2} {\foreach \j in {0,1,2}
   {\node [sc]  (\i\j) at (\j,\i) {\(\bullet\)} ;} ;} ;
 \draw (00.center) -- (10.center) -- (11.center) -- (12.center) -- (22.center) ;
 \end{scope}
 \begin{scope} [xshift = 8cm]
 \foreach \i in {0,1,2} {\foreach \j in {0,1,2}
   {\node [sc]  (\i\j) at (\j,\i) {\(\bullet\)} ;} ;} ;
 \draw (00.center) -- (01.center) -- (11.center) -- (22.center) ;
 \end{scope}
 \begin{scope} [xshift = 12cm]
 \foreach \i in {0,1,2} {\foreach \j in {0,1,2}
   {\node [sc]  (\i\j) at (\j,\i) {\(\bullet\)} ;} ;} ;
 \draw (00.center) -- (01.center) -- (11.center) -- (21.center) -- (22.center) ;
 \end{scope}
 \begin{scope} [xshift = 16cm]
 \foreach \i in {0,1,2} {\foreach \j in {0,1,2}
   {\node [sc]  (\i\j) at (\j,\i) {\(\bullet\)} ;} ;} ;
 \draw (00.center) -- (01.center) -- (12.center) -- (22.center) ;
 \end{scope}
 \begin{scope} [xshift = 20cm]
 \foreach \i in {0,1,2} {\foreach \j in {0,1,2}
   {\node [sc]  (\i\j) at (\j,\i) {\(\bullet\)} ;} ;} ;
 \draw (00.center) -- (01.center) -- (11.center) -- (12.center) -- (22.center) ;
 \end{scope}
 \end{scope}
\end{tikzpicture}\end{equation}

For significantly bigger values of \((p,q)\), only a program can produce the
corresponding filling sequences. A short Lisp program (45 lines) can produce
the justifying list for reasonably small values of \(p\) and \(q\). For
example, if \(p = q =8\), the filling sequence is made of 265,728 paths, a list
produced in 4 seconds on a modest laptop. But if \(p = q = 10\), the number of
paths is 8,097,452; and the same laptop is then out of memory. A typical
behaviour in front of exponential complexity: the necessary number of paths is
\(> 3^p\) if \(p=q>1\).

\subsection{If homology is enough.} \label{03005}

The proposed proof of Theorem~\ref{79523} is elementary but a little technical.
If you are only interested by the \emph{homological} Eilenberg-Zilber theorem,
a simple proof analogous to Proposition~\ref{77213}'s is sufficient.

\begin{prp} --- \label{79914}
The relative chain complex \(C_\ast(H\Delta^{p,q}, \partial\Delta^{p,q})\)
admits a W-reduction to the null complex.
\end{prp}

\prof This relative chain complex is generated by all the interior simplices
\(\sigma\) of \(\Delta^{p,q}\) except the last one. Representing such a simplex
\(\sigma\) by the corresponding s-path~\(\pi\) allows us to divide all these
simplices into two disjoint sets \(D\) and \(R\). The rule is the following:
you run the examined path \(\pi\) \emph{from} \((p,q)\) \emph{backward} to
\((0,0)\) in the lattice \(\up \times \uq\) and you are interested by the
\emph{first} ``event'':
\begin{fenum}
\item
Either you run a diagonal elementary step\ \ \begin{tikzpicture} [->, scale =
0.25, thick] \draw (1,1) -- (0,0) ; \end{tikzpicture}\ \ , in which case the
path \(\pi\) is in \(D\);
\item
Or you pass a \(+90^\circ\) bend\ \ \begin{tikzpicture} [->, scale = 0.25,
thick] \draw (1,1) -- (0,1) -- (0,0) ; \end{tikzpicture}\ \ (not a
\(-90^{\circ}\) bend\ \ \begin{tikzpicture} [->, scale = 0.25, thick] \draw
(1,1) -- (1,0) -- (0,0) ; \end{tikzpicture}\ \ ) in which case the path \(\pi\)
is in~\(R\).
\end{fenum}
\begin{equation}\begin{tikzpicture} [baseline = 0.5cm, scale = 0.4]
 \foreach \i in {0,1,2,3} \foreach \j in {0,1,2,3}
   {\coordinate (\i\j) at (\i,\j) ; \node at (\i\j) {\footnotesize\(\bullet\)}
   ;}
 \draw [thick] (00) -- (01) -- (11) (22) -- (32) -- (33) ;
 \draw [ultra thick, dotted] (11) -- (22) ;
 \node [anchor = 0] at (-0.5,1.5) {\(\pi_1 =\)} ;
 \node [anchor = 180] at (3.5,1.5) {\(\in D\)} ;
\begin{scope} [xshift = 15cm]
 \foreach \i in {0,1,2,3} \foreach \j in {0,1,2,3}
   {\coordinate (\i\j) at (\i,\j) ; \node at (\i\j) {\footnotesize\(\bullet\)}
   ;}
 \draw [thick] (00) -- (01) -- (11) (22) -- (32) -- (33) ;
 \draw [ultra thick, dotted] (11) -- (12) -- (22) ;
 \node [anchor = 0] at (-0.5,1.5) {\(\pi_2 =\)} ;
 \node [anchor = 180] at (3.5,1.5) {\(\in R\)} ;
 \end{scope}
\end{tikzpicture}
\end{equation}

The above figure displays one example in either case; the deciding ``event'' is
signalled by dotted lines. Observe \(\partial_3 \pi_2 = \pi_1\). More generally
it is clear the operator assigning to every path of \(R\) the face in \(D\)
corresponding to the \(+90^{\circ}\) bend is a bijection organizing all these
paths by pairs defining a discrete vector field, a good candidate to construct
the announced W-reduction: the unique path without any event fortunately is the
last simplex!

There remains to prove this vector field is admissible.  It is a consequence of
the organization of this vector field as a filling sequence given in
Section~\ref{68696}, but it is possible to prove it directly and simply.

The example of the vector \((\pi_1, \pi_2)\) above is enough to understand. We
have to consider the faces of \(\pi_2\) which are \emph{sources}, in other
words which are in \(D\), therefore in particular \emph{interior}, and
\emph{different} from \(\pi_1\). In general at most two faces satisfy these
requirements, here these faces \(\pi_3 = \partial_2 \pi_2\) and \(\pi_4 =
\partial_5 \pi_2\):
\begin{equation}\begin{tikzpicture} [baseline = 0.5cm, scale = 0.4]
 \foreach \i in {0,1,2,3} \foreach \j in {0,1,2,3}
   {\coordinate (\i\j) at (\i,\j) ; \node at (\i\j) {\footnotesize\(\bullet\)}
   ;}
 \draw [thick] (00) -- (01) -- (12) -- (22) -- (32) -- (33) ;
 \node [anchor = 0] at (-0.5,1.5) {\(\pi_3 = \partial_2 \pi_2 =\)} ;
 \node [anchor = 180] at (3.5,1.5) {\(\in D\)} ;
\begin{scope} [xshift = 15cm]
 \foreach \i in {0,1,2,3} \foreach \j in {0,1,2,3}
   {\coordinate (\i\j) at (\i,\j) ; \node at (\i\j) {\footnotesize\(\bullet\)}
   ;}
 \draw [thick] (00) -- (01) -- (11) -- (12) -- (22) -- (33) ;
 \node [anchor = 0] at (-0.5,1.5) {\(\pi_4 = \partial_5 \pi_2 =\)} ;
 \node [anchor = 180] at (3.5,1.5) {\(\in D\)} ;
\end{scope}
\end{tikzpicture}
\end{equation}

This gives a Lyapunov function, see Definition~\ref{92078}. If \(\pi \in D\),
decide \(L(\pi)\) is the number of points of the \((p,q)\)-lattice strictly
\emph{above} the path. Observe \(L(\pi_1) = 5\) while \(L(\pi_3) = L(\pi_4) =
4\). A reader having reached this point of the text will probably prefer to
play in ending the proof by himself.\qed

\subsection{The Eilenberg-Zilber W-reduction.}

\subsubsection
  [\texorpdfstring{The $(p,q)$-Eilenberg-Zilber reduction.}{The (p,q)
  Eilenberg-Zilber reduction}]
  {The \boldmath$(p,q)$-Eilenberg-Zilber reduction.}

Proposition~\ref{79914} can be arranged to produce a sort of top-dimensional
Eilenberg-Zilber W-reduction for the prism \(\Delta^{p,q}\).

\begin{prp} ---
The discrete vector field used in Proposition~\ref{79914} induces a W-reduction
\(C_\ast(\Delta^{p,q}) \rrdc C_\ast^c(\Delta^{p,q})\) where in particular the
chain group \(C_{p+q}(\Delta^{p,q})\) of rank \begin{tikzpicture} [yscale =
0.15, inner sep = 0pt, font = \scriptsize, baseline = (b.base)]
 \node (a) at (0,1) {\(p+q\)} ; \node (b) at (0,-1) {\(p\,,\,q\)} ;
 \node [anchor = 0, yscale=2] at (a.west |- 0,0) {\emph{(}} ;
 \node [anchor = 180, yscale=2] at (a.east |- 0,0) {\emph{)}} ;
\end{tikzpicture}
is replaced by a critical chain group with a \emph{unique} generator, the
so-called last simplex \(\lambda_{p,q}\).
\end{prp}

\prof This vector field makes sense as well in this context as in
Proposition~\ref{79914} and the admissibility property remains valid. A
W-reduction is therefore generated, where in dimension \((p+q)\) the only
critical simplex is the last simplex \(\lambda_{p,q}\).\qed

\subsubsection{Products of simplicial sets.} \label{49366}

Let us recall the product \(X \times Y\) of two simplicial sets \(X\) and \(Y\)
is very simply defined. These simplicial sets \(X\) and \(Y\) are nothing but
contravariant functors \(\uD \rightarrow \textrm{\underline{Sets}}\), and the
simplicial set \(X \times Y\) is the \emph{product functor}. In particular \((X
\times Y)_p = X_p \times Y_p\) and if \(\alpha: \up \rightarrow \uq\) is a
\(\uD\)-morphism, then \(\alpha^\ast_{X \times Y} = \mbox{\(\alpha^\ast_X
\times \alpha^\ast_Y:\)}\ \mbox{\((X \times Y)_q\)} \rightarrow \mbox{\((X
\times Y)_p\)}\). It is not obvious when seeing this definition for the first
time this actually corresponds to the standard notion of topological product
but, except in esoteric cases when simplex sets are not countable, the
topological realization of the product is homeomorphic to the product of
realizations. In these exceptional cases, the result remains true under the
condition of working in the category of compactly generated
spaces~\cite{stnr2}.

A \(p\)-simplex of the product \(X \times Y\) is therefore a \emph{pair} \(\rho
= (\sigma, \tau)\) of \(p\)-simplices of \(X\) and \(Y\), and it is important
to understand when this simplex is degenerate or not. Taking account of the
Eilenberg-Zilber Lemma~\ref{34608}, we prefer to express both components of
this pair as the degeneracy of some non-degenerate simplex, which produces the
expression \((\eta_{i_{s-1}} \cdots\, \eta_{i_{0}} \sigma', \eta_{j_{t-1}}
\cdots\, \eta_{j_{0}} \tau')\) for our \(p\)-simplex \(\rho\),
where~\(\sigma'\) (resp. \(\tau'\)) is a non-degenerate \((p-s)\)-simplex of
\(X\) (resp. \((p-t)\)-simplex of \(Y\)); also the sequences \(i_\ast\) and
\(j_\ast\) must be \emph{strictly increasing} with respect to their indices.
Then the algebra of the elementary degeneracies \(\eta_i\) shows the simplex
\(\rho = (\sigma, \tau)\) is non-degenerate if and only if the intersection
\(\{i_{s-1}, \ldots, i_0\} \cap \{j_{t-1}, \ldots, j_0\}\) is empty.

The next definition is a division of all the non-degenerate simplices of the
product \(Z = X \times Y\) into three parts: the target simplices \(Z_\ast^t\),
the source simplices \(Z_\ast^s\) and the critical simplices \(Z_\ast^c\). This
division corresponds to a discrete vector field \(V\), the natural extension to
the whole product \(Z = X \times Y\) of the vector field constructed in
Sections~\ref{68696} and~\ref{03005} for the top bidimension \((p,q)\) of the
prism \(\Delta^{p,q}\).

We must translate the definition of the vector field \(V\) in
Section~\ref{03005} into the language of non-degenerate product simplices
expressed as pairs of possible degeneracies. To prepare the reader at this
translation, let us explain the recipe which translates an s-path into such a
pair. Let us consider this s-path:
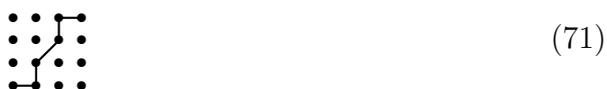
\begin{equation}\begin{tikzpicture} [scale = 0.3, font = \scriptsize, baseline = 5mm]
  \foreach \i in {0,...,3} \foreach \j in {0,...,3}
    {\coordinate (\i\j) at (\i,\j) ;
     \node at (\i\j) {\(\bullet\)} ;}
  \draw [thick] (00) -- (10) -- (11) -- (22) -- (23) -- (33) ;
\end{tikzpicture}\end{equation}

This s-path represents an interior 5-simplex of \(\Delta^{3,3}\) to be
expressed in terms of the maximal simplices \(\sigma, \tau \in \Delta^3_{3}\).
Run this s-path from \((0,0)\) to \((3,3)\); every vertical elementary step,
for example from \((1,0)\) to \((1,1)\) produces a degeneracy in the first
factor, the index being the time when this vertical step is started, here~1.
Another vertical step starts from \((2,2)\) at time 3, so that the first factor
will be \(\eta_3\eta_1\sigma\). In the same way, examining the horizontal steps
produces the second factor \(\eta_4\eta_0\tau\). Finally our s-path codes the
simplex \((\eta_3 \eta_1 \sigma, \eta_4 \eta_0 \tau)\). The index 2 is missing
in the degeneracies, meaning that at time 2 the corresponding step is diagonal:
the lists of degeneracy indices are directly connected to the structure of the
corresponding s-path. We will say the \emph{degeneracy configuration} of this
simplex is \(((3,1), (4,0))\); a degeneracy configuration is a pair of disjoint
decreasing integer lists.

Conversely, reading the indices of the degeneracy operators in the canonical
writing \(\rho = (\eta_{i_{s-1}} \cdots\, \eta_{i_{0}} \sigma, \eta_{j_{t-1}}
\cdots\, \eta_{j_{0}} \tau)\) of a simplex \(\rho\) of \(X \times Y\)
unambiguously describes the corresponding s-path.

This process settles a canonical bijection between \(S_{p,q}\) and \(D_{p,q}\)
if:
\begin{fenum}
\item
The set \(S_{p,q}\) is the collection of all the \emph{interior} s-paths
running from \((0,0)\) to \((p,q)\) in the \((p,q)\)-lattice.
\item
The set \(D_{p,q}\) is the collection of all the configurations of degeneracy
operators which can be used when writing a non-degenerate simplex in its
canonical form \(\rho = (\eta_{i_{s-1}} \cdots\, \eta_{i_{0}} \sigma,
\eta_{j_{t-1}} \cdots\, \eta_{j_{0}} \tau)\), when \(\sigma \in X_p^{ND}\) and
\(\tau \in Y_q^{ND}\). A configuration is a pair of integer lists \(((i_{s-1},
\ldots, i_0), (j_{t-1}, \ldots, j_{0}))\) satisfying the various coherence
conditions explained before: \(p+s = q+t\), every component \(i_-\) and \(j_-\)
is in \([0\ldots(p+q-1)]\), both lists are disjoint,  and their elements are
increasing with respect to their respective indices.
\end{fenum}

\subsubsection{The Eilenberg-Zilber vector field.} \label{08795}

\begin{dfn} --- \label{12047}
\emph{Let \(\rho = (\eta_{i_{s-1}} \cdots\, \eta_{i_{0}} \sigma, \eta_{j_{t-1}}
\cdots\, \eta_{j_{0}} \tau)\) be a non-degenerate \(p\)-simplex of \(X \times
Y\) written in the canonical form. The degeneracy configuration \(((i_{s-1},
\ldots, i_0), (j_{t-1}, \ldots, j_{0}))\) is a well-defined element of
\(D_{p-s,p-t}\) which in turn defines an s-path \(\pi(\rho) \in S_{p-s,p-t}\).
Then \(\rho\) is a \emph{target} (resp. \emph{source}, \emph{critical}) simplex
if and only if the s-path \(\pi(\rho)\) has the corresponding property.}\qed
\end{dfn}

\begin{dfn} ---
\emph{If \(\rho = (\eta_{i_{s-1}} \cdots\, \eta_{i_{0}} \sigma, \eta_{j_{t-1}}
\cdots\, \eta_{j_{0}} \tau)\) is a \(p\)-simplex of \(X \times Y\) as in the
previous definition, the pair \((p-s,p-t)\) is called the \emph{bidimension of}
\(\rho\). It is the bidimension of the smallest prism of \(X \times Y\)
containing this simplex.}\qed
\end{dfn}

For example the diagonal of the square \(\Delta^{1,1}\) has the bidimension
\((1,1)\). The sum of the components of the bidimension can be bigger than the
dimension.

\begin{dfn} --- \label{52397}
\emph{ Let \(X \times Y\) be the product of two simplicial sets. The division
of the non-degenerate simplices of \(X \times Y\) according to
Definition~\ref{12047} into \emph{target} simplices, \emph{source} simplices
and \emph{critical} simplices, combined with the pairing described in the proof
of Proposition~\ref{79914}, defines the \emph{Eilenberg-Zilber vector field}
\(V_{X \times Y}\) of \(X \times Y\).}\qed
\end{dfn}

\begin{thr} \textbf{\emph{(Eilenberg-Zilber Theorem)}} --- \label{46504}
Let \(X \times Y\) be the product of two simplicial sets. The Eilenberg-Zilber
vector field \(V_{X \times Y}\) induces the Eilenberg-Zilber W-reduction:
\begin{equation}
EZ: C_\ast^{ND}(X \times Y) \rrdc C_\ast^{ND}(X) \otimes C_\ast^{ND}(Y).
\end{equation}
\end{thr}

The reader may wonder why all these technicalities to reprove a well-known
sixty years old theorem. Two totally different reasons.

On the one hand, the Eilenberg-Zilber reduction is time consuming when
concretely programmed. In particular, profiler examinations of the effective
homology programs show the terrible homotopy component \(h: C_\ast(X \times Y)
\rightarrow C_\ast(X \times Y)\) of the Eilenberg-Zilber reduction, rarely
seriously considered\footnote{With two notable exceptions. In the landmark
papers by Eilenberg and\ldots MacLane~\mbox{\cite{elmc1,elmc2}}, more useful
than the standard reference~\cite{elzl2}, a nice recursive description of this
homotopy operator is given. Forty years later (\(!\)), when a computer program
was at last available to make experiments, Julio Rubio found a closed formula
for this operator, proved by Fr\'ed\'eric Morace a little later~\cite{real2}.
We reprove this formula and others in the next section, by a totally elementary
process depending only on our vector field, independent of Eilenberg-MacLane's
and Shih's recursive formulas, reproved as well.}, is the kernel program unit
the most used in concrete computations. Our description of the Eilenberg-Zilber
reduction makes the corresponding program unit simpler and more efficient.

On the other hand, maybe more important, the \emph{same} (\(!\)) vector field
will be soon used to process in the same way the \emph{twisted} products,
leading to totally elementary \emph{effective} versions of the Serre and
Eilenberg-Moore spectral sequences. Maybe the same for the Bousfield-Kan
spectral sequence.

\prof The vector field \(V_{X \times Y}\) has a layer for every bidimension
\((p,q)\). The admissibility proof given in Proposision~\ref{79914} shows that
every V-path starting from a source simplex of bidimension \((p,q)\) goes after
a finite number of steps to sub-layers. The Eilenberg-Zilber vector field is
admissible.

If \(\sigma\) (resp. \(\tau\)) is a non-degenerate \(p\)-simplex of \(X\)
(resp. \(q\)-simplex of \(Y\)), we can denote by \(\sigma \times \tau\) the
corresponding prism in \(X \times Y\), made of all the simplices of bidimension
\((p,q)\) with respect to \(\sigma\) and \(\tau\). The collection of the
\emph{interior} simplices of this prism \(\sigma \times \tau\) is nothing but
an \emph{exact} copy of the collection of the \emph{interior} simplices of the
standard prism \(\Delta^{p,q}\). In particular only one \emph{interior}
critical cell in every prism. You are attending the birth of the tensor product
\(C_\ast^{ND}(X) \otimes C_\ast^{ND}(Y)\): exactly one generator \(\sigma
\otimes \tau\) for every prism \(\sigma \times \tau\), namely the last simplex
of this prism: \(\sigma \otimes \tau\ \textrm{``=''}\ \lambda_{p,q}(\sigma
\times \tau)\).

There remains to prove the small chain complex so obtained  is not only the
right graded module \(C_\ast^{ND}(X) \otimes C_\ast^{ND}(Y)\), but is endowed
by the reduction process with the right differential. This is a corollary of
the next section, devoted to a detailed study of the Eilenberg-Zilber vector
field.\qed

We have the right generators for the reduced complex, but do we have the right
differential? The proof is again a little technical. A small movie of what
happens for \(\Delta^{3,2}\) is enough. The last simplex \(\lambda_{3,2}\) is
this one:
\begin{equation}\begin{tikzpicture} [scale = 0.3, font = \scriptsize, baseline = 0.3cm]
 \foreach \i in {0,1,2,3} \foreach \j in {0,1,2}
   {\coordinate (\i\j) at (\i,\j) ; \node at (\i\j) {\(\bullet\)} ; }
   \draw [thick] (0,0) -- (3,0) -- (3,2) ;
\end{tikzpicture}
\end{equation}
We must use the first formula~(\ref{76423}): \(d' = d_{3,3} - d_{3,1}
d_{2,1}^{-1} d_{2,3}\), a good opportunity to illustrate how it works.

In the prism \(\Delta^{3,2}\), the last simplex \(\lambda_{3,2}\) of dimension
5 has six faces. Five of these faces have a smaller bidimension and are
\emph{critical} cells, one face, namely the \(\partial_3\)-face has the same
bidimension, but it is a \emph{source} cell.

For example the 1-face \(\partial_1 \lambda_{3,2}\) is represented by this
path:
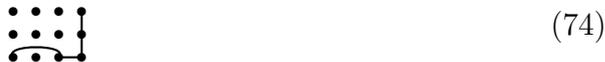
\begin{equation}\begin{tikzpicture} [scale = 0.3, font = \scriptsize, baseline = 0.3cm]
 \foreach \i in {0,1,2,3} \foreach \j in {0,1,2}
   {\coordinate (\i\j) at (\i,\j) ; \node at (\i\j) {\(\bullet\)} ; }
   \draw [thick] (0,0) to [in = 50, out = 130] (2,0) -- (3,0) -- (3,2) ;
\end{tikzpicture}
\end{equation}
of bidimension \((2,2)\), it is again a last simplex, but in a prism
\(\Delta^{2,2}\). For our last simplex \(\lambda_{3,2}(\sigma \times \tau)\),
we obtain \(\partial_1 \lambda_{3,2}(\sigma \times \tau) =
\lambda_{2,2}(\partial_1 \sigma \times \tau)\ \mbox{``=''}\ \partial_1 \sigma
\otimes \tau\), unless the face \(\partial_1 \sigma\) is degenerate in \(X\),
possible, in which case the 1-face does not contribute in the differential. The
same for the faces \(\partial_0\), \(\partial_2\), \(\partial_4\) and
\(\partial_5\), giving the contribution of the term \(d_{3,3}\) in our formula
for the differential \(d'\) to be computed. You can easily verify the right
signs are obtained.

Let us consider now the 3-face of \(\lambda_{3,2}\). It is the following source
cell:
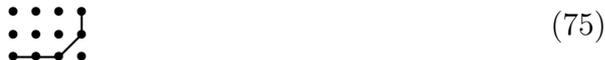
\begin{equation}\begin{tikzpicture} [scale = 0.3, font = \scriptsize, baseline = 0.3cm]
 \foreach \i in {0,1,2,3} \foreach \j in {0,1,2}
   {\coordinate (\i\j) at (\i,\j) ; \node at (\i\j) {\(\bullet\)} ; }
   \draw [thick] (0,0) -- (2,0) -- (3,1) -- (3,2) ;
\end{tikzpicture}
\end{equation}
The fact that this is a source cell implies the corresponding part of the
differential is given by the second term \(- d_{3,1} d_{2,1}^{-1} d_{2,3}\):
the last term \(d_{2,3}\) of this composition is used for a critical cell, our
last simplex \(\lambda_{3,2}\), index 3 in \(d_{2,3}\), producing a source
cell, index 2 in \(d_{2,3}\), the source cell displayed above.

After \(d_{2,3}\) we must apply the inverse \(d_{2,1}^{-1}\), really the key
point in all these calculations. We must apply the recursive
formula~(\ref{08522}). This source cell is paired with this target cell:
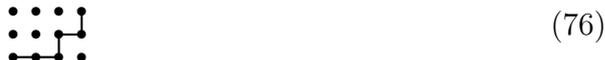
\begin{equation}\begin{tikzpicture} [scale = 0.3, font = \scriptsize, baseline = 0.3cm]
 \foreach \i in {0,1,2,3} \foreach \j in {0,1,2}
   {\coordinate (\i\j) at (\i,\j) ; \node at (\i\j) {\(\bullet\)} ; }
   \draw [thick] (0,0) -- (2,0) -- (2,0) -- (2,1) -- (3,1) -- (3,2) ;
\end{tikzpicture}
\end{equation}
The game rule now is the following: please consider the faces of this cell. The
critical faces contribute, think of the component \(d_{3,1}\) in the formula
for \(d'\), but here no critical face. The target faces, here the faces of
index 0, 1 and 5 are thrown away. There remain three source faces which lead to
continue to apply the recursive formula~(\ref{08522}). The 3-face is not to be
considered, for we \emph{come from this face} by the vector field: look at the
term ``\(-\{\sigma\}\)'' in the sum index of the formula~(\ref{08522}). The
source faces 2 and 4 are:
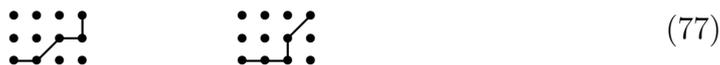
\begin{equation}\begin{tikzpicture} [scale = 0.3, font = \scriptsize, baseline = 0.3cm]
 \foreach \i in {0,1,2,3} \foreach \j in {0,1,2}
   {\coordinate (\i\j) at (\i,\j) ; \node at (\i\j) {\(\bullet\)} ; }
   \draw [thick] (0,0) -- (1,0) -- (2,1) -- (3,1) -- (3,2) ;
 \begin{scope} [xshift = 10cm]
 \foreach \i in {0,1,2,3} \foreach \j in {0,1,2}
   {\coordinate (\i\j) at (\i,\j) ; \node at (\i\j) {\(\bullet\)} ; }
 \draw [thick] (0,0) -- (2,0) -- (2,1) -- (3,2) ;
 \end{scope}
\end{tikzpicture}
\end{equation}
and the corresponding target cells are:
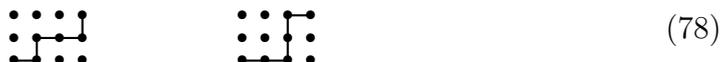
\begin{equation} \label{92825}
\begin{tikzpicture} [scale = 0.3, font = \scriptsize, baseline = 0.3cm]
 \foreach \i in {0,1,2,3} \foreach \j in {0,1,2}
   {\coordinate (\i\j) at (\i,\j) ; \node at (\i\j) {\(\bullet\)} ; }
   \draw [thick] (0,0) -- (1,0) -- (1,1) -- (2,1) -- (3,1) -- (3,2) ;
 \begin{scope} [xshift = 10cm]
 \foreach \i in {0,1,2,3} \foreach \j in {0,1,2}
   {\coordinate (\i\j) at (\i,\j) ; \node at (\i\j) {\(\bullet\)} ; }
 \draw [thick] (0,0) -- (2,0) -- (2,1) -- (2,2) -- (3,2) ;
 \end{scope}
\end{tikzpicture}
\end{equation}

Let us firstly examine the righthand s-path. The face 4 is not to be
considered, it is the source s-path we come from. The faces 0, 1, 2 and 3  are
target cells, to be thrown away. There remains the face 5:
\begin{equation}\begin{tikzpicture} [scale = 0.3, font = \scriptsize, baseline = 0.3cm]
 \foreach \i in {0,1,2,3} \foreach \j in {0,1,2}
   {\coordinate (\i\j) at (\i,\j) ; \node at (\i\j) {\(\bullet\)} ; }
   \draw [thick] (0,0) -- (2,0) -- (2,2) ;
\end{tikzpicture}
\end{equation}
which can also be written \(\lambda_{2,2}(\partial_3 \sigma \times \tau) =
\partial_3 \sigma \otimes \tau\), one of the still missing terms in the usual
boundary of \(\sigma \otimes \tau\). We think the reader is now able to
understand the end of the story for the left-hand term of the
figure~(\ref{92825}) as explained here:
\begin{equation}\begin{tikzpicture} [scale = 0.3, font = \scriptsize, baseline = 0.3cm]
 \foreach \i in {0,1,2,3} \foreach \j in {0,1,2}
   {\coordinate (\i\j) at (\i,\j) ; \node at (\i\j) {\(\bullet\)} ; }
   \draw [thick] (0,0) -- (1,0) -- (1,1) -- (2,1) -- (3,1) -- (3,2) ;
 \begin{scope} [xshift = 7cm]
 \foreach \i in {0,1,2,3} \foreach \j in {0,1,2}
   {\coordinate (\i\j) at (\i,\j) ; \node at (\i\j) {\(\bullet\)} ; }
 \draw [thick] (0,0) -- (1,1) -- (2,1) -- (3,1) -- (3,2) ;
 \end{scope}
 \begin{scope} [xshift = 14cm]
 \foreach \i in {0,1,2,3} \foreach \j in {0,1,2}
   {\coordinate (\i\j) at (\i,\j) ; \node at (\i\j) {\(\bullet\)} ; }
 \draw [thick] (0,0) -- (0,1) -- (3,1) -- (3,2) ;
 \end{scope}
 \begin{scope} [xshift = 21cm]
 \foreach \i in {0,1,2,3} \foreach \j in {0,1,2}
   {\coordinate (\i\j) at (\i,\j) ; \node at (\i\j) {\(\bullet\)} ; }
 \draw [thick] (0,1) -- (3,1) -- (3,2) ;
 \end{scope}
\end{tikzpicture}
\end{equation}
producing the last missing term \(\lambda_{3,1}(\sigma \times \partial_0 \tau)
= \sigma \otimes \partial_0 \tau\). The careful reader can object we have
``forgotten'' another source face of the first s-path above, namely:
\begin{equation} \label{57486}
\begin{tikzpicture} [scale = 0.3, font = \scriptsize, baseline = 0.3cm]
 \foreach \i in {0,1,2,3} \foreach \j in {0,1,2}
   {\coordinate (\i\j) at (\i,\j) ; \node at (\i\j) {\(\bullet\)} ; }
   \draw [thick] (0,0) -- (1,0) -- (1,1) -- (2,1) -- (3,2) ;
\end{tikzpicture}
\end{equation}
but it can be immediately thrown away, for the following reason: the only
possible critical cells which can appear in the final result are necessarily
locked in one of these boxes:
\begin{equation}\begin{tikzpicture} [scale = 0.3, font = \scriptsize, baseline = 0.3cm]
 \foreach \i in {0,1,2,3} \foreach \j in {0,1,2}
   {\coordinate (\i\j) at (\i,\j) ; \node at (\i\j) {\(\bullet\)} ; }
 \draw [thick] [thick] (-0.2,-0.2) -- (-0.2,1.2) -- (2.8,1.2) -- (2.8,2.2) -- (3.2,2.2)
    -- (3.2,-0.2) -- cycle ;
\begin{scope} [xshift = 10cm]
 \foreach \i in {0,1,2,3} \foreach \j in {0,1,2}
   {\coordinate (\i\j) at (\i,\j) ; \node at (\i\j) {\(\bullet\)} ; }
 \draw [thick] [thick] (-0.2,-0.2) -- (-0.2,0.2) -- (1.8,0.2) -- (1.8,2.2) -- (3.2,2.2)
   -- (3.2,-0.2) -- cycle ;
\end{scope}\end{tikzpicture}
\end{equation}
Now we observe the s-path (\ref{57486}) has escaped from these boxes; also the
recursive process always pushes a path in the north-west direction and it is
impossible the s-path (\ref{57486}) comes back inside one of our boxes and
contributes. The sceptical reader is invited to associate to an s-path the part
of the \((p,q)\)-\emph{rectangle} strictly above this s-path\footnote{This
argument will no longer be valid in the twisted case, for the 0-faces of the
product cells then become arbitrary.}.

Verifying the right signs are also obtained leads to the conclusion: in the
critical complex, the standard boundary formula \(d(\sigma \otimes \tau) = d
\sigma \otimes \tau + (-1)^{|\sigma|} \sigma \otimes d\tau\) is to be applied.
The critical chain complex produced by this W-reduction process is canonically
isomorphic to \(C_\ast^{ND}(X) \otimes C_\ast^{ND}(Y)\).\qed

It is natural to ask whether the traditional Alexander-Whitney,
Eilenberg-MacLane and Rubio-Morace formulas are obtained in our W-reduction.
The concrete programming work done for our Eilenberg-Zilber W-reductions gives
a positive \emph{experimental} answer, but we do not yet have a \emph{proof} of
this fact. Our decisions for the few allowed choices when defining the
Eilenberg-Zilber vector field were in fact done to obtain the traditional
formulas in some easy particular cases. The right method to obtain the general
hoped-for result is probably the following: the Alexander-Whitney,
Eilenberg-MacLane and Rubio-Morace~\cite[Section~6]{real2} formulas are
consequences of the recursive formulas detailed in~\cite{elmc1,elmc2}, and our
description of the vector field  can also be studied along the same recursive
process.

In fact we are not directly interested by \emph{general} formulas for the
components of our reductions. It happens most often the algorithm obtained from
the vector field is much faster than which is obtained from the various
\emph{general} known formulas. Please examine the formulas~(\ref{76423}). You
see the homotopy operator is null for the target and critical cells, while the
terrible general formula for this operator, see~\cite[Section~6]{real2}, does
not give any hint for such a property. In fact the vector field is the
\emph{optimal} algorithm for this homotopy operator: running the different
\(V\)-paths from some source cell cannot be avoided, and the vector field \(V\)
``knows'' that no path at all is to be considered if the starting cell is
target or critical.

The same for the Alexander-Whitney operator. Why this operator is null for a
target cell? The vector field knows this point. It is not hard to prove the
Alexander-Whitney operator is a sort of identity map for a critical cell, but
again this is obvious from the \(f_p\)-formula~(\ref{76423}).

It is a general property of this technique of vector fields: it is not a
machine producing by itself ``closed'' formulas for the studied operators. On
the contrary, most often the vector field is the shortest way to obtain the
right result for every particular case.

Nevertheless, we have proved above the reduced complex is canonically
isomorphic to the usual tensor product, and this is necessary for the
continuation of the story about the twisted products. Also it is well known the
naturality of the Eilenberg-MacLane section implies this operator is unique,
and it is obvious our Eilenberg-Zilber vector field is \emph{natural} with
respect to its arguments. So that the Eilenberg-MacLane formula is certainly a
consequence of the Eilenberg-Zilber vector field.

\subsection{Obtaining the classical Eilenberg-Zilber formulas.}

When explicit formulas for the Eilenberg-Zilber reduction are required, the
following formulas are used. The components of the Eilenberg-Zilber reduction
are traditionnally called \(AW\) (Alexander-Whitney), \(EML\)
(Eilenberg-MacLane) and \(SHI\) (Shih):
\begin{eqnarray}
 AW(x_p \times y_p) &=& \sum_{i=0}^p \partial_{p-i+1} \cdots \partial_p x_p \otimes
 \partial_0 \cdots \partial_{p-i-1} y_p
 \\
 EML(x_p \otimes y_q) &=& \sum_{(\eta_1,\eta_2)\in Sh(p,q)}
 (-1)^{\varepsilon(\eta_1, \eta_2)} (\eta_2 x_p, \eta_1 y_q)
 \\
 SHI(x_p \times y_p) &=& \sum_{0 \leq r \leq p-1,0 \leq s \leq
p-r-1,(\eta_1,\eta_2) \in Sh(s+1,r)} (-1)^{p-r-s+\varepsilon(\eta_1,\eta_2)}
 \\
 && \makebox[11cm][r]{ \(
 (\uparrow^{p-r-s}(\eta_2) \eta_{p-r-s-1}
 \partial_{p-r+1} \cdots
 \partial_{p} x_p, \uparrow^{p-r-s}(\eta_1) \partial_{p-r-s} \cdots \partial_{p-r-1}
 y_p)\)}\nonumber
\end{eqnarray}

The \(AW\) formula defines the standard coproduct of a simplicial chain
complex. Two symmetric formulas are possible, and any appropriate linear
combination of both is also possible. The \(EML\) formula is unique. Two
symmetric Shih formulas are possible for \emph{any} \(AW\) formula. We have
given here which seems to be the most traditional choices. It is clear four
symmetric possible choices are possible for the Eilenberg-Zilber vector field,
our choice producing the above formulas.

In the \(EML\) formula above, The set \(Sh(p,q)\) is made of all the
\((p,q)\)-\emph{shuffles} of \([0..(p+q-1)]\), that is, all the partitions of
these \(p+q\) integers in two increasing sequences of length \(p\) and \(q\).
Every shuffle \((\eta_1,\eta_2)\) can be understood as a permutation, producing
a signature \(\varepsilon(\eta_1,\eta_2)\). For example, the shuffle \(((0\,
1\,5), (2\,3\,4))\) in the case \(p = q = 3\) produces the term \(-(\eta_4
\eta_3 \eta_2 x_3, \eta_5 \eta_1 \eta_0 y_3)\), for the permutation
\((0\,1\,5\,2\,\,3\,4)\) is negative.

The third formula is called \(SHI\) in~\cite{real2}, a reference to the
recursive formula given at~\cite[Page~25]{shih}. In fact this formula is
already at~\cite[Formula(~2.13)]{elmc2}, but because of the already used
\(EML\), because also of the beautiful work of Shih Weishu along these lines,
we continue to call it the \(SHI\) formula.

This section is devoted to a careful analysis of the Eilenberg-Zilber vector
field, leading to new proofs of all these formulas. Nothing more than a
combinatorial game with the s-paths of prisms, that is, with the degeneracy
operators.

\subsubsection{Naturality of the Eilenberg-Zilber reduction.}

We consider here four simplicial sets \(X\), \(Y\), \(X'\) and \(Y'\) and two
simplicial morphisms \(\phi: X \rightarrow X'\) and \(\psi: Y \rightarrow Y'\).
These morphisms induce a simplicial morphism \(\phi \times \psi: X \times Y
\rightarrow X' \times Y'\). Also the products \(X \times Y\) and \(X' \times
Y'\) carry their respective Eilenberg-Zilber vector fields \(V\) and \(V'\).

\begin{thr} ---
With these data, the morphisms \(\phi\) and \(\psi\) induce a natural morphism
between both Eilenberg-Zilber reductions:
\begin{eqnarray}
 \phi \times \psi: [\rho = (f,g,h): C_\ast(X \times Y) \rrdc C_*(X) \otimes
 C_\ast(Y)] \longrightarrow
 \\
 && \makebox[2cm][r]{\([\rho' = (f',g',h'): C_\ast(X' \times Y')
 \rrdc C_*(X') \otimes C_\ast(Y')]\)}\nonumber
\end{eqnarray}\end{thr}

The chain complexes are normalized. At this time of the process, we do not have
much information for the small chain complexes: we know the underlying graded
modules are (isomorphic to) those of \(C_*(X) \otimes C_\ast(Y)\) and \(C_*(X')
\otimes C_\ast(Y')\), but we do not yet know their differentials.

The goal is the following, as it is common in a simplicial environment: once
this naturality result is known, it is often enough to prove some desired
result in the particular case of some appropriate \emph{model}, maybe a prism
\(\Delta^{p,q}\), and then to use some obvious simplicial morphism to transfer
this result to an arbitrary product and obtain the general result.

\prof We just have to proof, taking account of Theorem~\ref{03005}, the
morphism \(\phi \times \psi\) is compatible with the respective vector fields
\(V\) and \(V'\).

\section{The twisted Eilenberg-Zilber W-reduction.}

\subsection{Twisted products.}

If \(X\) and \(Y\) are two simplicial sets, the usual definition of the product
\(X \times Y\) was recalled in Section~\ref{49366}. A \emph{fibration} is a
sort of ``twisted'' product, a notion having a major role in algebraic
topology. We follow here the terminology and the notations of~\cite[\S18]{may}.

\begin{dfn} ---
\emph{A \emph{simplicial morphism} \(f: X \rightarrow Y\) between two
simplicial sets \(X\) and \(Y\) is a collection \((f_p: X_p \rightarrow Y_p)_{p
\in \bN}\) of maps between simplex sets, compatible with the \(\uD\)-operators:
for every \(\uD\)-morphism \(\alpha: \up \rightarrow \uq\), the relation
\(\alpha^\ast f_q = f_p \alpha^\ast\) is satisfied.\qed}
\end{dfn}

Think this diagram must be commutative:
\begin{ctikzpicture} [scale = 1.5, baseline = 0.8cm]
 \foreach \i in {0,1} \foreach \j in {0,1}
   {\coordinate (\i\j) at (\i,\j) ;}
 \foreach \i/\ii in {0/X, 1/Y} \foreach \j/\jj in {0/q, 1/p}
   {\node (n\i\j) at (\i\j) {\(\ii_{\jj}\)} ;}
 \begin{scope} [font = \footnotesize, ->]
 \draw (n01) -- node [above] {\(f_p\)} (n11) ;
 \draw (n00) -- node [above] {\(f_q\)} (n10) ;
 \draw (n00) -- node [left] {\(\alpha^\ast\)} (n01) ;
 \draw (n10) -- node [right] {\(\alpha^\ast\)} (n11) ;
 \end{scope}
\end{ctikzpicture}

See Section~\ref{84511}. The face and degeneracy operators generate all the
\(\uD\)-morphisms, so that it is enough this compatibility condition is
satisfied for face and degeneracy operators.

\begin{dfn} ---
\emph{A \emph{simplicial group} \(G\) is a simplicial set provided with two
simplicial morphisms, a group law \(\mu_G: G \times G \rightarrow G\) and an
inversion map \(\iota_G: G \rightarrow G\); every homogeneous \(p\)-dimensional
pair \((\mu_{G,p}, \iota_{G,p})\) must satisfy the usual group axioms.\qed}
\end{dfn}

Every homogeneous simplex set \(G_p\) is endowed with a group structure, and
the collection of groups \((G_p)_p\) is compatible with face and degeneracy
operators. We most often simply write in a multiplicative way \(\mu_G(g,g') =
g.g'\) or \(gg'\) and \(\iota_G(g) = g^{-1}\), or sometimes in an additive form
\(\mu_G(g,g') = g + g'\) and \(\iota_G(g) = -g\) if the group law is abelian.

\begin{dfn} ---
\emph{A simplicial action \(\gamma\) of the simplicial group \(G\) on a
simplicial set \(X\) is a simplicial morphism \(\gamma: G \times X \rightarrow
X\) satisfying the usual axioms of a group action.\qed}
\end{dfn}

For every \(p \in \bN\), a group action \(\gamma_p: G_p \times X_p \rightarrow
X_p\) is defined and all these actions are compatible with face and degeneracy
operators. Most often we will not denote the action \(\gamma\) by a letter,
using simply the product notation: \(\gamma(g,x)\) will be simply denoted by
\(g.x\) or \(gx\).

\begin{dfn} ---
\emph{Let \(F\) and \(B\) be two simplicial sets, the \emph{fiber space} and
the \emph{base space} of the fibration to be defined. Let \(G\) be a simplicial
group, the \emph{structural group} and let \(G \times F \rightarrow F\) be some
action of the structural group on the fiber space. A \emph{twisting function}
is a collection of maps \((\tau_p: B_p \rightarrow G_{p-1})_{p > 0}\)
satisfying the conditions:
\begin{equation}
\begin{array}{rcll}
 \partial_0(\tau b) &=& \tau(\partial_0 b)^{-1} . \tau(\partial_1 b)
 \\
 \partial_i \tau (b) &=& \tau(\partial_{i+1} b) & i > 0
 \\
 \eta_i \tau(b) &=& \tau(\eta_{i+1} b) & i \geq 0
 \\
 e_p &=& \tau(\eta_0 b) & b \in B_p
\end{array}
\end{equation}
if \(e_p\) is the neutral element of the group \(G_p\). }
\qed\end{dfn}

\begin{dfn} ---
\emph{If \(G, F, B\) and \(\tau\) are as in the previous definition, the
\emph{twisted product} \(F \times_\tau B\) is the simplicial set defined as
follows. Every homogeneous simplex set \((F \times_\tau B)_p\) is \emph{the
same as} for the non-twisted product: \((F \times_\tau B)_p = (F \times_\tau
B)_p = F_p \times B_p\). Only the 0-face operator is different:
\begin{equation}
\partial_0(f,b) = (\tau(b).\partial_0 f, \partial_0 b).
\end{equation}
The other face and degeneracy operators therefore are simply \(\partial_i(f,b)
= (\partial_i f, \partial_i b)\) for \(i > 0\) and \(\eta_i(f,b) = (\eta_i f,
\eta_i b)\) for arbitrary relevant \(i\).\qed}
\end{dfn}

In other words, the twisting function \(\tau\) is used only to \emph{perturb}
the 0-face operator in the \emph{vertical} direction if, as usual, we think of
the base space as the horizontal component of the (twisted) product and the
fiber space as the vertical component.

Unfortunately, our reference book~\cite{may} does not give any easily
understandable example. The reader could profitably examine the
notes~\cite{srgr01}, in particular Sections~7 and 12.

\subsection{The twisted Eilenberg-Zilber vector field.}

The notion of twisted product, due to Daniel Kan~\cite[Section 6]{kan},
recalled in the previous section, is really a \emph{miracle}. The point is that
only the 0-face operator is modified by the twisting process, while the
Eilenberg-Zilber vector field described in Section~\ref{08795} invokes
``vectors'' \((\sigma, \sigma')\) where the face-index~\(i\) in the regular
face relation \(\partial_i \sigma' = \sigma\) always satisfies \(i > 0\): the
0-face is never concerned in this vector field.

\begin{thr} ---
Let \(F \times_\tau B\) be a twisted product as defined in the previous
section. Then the Eilenberg-Zilber vector field defined in Theorem~\ref{46504}
for the non-twisted product \(F \times B\) can be used as well for the twisted
product \(F \times_\tau B\). It is admissible and therefore defines a
homological reduction:
\begin{equation}
 TEZ: C_\ast(F \times_\tau B) \rrdc C_\ast(F) \otimes_t C_\ast(B).
\end{equation}
The small chain complex of this reduction has the same underlying graded module
as \(C_\ast(F) \otimes C_\ast(B)\); only the differential is modified, which is
indicated by the index~\(t\) of the twisted product symbol `\(\otimes\)': it is
a \emph{twisted} tensor product.
\end{thr}


\prof We reuse the terminology and the notations of Section~\ref{25115}. Every
\(r\)-simplex \(\rho = (\phi, \beta) = (\eta_{i_{s-1}} \cdots\, \eta_{i_{0}}
\phi', \eta_{j_{t-1}} \cdots\, \eta_{j_{0}} \beta')\) of the non-twisted
cartesian product \mbox{\(F \times B\)} is as well an \(r\)-simplex of the
twisted product \(F \times_\tau B\). We decide the degeneracy configuration
\(((i_{s-1}, \ldots, i_0), (j_{t-1}, \ldots, j_{0}))\) determines the nature of
the simplex \(\rho\) in the vector field \(V_{F \times_\tau B}\) to be defined
on \(F \times_\tau B\), source, target or critical, exactly like in
Definition~\ref{52397}: the respective status of \(\rho\) with respect to the
vector fields \(V_{F \times B}\) (resp. \(V_{F \times_\tau B}\)) of \(F \times
B\) (resp. \(F \times_\tau B\)) are \emph{the same}.

In particular, if a pair \((\rho_1, \rho_2) = ((\phi_1, \beta_1), (\phi_2,
\beta_2))\) is a \emph{vector} of the Eilenberg-Zilber vector field \(V_{F
\times B}\) of \(F \times B\), we decide it is as well a vector ot the twisted
Eilenberg-Zilber vector field \(V_{F \times_\tau B}\) we are defining. The
incidence relation \(\rho_1 = \partial_i \rho_2\) is satisfied in \(F \times
B\) for a unique face index \(i\) satisfying \(0 < i < r\) if \(r\) is the
dimension of \(\rho_2\); so that the \emph{same} incidence relation is also
satisfied in the twisted cartesian product \(F \times_\tau B\), for only the
0-face is modified in the twisted product.

However the 0-face operator plays a role when studying whether a vector field
is \emph{admissible}. The admissibility proof of the vector field \(V_{X \times
Y}\) is based on the Lyapunov function \(L\) defined and used in
Proposition~\ref{79914} and on the partition of all the simplices of \(X \times
Y\) in layers indexed by the bidimension \((p,q)\), see Definition~\ref{79914}.
This proof amounts to defining a new Lyapunov function \(L(\rho_1) = (p, q,
L(\chi))\) if the bidimension of the simplex \(\rho_1\) is \((p, q)\),
certainly the same as the bidimension of \(\rho_2\), and if \(\chi\) is the
degeneracy configuration of \(\rho_1\); we use the lexicographic order to
compare the images of this function.

We claim the same Lyapunov function can be used for the twisted product \(F
\times_\tau B\). Let us assume \((\rho_1, \rho_2) \in V_{F \times_\tau B}\)
with \(\rho_1 = \partial_i \rho_2\) and \(\rho_3 = \partial_j \tau\) with \(j
\neq i\). If \(j > 0\), then the \(j\)-face \(\partial_j \tau\) is the same in
both products, twisted or non-twisted, so that \(L'(\rho_3)\) has the same
value whatever the product you consider, twisted or non-twisted.

Let \(\rho_2 = (\phi_2, \beta_2) = (\eta_{i_{s-1}} \cdots\, \eta_{i_{0}}
\phi_2', \eta_{j_{t-1}} \cdots\, \eta_{j_{0}} \beta_2')\) be the canonical
expression of~\(\rho_2 = (\phi_2, \beta_2)\) using the non-degenerate simplices
\(\phi'_2 \in F_{p}\) and \(\beta'_2 \in B_{q}\) if the bidimension of
\(\rho_2\) is \((p,q)\). The first formula below gives the \(0\)-face in \(F
\times B\) while the second one gives the \(0\)-face in \(F \times_\tau B\):
\begin{equation}
\begin{array}{rcrcl}
 \partial_0 \rho_2 &=& (\partial_0 \eta_{i_{s-1}} \cdots\, \eta_{i_{0}}
 \phi'_2,
 \partial_0 \eta_{j_{t-1}} \cdots\, \eta_{j_{0}} \beta'_2)
 \\
 \partial_{0,\tau} \rho_2 &=& (\tau(\beta_2)
 \partial_0 \eta_{i_{s-1}} \cdots\, \eta_{i_{0}} \phi'_2, \partial_0
 \eta_{j_{t-1}} \cdots\, \eta_{j_{0}} \beta'_2)
\end{array}
\end{equation}
It has been observed in Proposition~\ref{75744} the 0-face of an
\emph{interior} simplex of \(\Delta^{p,q}\) is always \emph{exterior}. In other
words, the 0-face of a non-degenerate simplex always has a stictly smaller
bidimension. The group operator \(\tau(\beta_2)\) is a simplicial
\emph{isomorphism}, so that the \emph{geometrical} dimensions of the
\emph{first} components of \(\partial_0 \rho_2\) and \(\partial_{0,\tau}
\rho_2\), that is, the dimensions of the corresponding non-degenerate simplices
given by the Eilenberg-Zilber lemma~\ref{34608}, are the same. This implies the
bidimension of \(\partial_0 \rho_2\) and \(\partial_{0,\tau} \rho_2\) are the
same, therefore stictly smaller than the bidimension \((p,q)\) of \(\rho_2.\)

This implies the Lyapunov function \(L'\) can be used for both products and
their respective vector fields \(V_{F \times B}\) and \(V_{F \times_\tau B}\).
This vector field \(V_{F \times_\tau B}\) therefore is \emph{admissible}.\qed

\section[\texorpdfstring{The [\boldmath$BG - \textrm{Bar}(G)$]
  W-reduction.}
  {The [BG - Bar(G)] W-reduction}]
{The [\boldmath$\Omega X - \textrm{Cobar}(X)$] W-reduction.} \label{00896}

\section[\texorpdfstring{The [\boldmath$\Omega X - \textrm{Cobar}(X)$]
  W-reduction.}
  {The [Omega X - Cobar(X)] W-reduction}]
{The [\boldmath$\Omega X - \textrm{Cobar}(X)$] W-reduction.} \label{40626}

\section{The Adams model of a loop space.}

\section{The Bousfield-Kan W-reduction (??).}


\end{document}